\documentclass[preprints,article,accept,moreauthors,pdftex]{Definitions/mdpi} 

\firstpage{1} 
\pubvolume{xx}
\issuenum{1}
\articlenumber{5}
\pubyear{2019}
\copyrightyear{2019}
\history{} 

\pdfoutput=1

\usepackage{verbatim} 
\usepackage{epsfig}
\usepackage{amssymb,array}
\usepackage{picinpar}
\usepackage{subcaption} 
\usetikzlibrary{shapes,arrows}
\usetikzlibrary{positioning,automata}
\usepackage[mathscr]{euscript}
\usepackage{physics}

\allowdisplaybreaks 
\definecolor{red}{rgb}{1,0.00,0.00}
\definecolor{blue}{rgb}{0,0.08,0.55}
\definecolor{green}{rgb}{0,0.55,0.08}
\definecolor{grey}{rgb}{0.87,0.87,0.87}

\newcommand{\ronertwo}[1]{\textcolor{black}{#1}} 
\newcommand{\ronerthree}[1]{\textcolor{black}{#1}} 
\newcommand{\supun}[1]{\textcolor{black}{#1}} 

\newcommand{\Aut}{\mathop{\rm Aut}\nolimits}

\newcommand{\inv}{^{-1}}
\newcommand{\om}{_\omega}
\newcommand{\omz}{_{\omega_0}}
\newcommand{\som}{_{T\omega}}

\newcommand{\ol}[1]{\overline{#1}}

\newcommand{\Sch}[1]{\mathcal{X}_{#1}^{\mathop{\rm Sch}\nolimits}}
\newcommand{\ltwo}{\ell^2(\Gamma)}
\newcommand{\ltwoG}{\ell^2(G)}
\newcommand{\spec}{\mathop{\rm Sp}\nolimits}
\newcommand{\jspec}{\mathop{\rm JSp}\nolimits}
\newcommand{\mass}[1]{\delta_{#1}}
\newcommand{\gen}[1]{\left< #1 \right>}

\newcommand{\A}{\mathcal A}
\newcommand{\GG}{\mathcal G}
\newcommand{\B}{\mathcal B}

\newcommand{\Gr}{\mathcal G}
\newcommand{\Hh}{\mathcal H}
\newcommand{\J}{\mathcal{J}}
\newcommand{\M}{\mathcal{M}}
\newcommand{\Sy}{\mathcal S}
\newcommand{\OO}{\mathcal{O}}
\newcommand{\LL}{\mathcal L}
\newcommand{\NN}{\mathcal N}
\newcommand{\X}{\mathcal X}
\newcommand{\K}{\mathcal K}

\newcommand{\R}{\mathbb{R}}
\newcommand{\CC}{\mathbb{C}}
\newcommand{\Z}{\mathbb Z}
\newcommand{\N}{\mathbb N}

\newcommand{\mkdschreier}[1]{\M^*_{Sch}(#1)}
\newcommand{\wt}[1]{\widetilde{#1}}

\DeclareMathOperator{\stab}{St}

\newcommand{\arr}{\longrightarrow}

\makeatletter
\newcommand*{\mtset}{%
	\nfss@text{%
		\sbox0{0}%
		\sbox2{/}%
		\sbox4{%
			\raise\dimexpr((\ht0-\dp0)-(\ht2-\dp2))/2\relax\copy2 %
		}%
		\ooalign{%
			\hfill\copy4 \hfill\cr
			\hfill0\hfill\cr
		}%
		\vphantom{0\copy4 }
	}%
}

\Title{Integrable  and  chaotic  systems associated  with  fractal  groups\footnote[0]{This article is accepted by \emph{Entropy} and is pending publication. To comply with \emph{Entropy} guidelines, the references are numbered in order of appearance in the text.}}

\Author{Rostislav Grigorchuk and Supun Samarakoon* \orcidS{}}

\AuthorNames{Rostislav Grigorchuk and Supun Samarakoon}

\address[1]{%
Department of Mathematics, Texas A\&M University, College Station, TX 77843, USA; grigorch@math.tamu.edu~(R.G.); sts@tamu.edu~(S.S.)}

\corres{Correspondence: sts@tamu.edu (S.S.)}

\abstract{Fractal  groups  (also called self-similar  groups)  is the class  of  groups  discovered by  the  first  author  in  the 80-s  of  the  last  century with  the  purpose  to  solve some famous  problems  in  mathematics,   including  the  question  raising to von Neumann  about  non-elementary  amenability (in the  association  with  studies  around   the  Banach-Tarski  Paradox)  and  John  Milnor's   question  on  the  existence  of  groups  of intermediate  growth  between  polynomial  and  exponential.
Fractal  groups arise  in  various  fields  of  mathematics,  including the theory of  random  walks, holomorphic  dynamics,  automata theory, operator  algebras,  etc.  They  have  relations  to  the  theory  of  chaos,  quasi-crystals, fractals,  and  random Schr\"odinger operators.
One  of important  developments is  the  relation  of  them  to  the  multi-dimensional  dynamics,  theory  of joint  spectrum  of pencil  of  operators,  and  spectral theory of  Laplace  operator on graphs.    The  paper  gives a quick  access   to  these  topics,  provide calculation  and  analysis  of multi-dimensional  rational  maps arising via  the Schur  complement in some  important  examples,   including  the  first  group of  intermediate  growth and its overgroup,  contains  discussion  of the  dichotomy  ``integrable-chaotic'' in the considered  model, and  suggests a possible probabilistic approach  to  the  study  of discussed problems.}


\begin{document}



\section{Introduction}

\begin{figure}[H]
\centering
\includegraphics[width=0.9\textwidth]{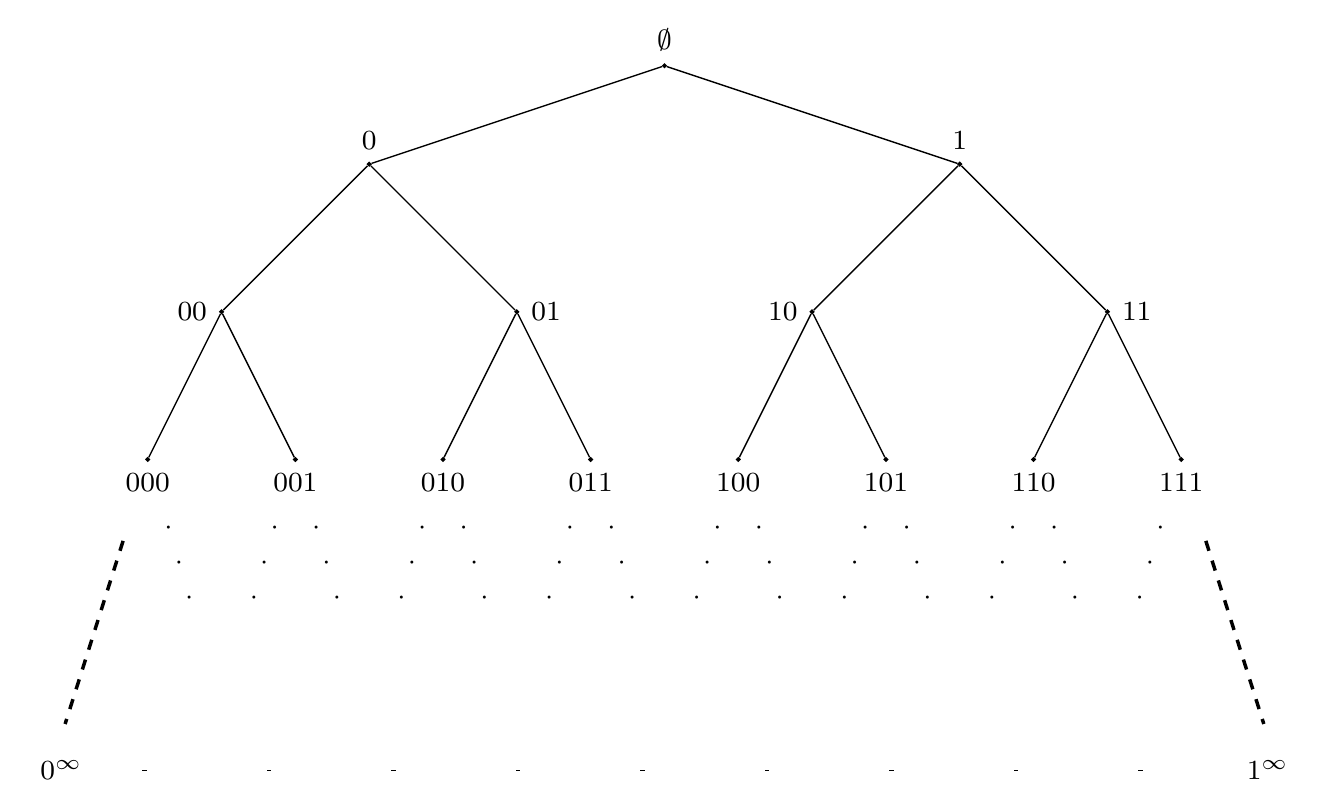}
\caption{Binary rooted tree, $T_2$, where the vertices are identified with $\{0,1\}^*$.}	\label{fig:T_2}
\end{figure}

\begin{figure}[H]
	\centering
	\begin{subfigure}[b]{0.3\textwidth}
		\includegraphics[width=\textwidth]{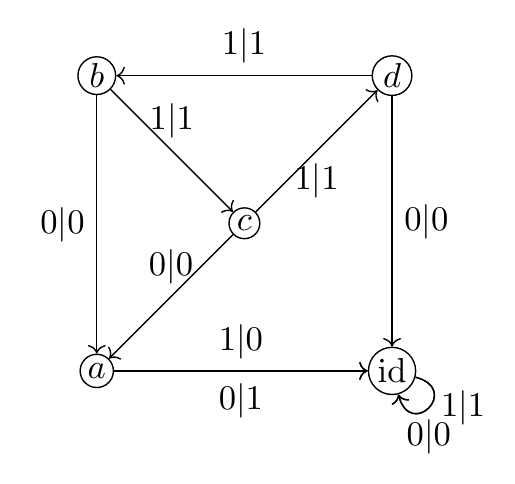}
		\caption{}
		\label{subfig:gri}
	\end{subfigure}
	\qquad
	\begin{subfigure}[b]{0.3\textwidth}
		\includegraphics[width=\textwidth]{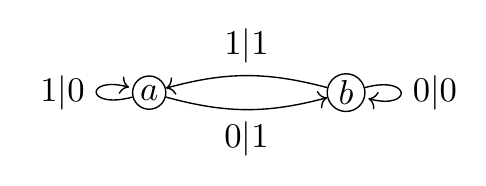}
		\caption{}
		\label{subfig:lamplighter}
	\end{subfigure}
	\\ \vspace*{8pt}
	\begin{subfigure}[b]{0.35\textwidth}
		\includegraphics[width=\textwidth]{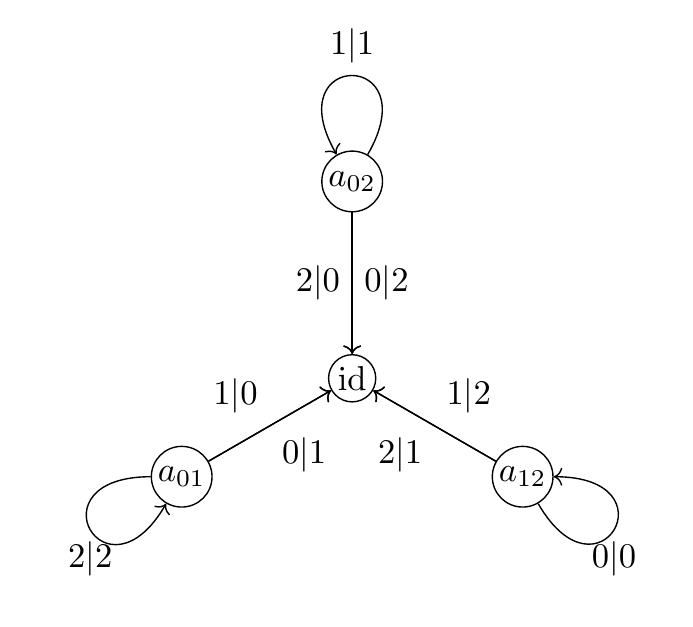}
		\caption{}
		\label{subfig:hanoi3}
	\end{subfigure}
	\qquad
	\begin{subfigure}[b]{0.35\textwidth}
		\includegraphics[width=\textwidth]{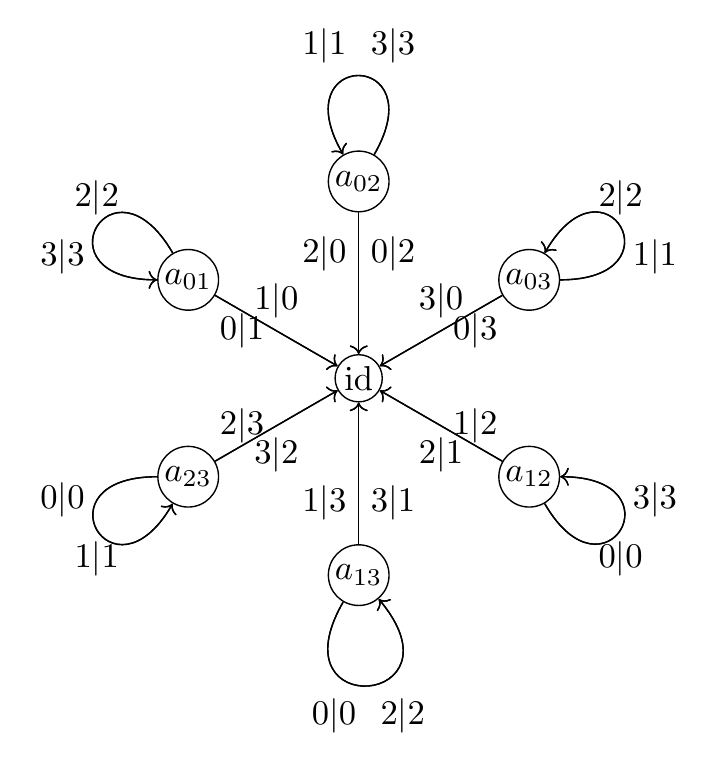}
		\caption{}
		\label{subfig:hanoi4}
	\end{subfigure}
	\\ \vspace{8pt}
	\begin{subfigure}[b]{0.24\textwidth}
		\includegraphics[width=\textwidth]{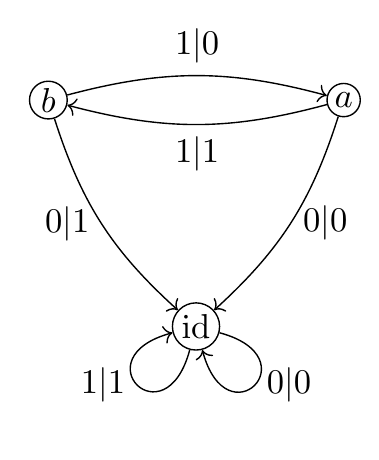}
		\caption{}
		\label{subfig:basilica}
	\end{subfigure}
	\qquad
	\begin{subfigure}[b]{0.29\textwidth}
		\includegraphics[width=\textwidth]{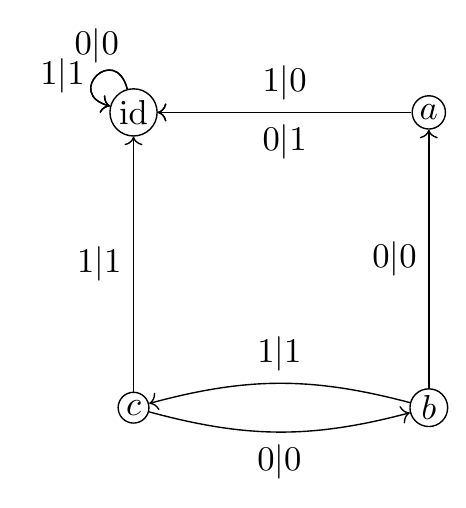}
		\caption{}
		\label{subfig:img_plus_i}
	\end{subfigure}
	\\ \vspace*{8pt}
	\begin{subfigure}[b]{0.3\textwidth}
		\includegraphics[width=\textwidth]{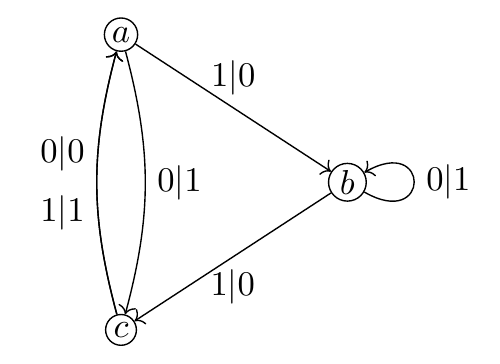}
		\caption{}
		\label{subfig:free3}
	\end{subfigure}
	\qquad
	\begin{subfigure}[b]{0.3\textwidth}
		\includegraphics[width=\textwidth]{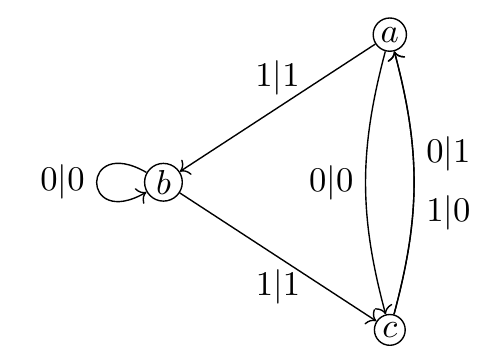}
		\caption{}
		\label{subfig:bellaterra}
	\end{subfigure}
	\caption{Examples of finite automata generating (\textbf{a}) Grigorchuk group $\Gr$, (\textbf{b}) Lamplighter group, (\textbf{c}) Hanoi tower group $\Hh^{(3)}$, (\textbf{d}) Hanoi tower group $\Hh^{(4)}$, (\textbf{e}) Basilica group, (\textbf{f}) $IMG(z^2+i)$, (\textbf{g}) free group $F_3$ of rank three, and (\textbf{h}) $\Z_2 \ast \Z_2 \ast \Z_2$.}\label{fig:automata}
\end{figure}

Fractal groups are groups acting on self-similar objects in \ronerthree{a} self-similar way. \ronerthree{The term ``fractal group'' was used for the first time in \cite{BG00}} and then appeared in \cite{BGN03}. Although there is no rigorous definition of \ronerthree{a} fractal group (like there is no rigorous definition of \ronerthree{a} fractal set), there is a definition of \ronerthree{a} self-similar group (see Definition \ref{def:self-similar gp}). Self-similar groups act by automorphisms on regular rooted trees (like a binary rooted tree shown in {Figure} \ref{fig:T_2}). 
Such trees are among the most natural and often used self-similar objects. The properties  of self-similar groups  and their structure resemble the self-similarity properties of the trees and their boundaries. The nicest examples come from finite Mealy type automata, like automata presented by {Figure} \ref{fig:automata}. 

Moreover, there are several ways to associate  geometric objects of fractal type with a self-similar group. This includes limits of {Schreier graphs} 
\cite{BG00,BG00Hecke,GNS00,Vor12,GKN12}, limit spaces and limit solenoids of Nekrashevych \cite{Nek09,Nek01,Nek03,BHN06}, quasi-crystals \cite{GLN17bSchreier}, Julia sets \cite{BG00,Nek05}, etc.

Self-similar groups were used to solve several outstanding problems in different areas of mathematics. They provide an elegant contribution to the general Burnside problem \cite{Gri80}, to the J.~Milnor problem on growth \cite{Gri83,Gri84}, to the von Neumann - Day problem on non-elementary amenability \cite{Gri84,Gri98}, to the Atiyah problem in $L^2$-Betti numbers \cite{GLSZ00}, etc. Self-similar groups have applications in many areas of mathematics such as dynamical systems, operator algebras, random walks, spectral theory of groups and graphs, geometry and topology, computer science, and many more (see the surveys \cite{BGN03,Gri05,Gri11,GLN17bSchreier,GN07,Gri14Mil,GNS00,GNS15} and the monograph \cite{Nek05}).

Multi-dimensional rational maps appear in \ronerthree{the} study of spectral properties of graphs and unitary representations of groups (including representations of Koopman type). The spectral theory of such objects  is closely related to the theory of joint spectrum of \ronerthree{a} pencil of operators in a Hilbert (or more generally in a Banach) space and is implicitly considered in \cite{BG00} and explicitly outlined in \cite{Yan09}.

There is some mystery regarding how multi-dimensional rational maps appear in the context of self-similar groups. There are basic examples like the first group $\Gr$ of intermediate growth from \cite{Gri80,Gri83}, the groups called ``Lamplighter'', ``Hanoi'', ``Basilica'', Spinal Groups, GGS-groups, etc. The indicated classes of groups produce a large family of such maps, part of which is presented by examples \eqref{eq:F map 1stGri}--\eqref{eq:G map 1stGri} and \eqref{eq:5dim F map 1stGri}--\eqref{eq:map img con}.

These maps are \ronerthree{\emph{very special}} and \ronerthree{\emph{quite degenerate}} as claimed by N. Sibony and M. Lyubich, respectively. Nevertheless, they are interesting and useful, as, on the one hand, they are \ronerthree{responsible} for the associated spectral problems, on the other hand, they give a lot of \ronerthree{material} for people working in dynamics, being quite different from the maps that were considered before.

\begin{figure}[H]
\centering
\begin{subfigure}[b]{0.45\textwidth}
	\includegraphics[width=\textwidth]{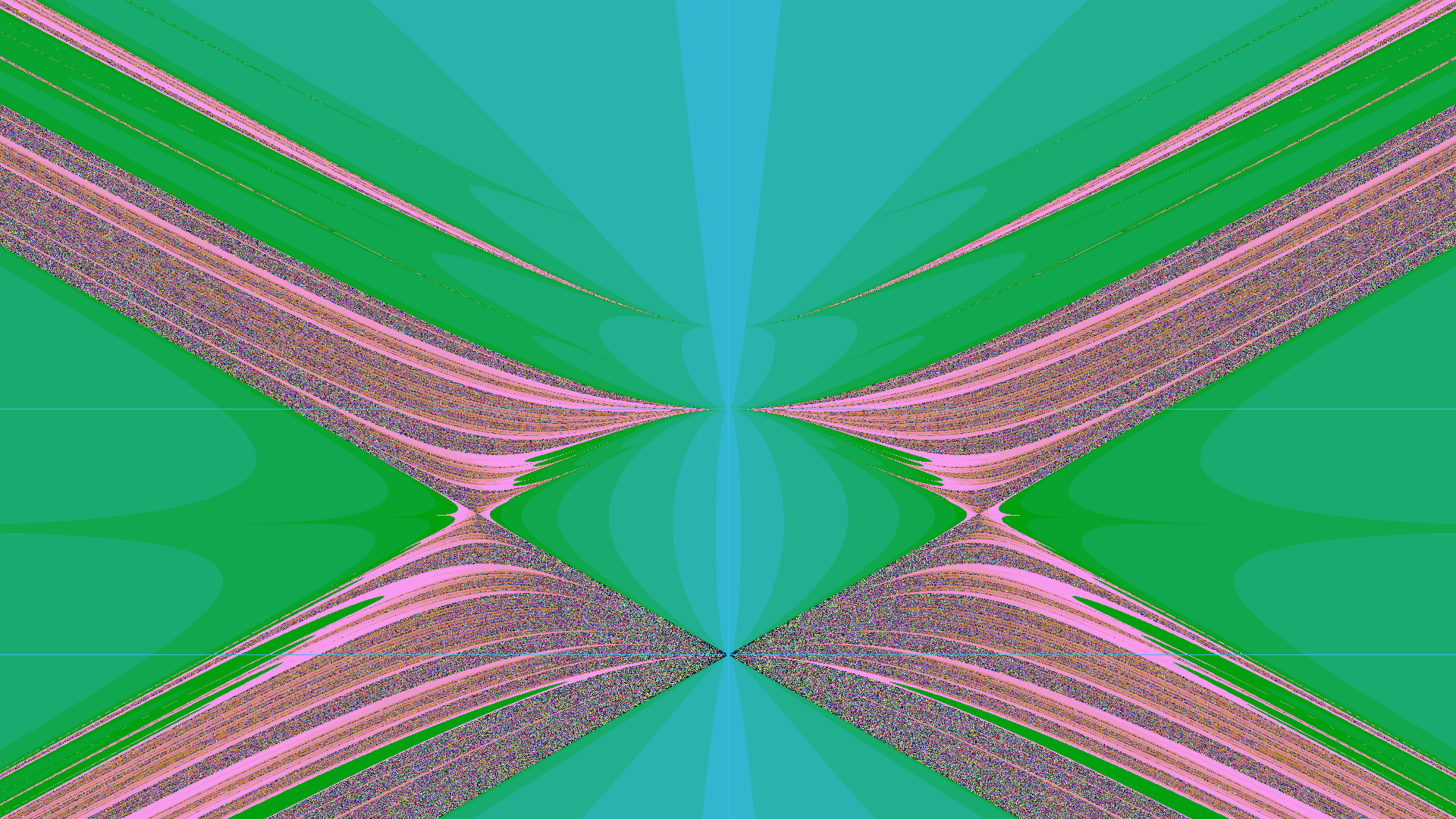}
	\caption{}\label{subfi:dynamical map 012}
\end{subfigure}
\qquad
\begin{subfigure}[b]{0.45\textwidth}
	\includegraphics[width=\textwidth]{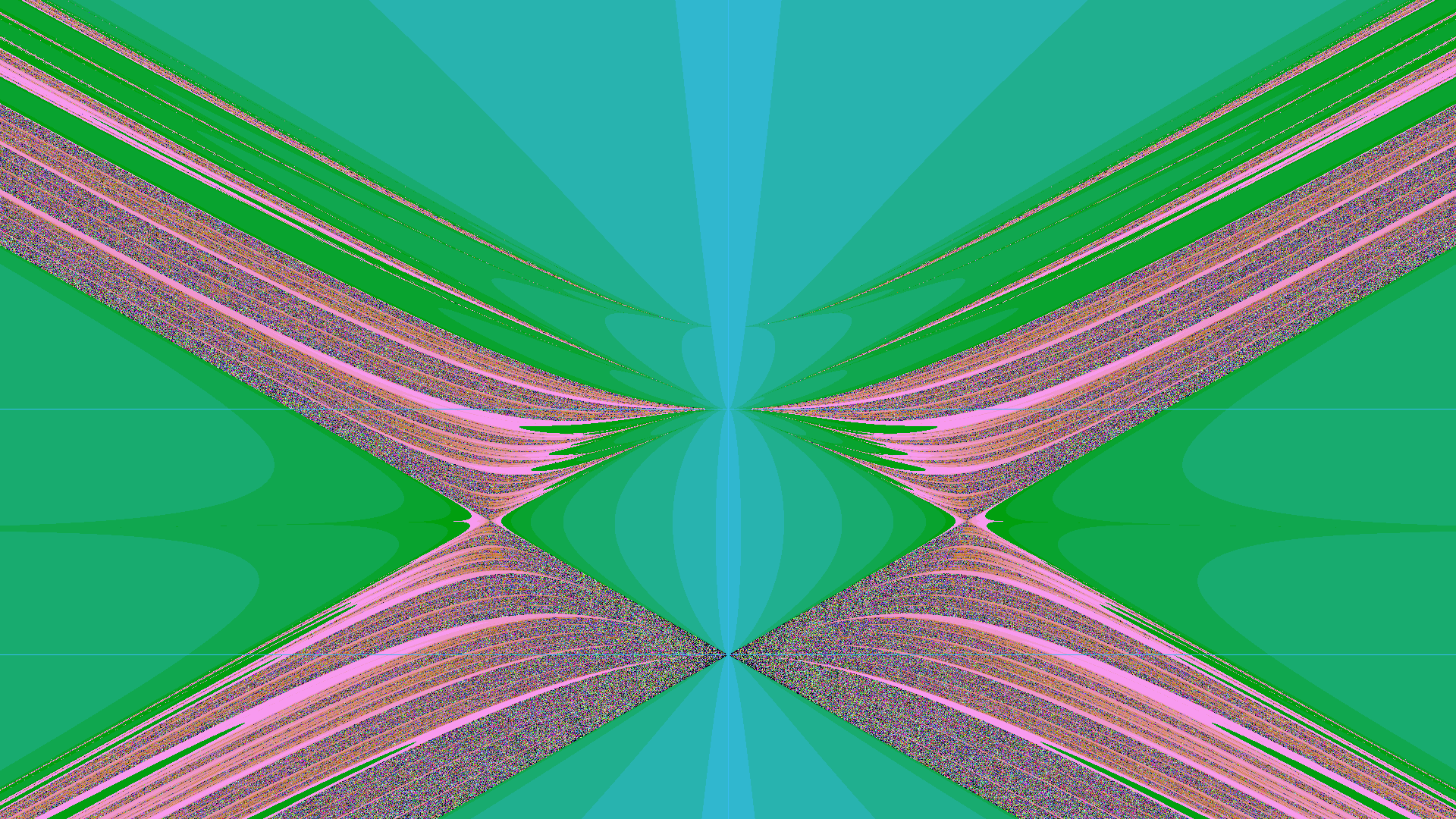}
	\caption{}
\end{subfigure}
\\ \vspace*{8pt}
\begin{subfigure}[b]{0.45\textwidth}
	\includegraphics[width=\textwidth]{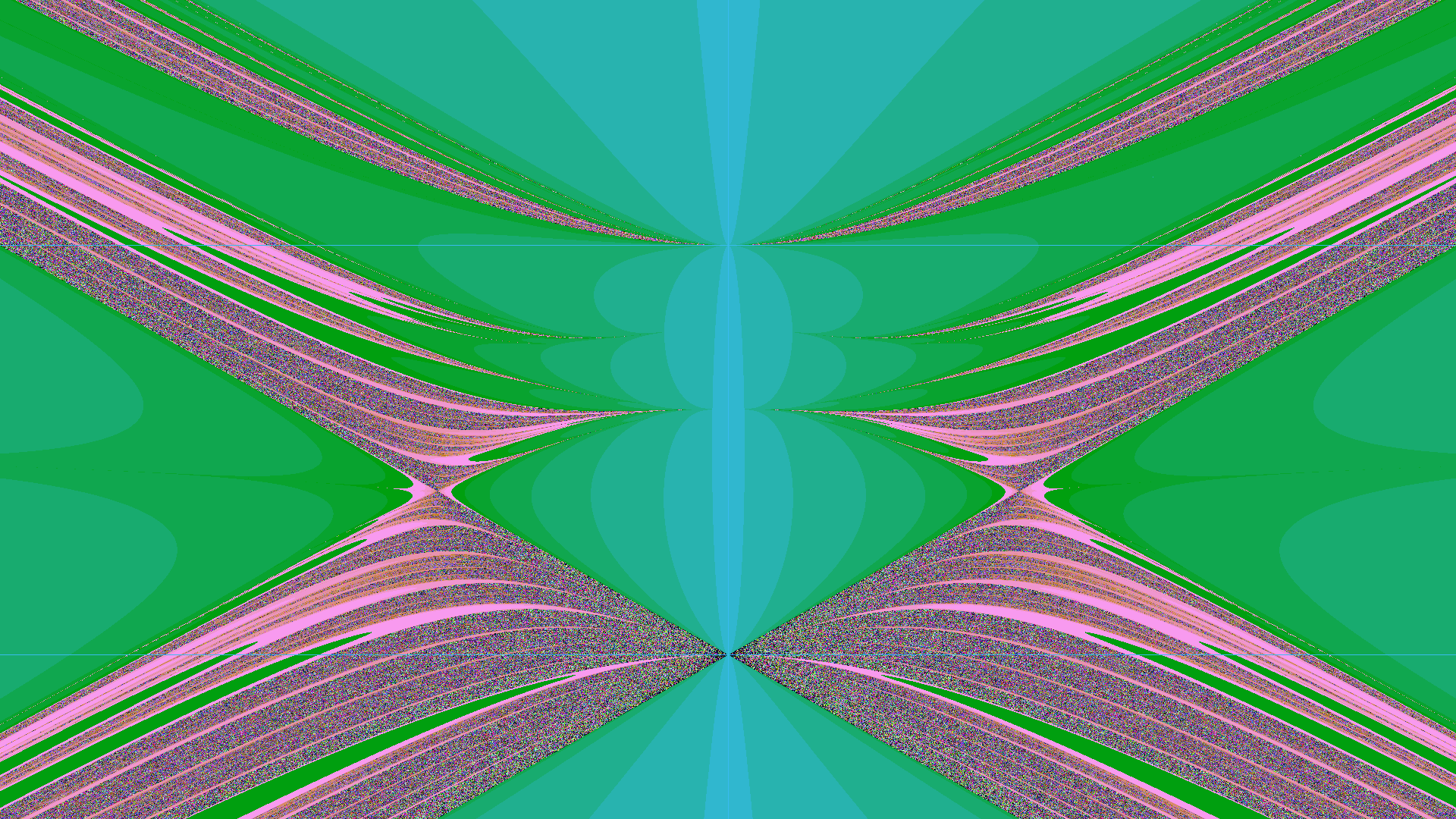}
	\caption{}
\end{subfigure}
\qquad
\begin{subfigure}[b]{0.45\textwidth}
	\includegraphics[width=\textwidth]{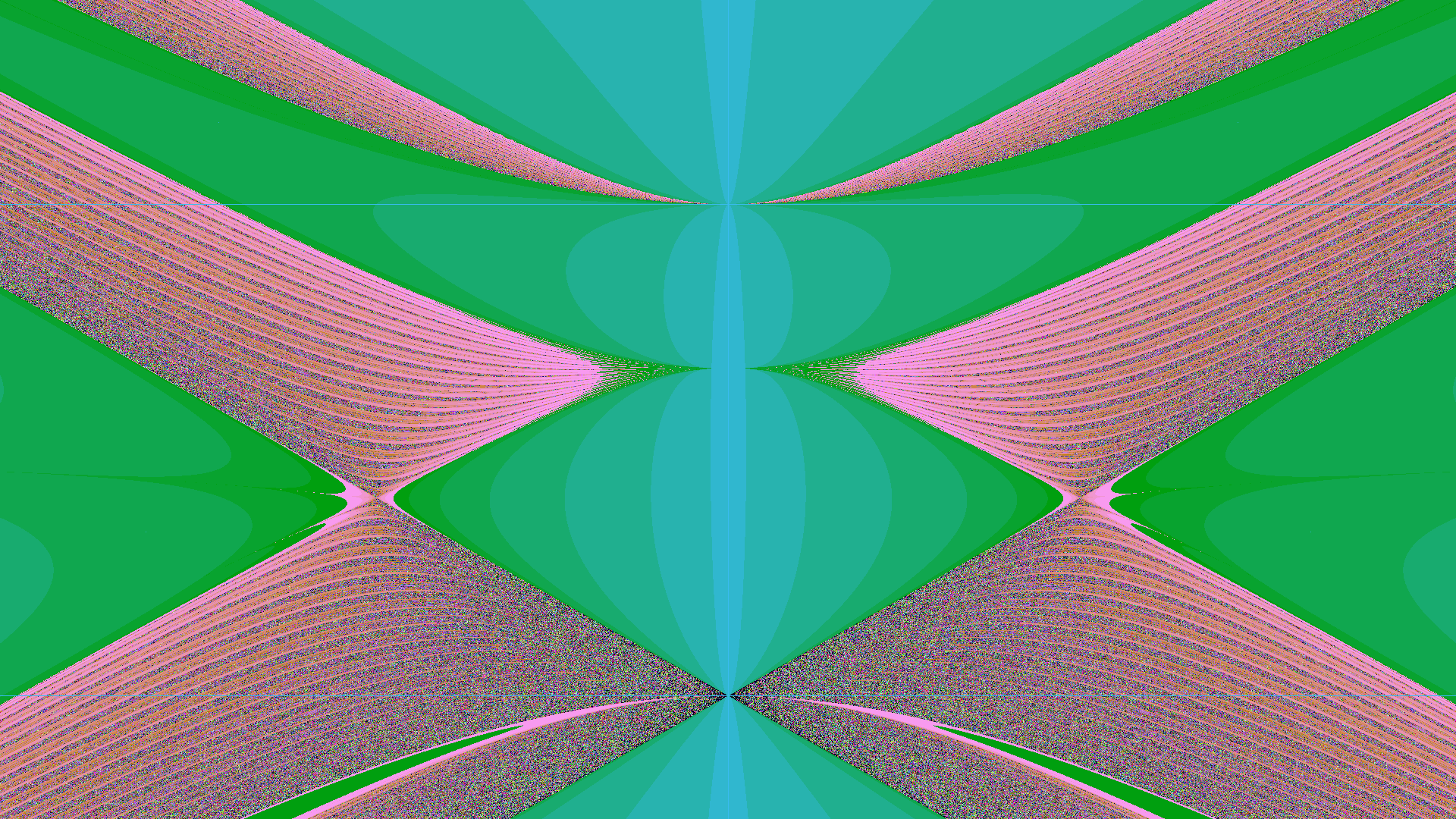}
	\caption{}
\end{subfigure}
\vspace*{8pt}
\caption{Dynamical pictures of $F_{\omega_{n-1}} \circ \hdots \circ F_{\omega_0}$ for (\textbf{a}) $\omega = (012)^\infty$ and $(y,z,u) = (1,2,3)$, (\textbf{b}) $\omega = (01)^\infty$ and $(y,z,u) = (1,2,3)$, (\textbf{c}) a random $\omega$ and $(y,z,u) = (1,2,3)$, and (\textbf{d}) a random $\omega$ and $(y,z,u) = (1,3,3)$.}\label{fig:random_maps}
\end{figure}

Some of them demonstrate features of integrability, which means that they \ronerthree{semiconjugate} to lower-dimensional maps, while the others do not seem to have integrability features and their dynamics (at least on an experimental level) demonstrate the chaotic behavior presented, for instance, by {Figure}~\ref{fig:random_maps}. 

The phenomenon of integrability, discovered in the basic examples including the groups $\Gr$, Lamplighter, and Hanoi, is thoroughly investigated by M-B. Dang, M. Lyubich and the first author in \cite{DGL20}. For examples of intermediate complexity, the criterion found in \cite{DGL20} based on the fractionality of the dynamical degree shows non-integrability; for instance, this is the case for  the Basilica map \eqref{eq:map Basilica}. More complicated cases of maps, like the Basilica map, or higher-dimensional $\Gr$-maps given by \eqref{eq:5dim F map 1stGri} and \eqref{eq:5dim G map 1stGri} still wait for their resolution. An interesting phenomenon discovered in \cite{GLN18} is the relation of self-similar groups with quasi-crystals and random Schr\"odinger operators.

\ronertwo{This article surveys and explores the use of self-similar groups} in the dynamics of multi-dimensional rational maps and provides a panorama of ideas, methods, and applications of fractal groups.

\section{Basic Examples}
We begin with several examples of rational maps that arise from fractal groups. Only minimal information is given about each case. The main examples are related to the group $\Gr$ given by presentation \eqref{eq:1st gp presentation}, overgroup $\wt\Gr$ given by matrix recursions \eqref{eq:matrix map overgp}, and generalized groups $\Gr\om$ given by \eqref{eq:gen_gri_gp}. We begin with dimension 2 and then consider the higher-dimensional case given by \eqref{eq:5dim F map 1stGri}, \eqref{eq:5dim G map 1stGri}. The justification is given in Section \ref{sec:schur_calculation}.

\subsection{Grigorchuk Group}

Consider two maps
\begin{equation}\label{eq:F map 1stGri}
F\colon\left(\begin{array}{c}
x\\
y
\end{array}\right)
\mapsto
\left(\begin{array}{c}
\frac{2x^2}{4-y^2}\\[3mm]
y+\frac{x^2y}{4-y^2}
\end{array}\right),
\end{equation}
\begin{equation}\label{eq:G map 1stGri}
G\colon\left(\begin{array}{c}
x\\
y
\end{array}\right)
\mapsto
\left(\begin{array}{c}
\frac{2(4-y^2)}{x^2}\\[3mm]
-y-\frac{y(4-y^2)}{x^2}
\end{array}\right).
\end{equation}
They come from the group of intermediate growth, between polynomial and exponential, \cite{Gri80, Gri84}

\begin{equation}\label{eq:1st gp presentation}
\mathcal G= \gen{ a,b,c,d | 1=a^2=b^2=c^2=d^2=bcd=\sigma^k((ad)^4)=\sigma^k((adacac)^4), k=0,1,2,\dots },
\end{equation}  
where
\[\sigma: a\rightarrow aca, b \rightarrow d, c\rightarrow b, d\rightarrow c\]
is a substitution. The maps $F$ and $G$ are related by $H \circ F = G$ and $H \circ G = F$, where $H$ is the involutive map (i.e., $H \circ H = id$)
\begin{equation}\label{eq:H map}
H\colon\left(\begin{array}{c}
x\\
y
\end{array}\right)
\mapsto
\left(\begin{array}{c}
\frac{4}{x}\\[3mm]
-\frac{2y}{x}
\end{array}\right).
\end{equation}
The point of interest is the dynamics of $F,G$ acting on $\R^2, \CC^2$ or their projective counterparts and the dynamics of the subshift $(\Omega_{\sigma},T)$ generated by the substitution $\sigma$ (which is briefly discussed in the Section \ref{sec:dynamical system}).

\begin{figure}[H]
\centering
\begin{subfigure}[b]{0.4\textwidth}
	\includegraphics[width=\textwidth]{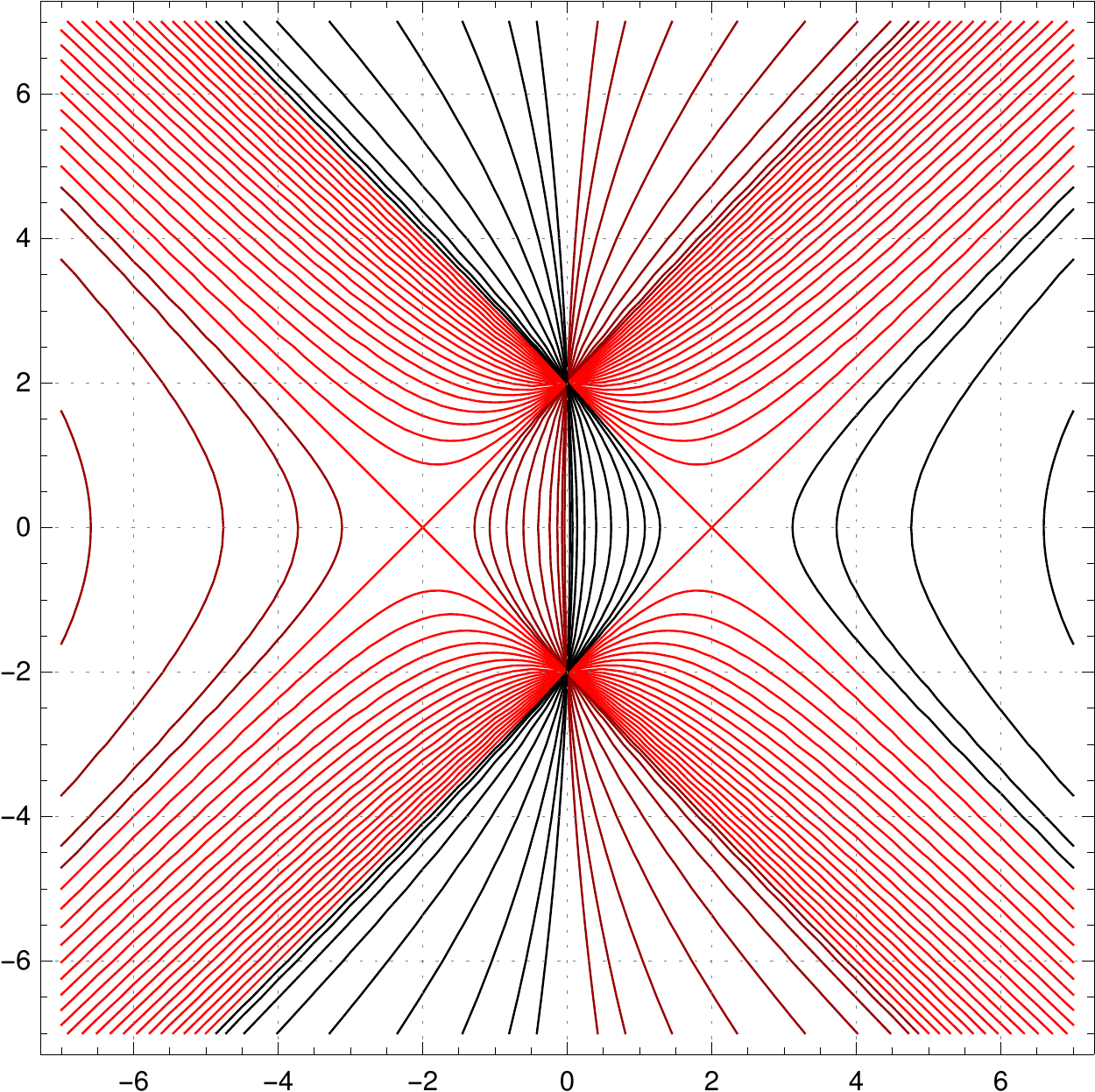}
	\caption{}
	\label{subfig:hor_hyperbola_all_var_G}
\end{subfigure}
\qquad
\begin{subfigure}[b]{0.4\textwidth}
	\includegraphics[width=\textwidth]{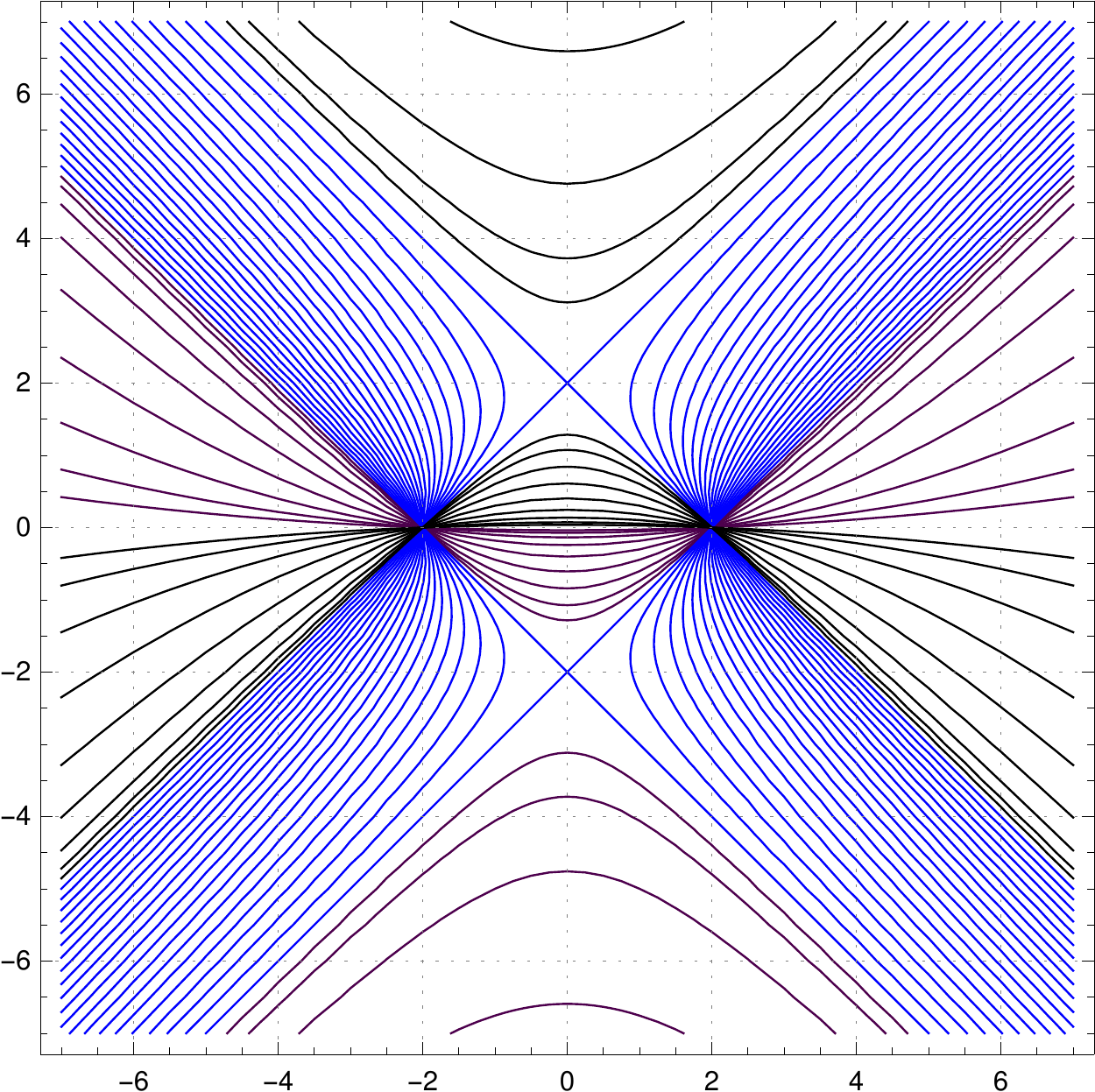}
	\caption{}
	\label{subfig:ver_hyperbola_all_var_G}
\end{subfigure}
\caption{Foliation of $\R^2$ by (\textbf{a}) horizontal hyperbolas $\mathcal{F}_\theta$ where, maroon, red and black corresponds to $\theta <-1, \theta \in [-1,1]$ and $\theta >1$, respectively, and (\textbf{b}) vertical hyperbola $\mathcal{H}_\eta$ where, purple, blue and black corresponds to $\eta <-1, \eta \in [-1,1]$ and $\eta >1$, respectively.}\label{fig:hyperbola_all_var_G}
\end{figure}

The map $F$ demonstrates features of an integrable map as it has two \emph{almost transversal} families of \emph{horizontal} hyperbolas $\mathcal  F_{\theta}=\{(x,y) \colon 4+x^2-y^2-4\theta x = 0 \}$ and \emph{vertical} hyperbolas ${\mathcal H_{\eta} =\{(x,y)\colon4-x^2+y^2-4\eta y=0\}}$, shown in Figure \ref{fig:hyperbola_all_var_G}. The first family $\{\mathcal{F}_\theta \}$ is invariant as a family and $F^{-1}(\mathcal{F}_\theta) = \mathcal{F}_{\theta_1} \sqcup  \mathcal{F}_{\theta_2}$,  where  $\theta_1,\theta_2$  are  preimages  of  $\theta$ under the Chebyshev  map $\alpha \colon z \mapsto 2z^2-1$ (also known as the Ulam - von Neumann map),  and the family $\{ \mathcal{H}_\eta \}$ consists of invariant curves.

The map $\pi = \varphi \times \psi$, where
\begin{align*}
	\psi(x,y) = & \frac{4+x^2-y^2}{4x}, \\ 
	\varphi(x,y) = & \frac{4-x^2+y^2}{4y},
\end{align*}
semiconjugates $F$ to the map $id \times \alpha$ and as the dynamics of $\alpha$ is well understood, some additional arguments lead to;

\begin{Theorem}[Equidistribution Theorem \cite{DGL20}]
	Let  $\Gamma$ and  $S$ be  two irreducible  algebraic curves in $\mathbb{C}^2$ in coordinates $(\varphi,\psi)$ such  that  $\Gamma$  is  not a vertical hyperbola  while  $S$  is not a horizontal hyperbola.  Then
	\[\frac{1}{2^n}[(F^n)^{\ast}\Gamma \cap S]\xrightarrow{n\to \infty} (deg\Gamma)\cdot( degS)\cdot \omega_S,\]
	where  $\omega_S$ is  the  restriction of  the  $1$-form  $\omega=\frac{d\psi}{\pi\sqrt{1-\psi^2}}$   to  $S$, $[\cdot]$ is  the counting  measure, and $F^n$ denotes the $n$-th iteration of $F$.
\end{Theorem}

\begin{figure}[H]
\centering
\begin{subfigure}[b]{0.29\textwidth}
	\includegraphics[width=\textwidth]{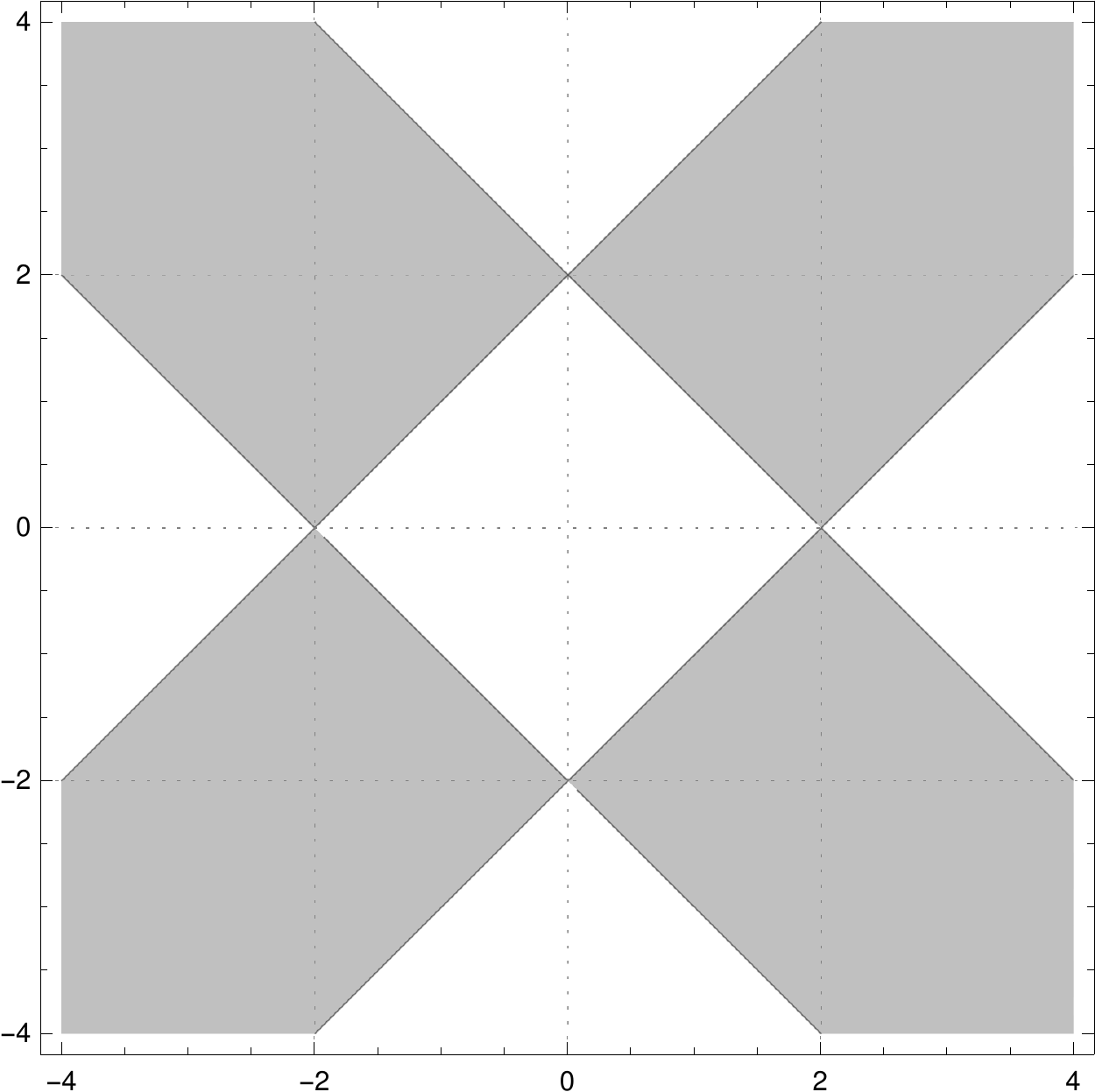}
	\caption{}
	\label{subfig:all_hyperbola_interval_var_G}
\end{subfigure}
\qquad
\begin{subfigure}[b]{0.29\textwidth}
	\includegraphics[width=\textwidth]{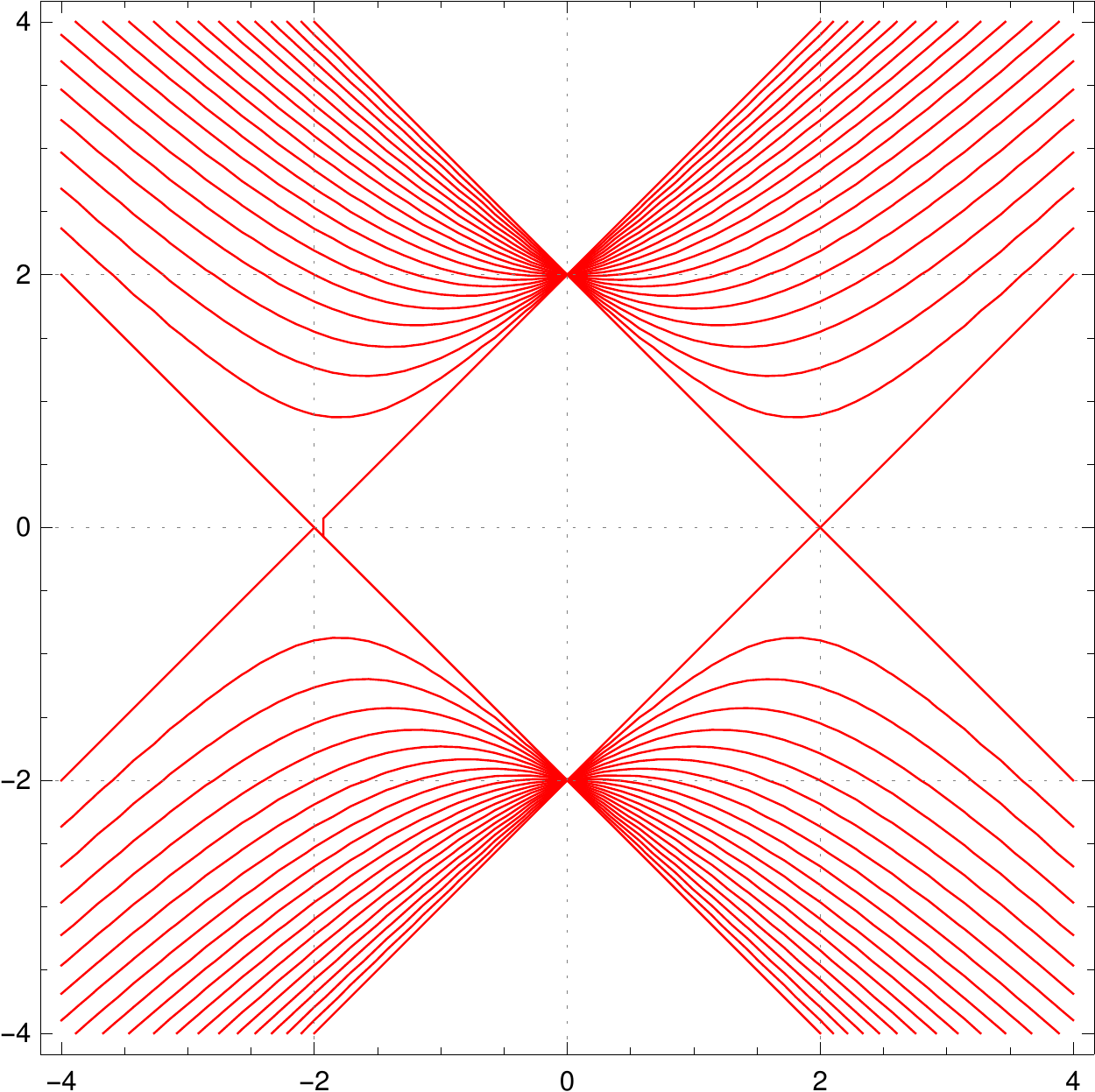}
	\caption{}
	\label{subfig:hor_hyperbola_interval_var_G}
\end{subfigure}
\qquad
\begin{subfigure}[b]{0.29\textwidth}
	\includegraphics[width=\textwidth]{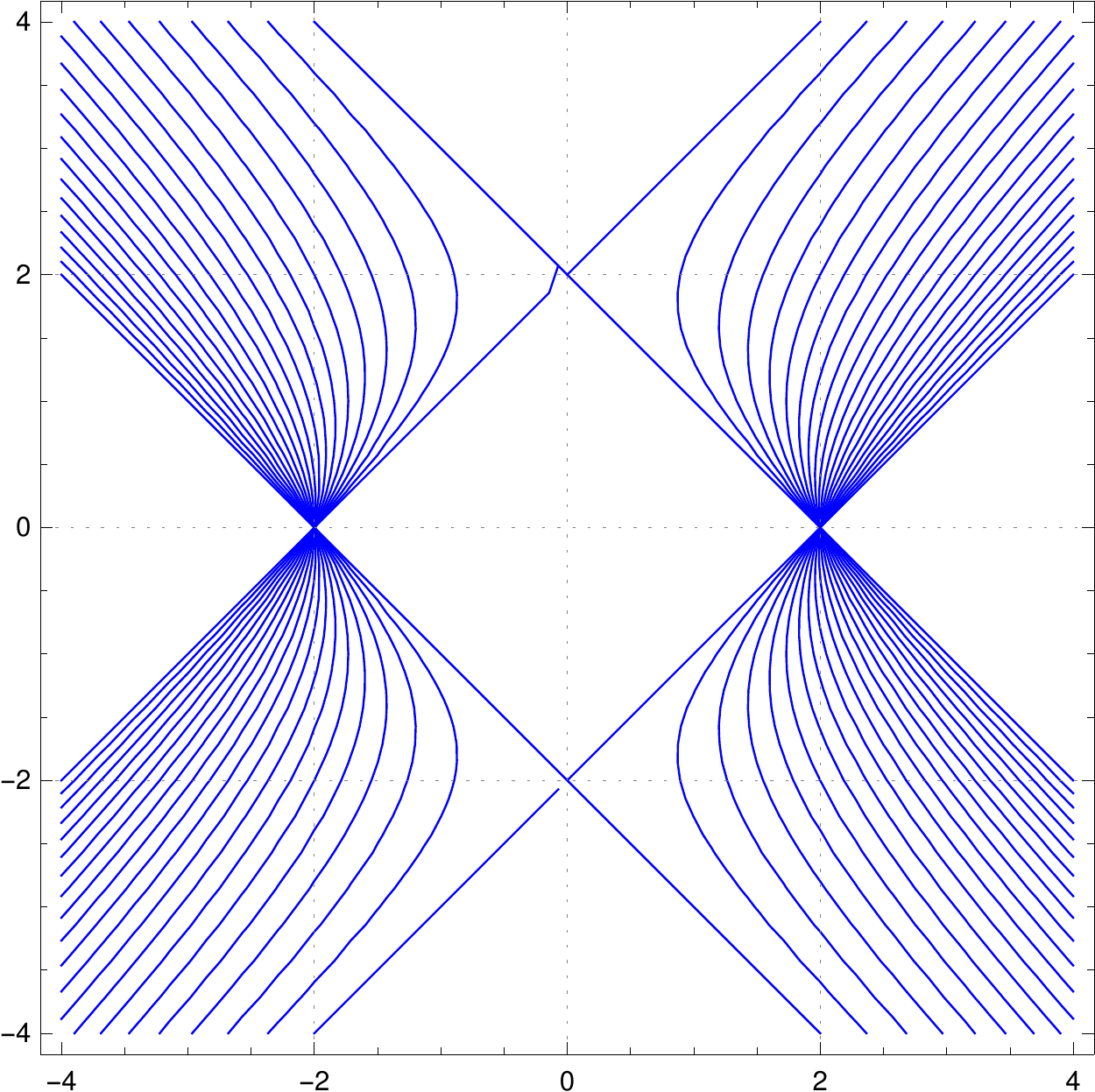}
	\caption{}
	\label{subfig:ver_hyperbola_interval_var_G}
\end{subfigure}
\caption{(\textbf{a}) The ``cross'' $\K$, (\textbf{b}) foliation of $\K$ by real slices of horizontal hyperbolas $\mathcal{F}_\theta$ ($\theta \in [-1,1]$), and (\textbf{c}) foliation of $\K$ by real slices of vertical hyperbolas $\mathcal{H}_\eta$ ($\eta \in [-1,1]$).}\label{fig:hyperbola_interval_var_G}
\end{figure}

The set $\K$ shown in Figure \ref{subfig:all_hyperbola_interval_var_G} (we will call this set the ``cross'') is of special interest for us as it represents \ronerthree{the} joint spectrum of several families of operators associated with \ronerthree{the} element ${m(x,y) = -xa +b+c+d - (y+1)1}$ of the group algebra $\R[\Gr]$ \cite{BG00,GN07,DG17}.  It can be foliated by the hyperbolas $\mathcal{F}_\theta, -1 \leq \theta \leq 1$ as shown in Figure \ref{subfig:hor_hyperbola_interval_var_G} (or by hyperbolas $\mathcal{H}_\eta, -1 \leq \eta \leq 1$ shown in Figure~\ref{subfig:ver_hyperbola_interval_var_G}). The $F$-preimages of the border line $x+y=2$ constitutes a dense family of curves for $\K$ (the same is true for $G$-preimages) and $\K$ is completely invariant set for $F$ or $G$ (i.e., $F\inv(\K) \subset \K$ and $F(\K) \subset \K$, so $F(\K) = \K$).
The map $F$ is comprehensively investigated in \cite{DGL20} (its close relative is studied in \cite{GY17} and \cite{GY20} from a different point of view) and serves as a \ronerthree{basis} for the integrability theory developed there. The map $G$ happens to be more complicated and its study is ongoing.

$F$ and $G$ are low-dimensional relatives of $\CC^5 \to \CC^5$ maps
\begin{equation}\label{eq:5dim F map 1stGri}
\wt{F} \colon 
\left(
\begin{array}{c}
x \\
y \\
z \\
u \\
v \\
\end{array}
\right)
\mapsto 
\left(
\begin{array}{c}
\frac{x^2(y+z)}{(v+u+y+z)(v+u-y-z)} \\
u \\
y \\
z \\
v -  \frac{x^2(v+u)}{(v+u+y+z)(v+u-y-z)}  \\
\end{array}
\right),
\end{equation}
\begin{equation}\label{eq:5dim G map 1stGri}
\wt{G} \colon 
\left(
\begin{array}{c}
x \\
y \\
z \\
u \\
v \\
\end{array}
\right)
\mapsto 
\left(
\begin{array}{c}
y+z \\
-  \frac{x^2(2vyz - u(v^2-u^2+y^2+z^2))}{(v+u+y+z)(v-u+y-z)(v+u-y-z)(v-u-y+z)} \\
- \frac{x^2(2vuz - y(v^2+u^2-y^2+z^2))}{(v+u+y+z)(v-u+y-z)(v+u-y-z)(v-u-y+z)} \\
- \frac{x^2(2vuy - z(v^2+u^2+y^2-z^2))}{(v+u+y+z)(v-u+y-z)(v+u-y-z)(v-u-y+z)} \\
v+u -  \frac{x^2(2uyz - v(-v^2+u^2+y^2+z^2))}{(v+u+y+z)(v-u+y-z)(v+u-y-z)(v-u-y+z)}  \\
\end{array}
\right),
\end{equation}
that come from the 5-parametric pencil $xa + yb+zc+ud +v1 \in \R[\Gr]$ \cite{GN07}. The way how they were computed is explained in Section \ref{sec:schur_calculation}.

It is known that there is a subset $\Sigma \subset \R^5$ \ronerthree{which is} both $\wt{F}$ and $\wt{G}$-invariant, \ronerthree{and} the sections of which, by the lines parallel to the direction when $y=z=u$ are unions of two intervals (or an interval), while in all other directions it is a Cantor set of the Lebesgue measure zero. This follows from the results of D.~Lenz, T Nagnibeda and the first author \cite{GLN18} and is based on the use of the substitutional dynamical system $(\Omega_\sigma,T)$ determined by substitution $\sigma$ (see Section \ref{sec:dynamical system} and also \cite{GLN17aCombinatorics,GLN17bSchreier,GLN18}).

The integrability criterion found in \cite{DGL20} clarifies the roots of integrability of the maps \eqref{map lamplighter}, \eqref{map Hanoi} presented in the next two examples.

\subsection{Lamplighter Group}

The Lamplighter  map
\begin{equation}\label{map lamplighter}
\left(\begin{array}{c}
x\\
y
\end{array}\right)
\mapsto
\left(\begin{array}{c}
\dfrac{x^2-y^2-2}{y-x}\\[3mm]
\dfrac{2}{y-x}
\end{array}\right)
\end{equation}
comes from the Lamplighter group $\LL = \Z_2 \wr \Z = \left(\oplus_{\Z} \Z_2\right) \rtimes \Z$ (where, $\wr$ and $\rtimes$ denote the wreath product and the semidirect product, respectively) realized as a group generated by the automaton in {Figure} 
\ref{subfig:lamplighter} (generation of groups by automata is explained in the next section). It has a family of invariant lines $l_c \equiv x+y=c$, and its restriction to $l_c$ is the M\"obius map represented by the matrix
\[	\left(\begin{array}{cc}	c & -\dfrac{c^2}{2} - 1\\	1 & - \dfrac{c}{2}	\end{array}\right)	\in SL(2,\R).	\]

The Lamplighter map was used in \cite{GZ01} and \cite{GS19} to describe unusual spectral properties of the Lamplighter group, and to answer the Atiyah question on $L^2$-Betti numbers \cite{GLSZ00}.

\subsection{Hanoi Group}
\begin{figure}[H]
	\centering
	\includegraphics[width=0.4\textwidth]{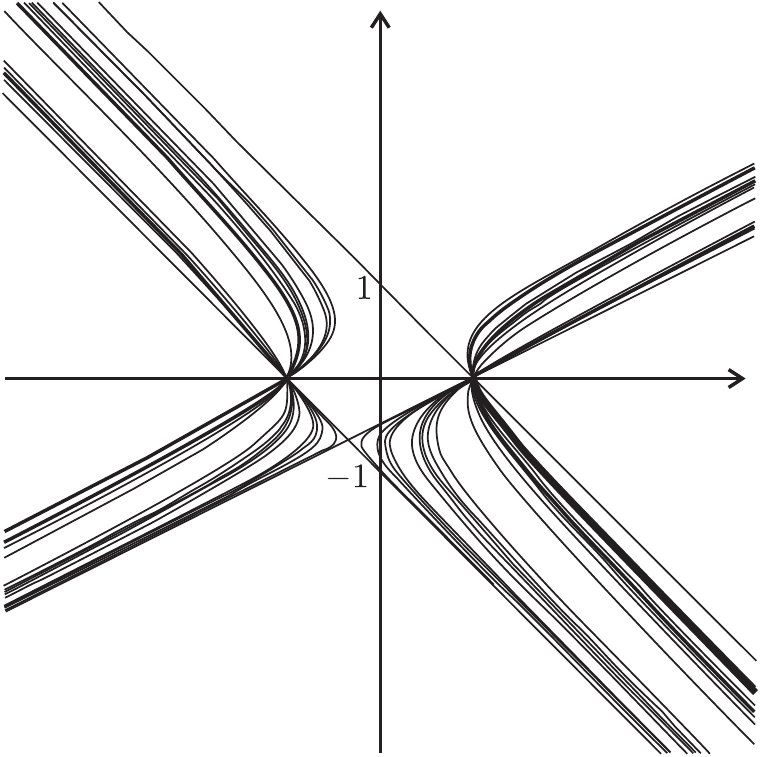}
	\caption{Joint spectrum of $\Hh^{(3)}$}\label{fig:hanoi spectrum}
\end{figure}

The Hanoi  map
\begin{equation}\label{map Hanoi}
\left(\begin{array}{c}
x\\
y
\end{array}\right)
\mapsto
\left(\begin{array}{c}
x-\dfrac{2(x^2-x-y^2)y^2}{(x-y-1)(x^2+y-y^2-1)}\\[3mm]
\dfrac{(x+y-1)y^2}{(x-y-1)(x^2+y-y^2-1)}
\end{array}\right)
\end{equation}
comes   from  the Hanoi  group $\mathcal H^{(3)}$, associated with  the Hanoi  Towers  Game on three  pegs \cite{GS06,GS07}. It is the  group generated by the automaton in {Figure} \ref{subfig:hanoi3}. 
As shown in \cite{GS07}, this map is semiconjugate by ${\psi \colon \R^2 \to \R}$, $\psi(x,y) = \dfrac{1}{y} (x^2 -1 -xy -2y^2)$ to the map ${\beta \colon x \mapsto x^2 -x -3}$ and has an invariant set $\Sigma$ (the joint spectrum) shown in Figure \ref{fig:hanoi spectrum}.

The set $\Sigma$ is the closure of the family of hyperbolas ${x^2 -1 -xy -2y^2 - \theta y = 0}$ when ${\displaystyle \theta \in \bigcup_{i=0}^\infty \beta^{-i}(\theta) \cup \bigcup_{i=0}^\infty \beta^{-i}(-2)}$ and the intersection of $\Sigma$ by vertical lines is a union of a countable set of isolated points and a Cantor set to which this set of points accumulates. Theorem \ref{thm:Hanoi_spectrum}, stated later, describes the intersection of $\Sigma$ by horizontal line $y=1$.

In \cite{GS06}, \emph{higher-dimensional} Hanoi groups $\Hh^{(n)}$, $n\geq 4$ are also introduced (the automaton generating $\Hh^{(4)}$ is presented by {Figure} \ref{subfig:hanoi4}), 
but unfortunately so far it was not possible to associate any maps with these groups.

\subsection{Basilica Group}
\begin{figure}[h]
	\centering
	\includegraphics[width=0.6\textwidth]{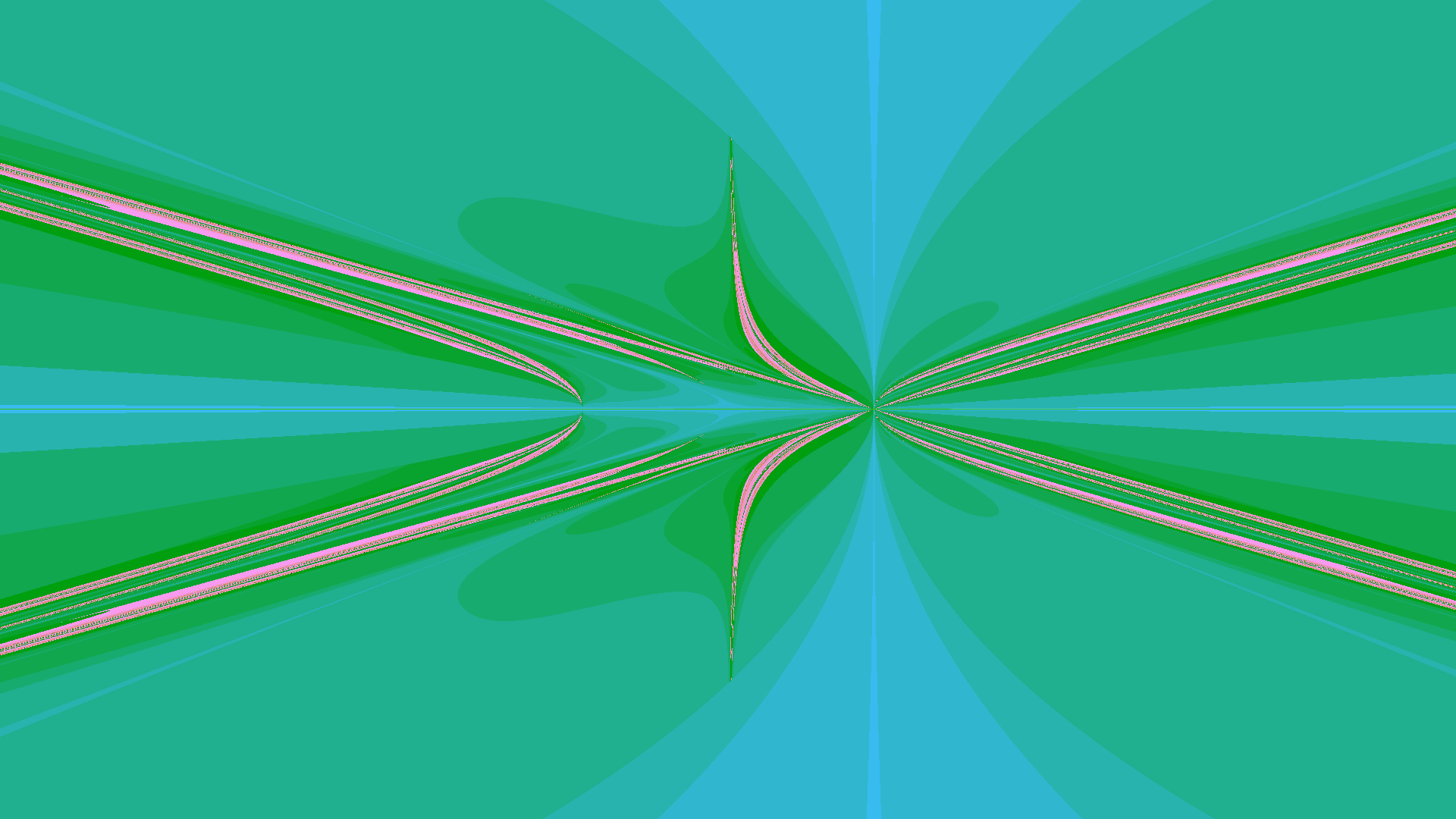}
	\caption{Dynamical picture of the Basilica map $B$.}\label{fig:basilica map}
\end{figure}

The Basilica  map
\begin{equation}\label{eq:map Basilica}
\left(\begin{array}{c}
x\\
y
\end{array}\right)
\mapsto
\left(\begin{array}{c}
-2+\dfrac{x(x-2)}{y^2}\\[3mm]
\dfrac{2-x}{y^2}
\end{array}\right)
\end{equation}
comes   from the  Basilica  group, introduced in \cite{GZ02Ontorsionfree,GZ02Spectral} as \ronerthree{the} group generated  by the automaton in {Figure}~\ref{subfig:basilica}. 
The map in \eqref{eq:map Basilica} is much more complicated than the other 2-dimensional maps described above. The dynamical picture \ronerthree{representing points of bounded orbits (a kind of a Julia set)} of the Basilica map is presented by Figure~\ref{fig:basilica map}.

The  Basilica group  \supun{can also} be defined as the Iterated  Monodromy  group  $IMG(z^2-1)$ of  the polynomial  $z^2-1$ \cite{BGN03,Nek05}. It  is the  first  example  of  amenable  but not  subexponentially  amenable  group, as shown in \cite{BV05}. For spectral properties of this group see \cite{BGJRT20,GZ02Spectral}.

\subsection{Iterated Monodromy Group of $z^2 + i$}

The group $IMG(z^2+i)$ is the group generated by the automaton in {Figure} \ref{subfig:img_plus_i}. 
The map 
\begin{equation}\label{eq:map img}
\left(\begin{array}{c}
y\\z\\
\lambda
\end{array}\right)\mapsto \left(\begin{array}{c}
\dfrac
zy\\[3mm]
\dfrac{1}{-y^2+z^2-2z\lambda+\lambda^2}\\[3mm]
\dfrac{-\lambda y^2+\lambda
	z^2-2z\lambda^2+\lambda^3+z-\lambda}{y(-y^2+z^2-2z\lambda+\lambda^2)}
\end{array}\right),
\end{equation}
introduced in \cite{GSS07}, is responsible for the spectral problem associated with this group. It is conjugate to a \ronerthree{simpler} map 
\begin{equation}\label{eq:map img con}
\left(\begin{array}{c}
y\\z\\
\lambda
\end{array}\right)\mapsto \left(\begin{array}{c}
\dfrac
zy\\[3mm]
\dfrac \lambda y (-2+y\lambda)\\[3mm]
\dfrac 1 \lambda (-y -y\lambda^2 -\lambda)
\end{array}\right),
\end{equation}
but, basically this is all that is known about this map.


\section{Self-Similar Groups}\label{sec:self similar groups}

Self-similar groups arise from actions on $d$-regular rooted trees $T_d$ (for $d\geq 2$), while self-similar (operator) algebras arise from $d$-similarities $\psi \colon H \to H^d = H \oplus \hdots  \oplus H$ on \ronertwo{an} infinite dimensional Hilbert space $H$, where $\psi$ is an isomorphism.

Let us begin with \ronerthree{the} definition of a self-similar group. Let $X = \{x_1, \hdots,x_d\}$ be an alphabet, $X^*$ be the set of finite words (including the empty word $\mtset$) ordered lexiographically (assuming ${x_1 < x_2< \hdots < x_d}$), and let $T_d$ be a $d$-regular rooted tree with the set of vertices $V$ identified with $X^*$ and set of edges $E = \{(w,wx) \mid w \in X^*, x \in X\}$. The Figure \ref{fig:T_2} shows the binary rooted tree when $X = \{0,1\}$. Usually we will omit the index $d$ in $T_d$. The root vertex corresponds to the empty word.

From a geometric point of view, the boundary $\partial T$ of the tree $T$  consists of infinite paths (without back tracking) joining the root vertex with infinity. It can be identified with the set $X^\N$ of infinite words (sequences) of symbols from $X$ and equipped with the Tychonoff product topology which makes it homeomorphic to a Cantor set. Let $\Aut T$ be the group of automorphisms of $T$ (i.e., of bijections on $V$ that \ronertwo{preserve} the tree structure). The cardinality of $\Aut T$ is $2^{\aleph_0}$ and this group supplied with a natural topology is a profinite group (i.e., compact totally disconnected topological group, or a projective limit $\displaystyle{\lim_{n \to \infty} \Aut T_{[n]}}$ of group of automorphisms of finite groups $\Aut T_{[n]}$, where $T_{[n]}$ is \ronertwo{the} finite subtree of $T$ from the root until the $n$-th level, for $n = 1,2,\hdots$).

Symmetric group $\Sy_d$ ($\cong \Sy (X))$ of permutations on $\{1,2,\hdots,d\}$ naturally acts on $X$ and on $V$ by $\sigma(xw) = \sigma(x)w$, for $w \in V$ and $\sigma \in \Sy_d$. \ronerthree{That is,} a permutation $\sigma$ permutes vertices of the first level according to its action on $X$ and no further action below first level. For $v \in V$ let $T_v$ be the subtree of $T$ with the root at $v$.

For each $v \in V$, the subtree $T_v$ is naturally isomorphic to $T$ and the corresponding isomorphisms $\alpha_v 
\colon T \to T_v$ constitute a canonical system of self-similarities of $T$. Any automorphism $g \in \Aut T$ can be described by a permutation $\sigma$ showing how $g$ acts on the first level and a $d$-tuple of automorphisms $g_{x_1}, \hdots, g_{x_d}$ of trees $T_{x_1}, \hdots, T_{x_d}$ showing how $g$ acts below the first level. As $T \cong T_{x_i}$ this leads to the isomorphism
\begin{equation}\label{eq:tree identification}
\Aut T \stackrel{\psi}{\cong} \left( \Aut T \times \hdots \times \Aut T \right) \rtimes \Sy_d,
\end{equation}
where $\rtimes$ \ronertwo{denotes the} operation of semidirect product (recall that if $N$ is a normal subgroup in a group $G$, $H$ is a subgroup in $G$, $N\cap H = \{e\}$, and $NH = G$, then $G = N \rtimes H$). Another interpretation of \ronerthree{the} isomorphism \eqref{eq:tree identification} is
\begin{equation}\label{eq:tree wreath identification}
\Aut T \cong \Aut T \wr_X \Sy_d 
\end{equation}
where $\wr_X$ denotes \ronertwo{the} permutational wreath product \cite{GNS00,Nek05}. According to \eqref{eq:tree identification}, for $g \in \Aut T$,
\begin{equation}\label{eq:wreath recursions}
\psi(g) = (g_1,\hdots,g_d)\sigma.
\end{equation}

\ronertwo{Relations} of this sort are called wreath recursions and elements $g_x, x\in X$ are called sections.

If $V_n = X^n = X \times  \hdots \times X$ denotes \ronerthree{the} $n$-th level of the tree, then every $g \in \Aut T$ preserves the level $V_n$, for $n = 1,2,\hdots$. Thus the maximum possible transitivity of a group $G \leq \Aut T$ is the level transitivity.

\begin{Definition}\label{def:self-similar gp}
	A group $G$ acting on a tree $T (=T_d)$ by automorphism is said to be self-similar if for all $g \in G, x \in X$ the section $g_x$ coming from wreath recursion \eqref{eq:wreath recursions} belongs to $G$ after identification of $T_x$ (on which $g_x$ acts) with $T$ using identifications $\alpha_x^{-1} \colon T_x \to T$.
\end{Definition}

An alternative way to define self-similar groups is via Mealy automata (also known as the transducers or the sequential machines. See \cite{AS03,BLRS09} for more applications of automata).

A non-initial Mealy automaton $\A = \gen{Q,X,\pi,\lambda }$ consists of \ronerthree{a} finite alphabet $X=\{x_1,\hdots,x_d\}$, a set $Q$ of states, a transition function $\pi \colon Q\times X \to Q$, and an output function $\lambda \colon Q \times X \to X$. Selecting a state $q \in Q$ as initial, produces the initial automaton $\A_q$. The functions $\pi$ and $\lambda$ naturally extends to $\pi \colon Q \times X^* \to Q$ and $\lambda \colon Q \times X^* \to X^*$ via inductive definitions
\begin{align*}
	\pi(q,xw) & = \pi (\pi(q,x),w), \\
	\lambda(q,xw) & = \lambda(q,x) \lambda(\pi(q,x),w),
\end{align*}
for all $w \in X^*$. Thus the initial automaton $\A_q$ determines the maps
\begin{align}
	X^n & \to X^n, \quad n=1,2, \hdots \label{eq:automaton mao n}\\
	X^\N & \to X^\N \label{eq:automaton map natural}  
\end{align}
which we will denote also by $\A_q$ (or sometimes even by $q$). Moreover, $\A_q$ induces an endomorphism of the tree $T = T(X)$ via identification of $V$ with $X^*$. The automaton $\A$ is said to be finite if $|Q| <\infty$.

The initial automaton $\A_q$ can be schematically viewed as the sequential machine (or the transducer) shown in Figure~\ref{fig:transducer}. At the zero moment $n=0$, the automaton $\A_q$ is in the initial state $q=q_0$, reads the symbol $x_0$, produces the output $y_0 = \lambda(q,x_0)$, and moves to the state  $q_1 = \pi(q,x_0)$. Then $\A_q$ continues to operate with input symbols in the same fashion until \ronerthree{reading} the last symbol $x_n$.
\tikzstyle{decision} = [diamond, draw, fill=blue!20, 
text width=4.5em, text badly centered, node distance=3cm, inner sep=0pt]
\tikzstyle{block} = [rectangle, draw,  
text width=5em, text centered,  minimum height=4em]
\tikzstyle{line} = [draw, -latex']
\tikzstyle{cloud} = [ellipse, node distance=4cm,
minimum height=2em]
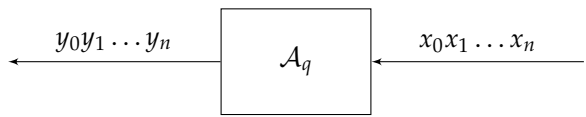
\begin{figure}[H]
	\centering
	\begin{tikzpicture}[node distance = 2cm, auto]
	\node [block] (init) {$\A_q$};
	\node [cloud, left of=init] (expert) {};
	\node [cloud, right of=init] (system) {};
	\path [line] (init) -- node [above]{$y_0 y_1 \hdots y_n$}(expert);
	\path [line] (system) -- node [above]{$x_0x_1 \hdots x_n$}(init);
	\end{tikzpicture}
	\caption{Transducer or sequential machine.}\label{fig:transducer}
\end{figure}

An automaton  of this type is called a synchronous automaton. Asychronous automata can also be defined and used in group theory and coding as explained in \cite{GNS00} and \cite{Lev64}.

An automaton $\A$ is invertible if for any $q \in Q$, the map $\pi(q,\cdot) \colon X \to X$ is a bijection, i.e., $\pi(q,\cdot)$ is an element $\sigma_q$ of the symmetric group $\Sy(X)$. Invertibility of $\A$ implies that for any $q\in Q$ the initial automaton $\A_q$ induces an automorphism of the tree $T=T(X)$. The compositions $\A_q \circ \B_s$ of maps $\A_q, \B_s \colon X^n \to X^n$, where $\B_s = \gen{S,X,\pi',\lambda'}$ is another automaton over the same alphabet, is the map determined by the automaton $\mathcal{C}_{(q,s)} = \A_q \circ \B_s$ with the set of states $Q \times S$ and the transition and \supun{output} functions determined by $\pi, \pi', \lambda, \lambda'$ in the obvious way ({see Figure} \ref{fig:composition}).
If $\A$ is an invertible automaton, then for any $q \in Q$, the inverse map also is determined by an automaton, which will be denoted by  $\A^{-1}_q$.
\tikzstyle{block} = [rectangle, draw, text width=5em, text centered,  minimum height=4em]
\tikzstyle{line} = [draw, -latex']
\tikzstyle{cloud} = [ellipse]
\begin{figure}[H]
	\centering
	\begin{tikzpicture}[node distance = 2cm, auto]
	\node [block] (as) {$\A_q$};
	\node [cloud, left of=as, node distance=2cm] (lmid) {};
	\node [cloud, left of=as, node distance=4cm] (lend) {};
	\node [block, right of=as, node distance=5cm] (bs) {$\B_s$};
	\node [cloud, right of=bs, node distance=2cm] (rmid) {};
	\node [cloud, right of=bs, node distance=4cm] (rend) {};
	\node [cloud, right of=as, node distance=2.5cm] (mid) {};
	\node [block, below of=mid] (comp) {$\A_q \circ \B_s$};
	\path [line] (as) -- node [above]{$z_0 z_1 \hdots z_n$}(lend);
	\path [line] (bs) -- node [above]{$y_0 y_1 \hdots y_n$}(as);
	\path [line] (rend) -- node [above]{$x_0 x_1 \hdots x_n$}(bs);
	\path [line] (comp) -| (lmid);
	\path [line] (rmid) |- (comp);
	\end{tikzpicture}
	\caption{{Composition of automata.}}\label{fig:composition} 
\end{figure}
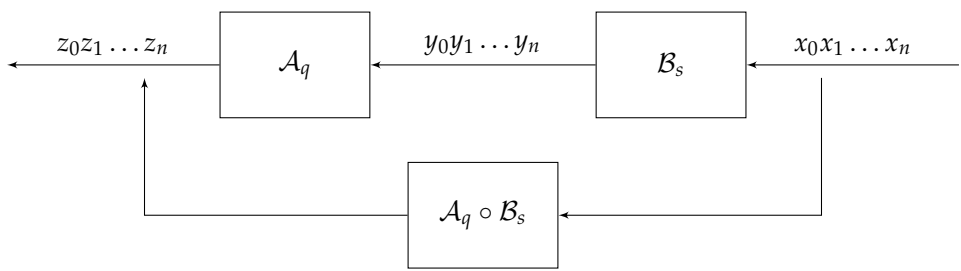

The above discussion shows that for each $m=2,3,\hdots$, we have a semigroup $\mathcal{FAS}(m)$ of finite initial automata over the alphabet on $m$ letters. \ronerthree{We can also} define a group $\mathcal{FAG}(m)$ of finite invertible initial automata. The group $\mathcal{FAG}(m)$ naturally \ronerthree{embeds} in $\mathcal{FAG}(m+1)$. These groups are quite complicated, contain many remarkable subgroups, and depend on $m$. At the same time in the asynchronous case there is only one (up to isomorphisms) group, introduced in \cite{GNS00}, called the group of rational homeomorphisms of a Cantor set. This group, for instance, contains famous R.~Thompson's groups $F,T$ and $V$. In fact, the elements in $\mathcal{FAS}(m)$ and $\mathcal{FAG}(m)$ are classes of equivalence of automata, usually presented by the minimal automaton. The classical algorithm of minimization of automata solves the word problem in $\mathcal{FAS}(m)$ and $\mathcal{FAG}(m)$.

A convenient way to present finite automata is by diagrams, of the type shown on Figure \ref{fig:automata}.
The nodes (vertices) of a such diagram correspond to the states of $\A$, each state $q \in Q$ has $|X|$ outgoing edges of the form
\begin{center}
	\begin{minipage}{0.3\textwidth}
		\includegraphics[width=\linewidth]{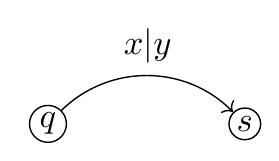}
	\end{minipage}\
	\qquad or
	\begin{minipage}{0.3\textwidth}
		\includegraphics[width=\linewidth]{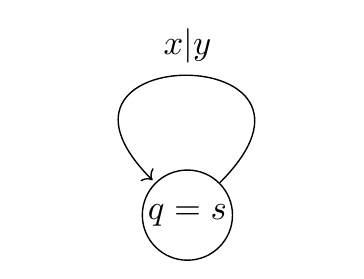}
	\end{minipage}
\end{center}
indicating that if current state is $q$ and the input symbol is $x$, then the next state will be $s$ and the output will be $y$. \ronerthree{This way} we describe the transition and the output functions simultaneously.

If $\A$ is an invertible automaton, then we define $S(\A)$ and $G(\A)$, the semigroup and  the group generated by $\A$;
\[S(\A) = \gen{\A_q \mid q\in Q}_{sem.}\]
is the semigroup generated by initial automata $\A_q, q \in Q$ and 
\begin{align*}
	G(\A) & = \gen{\A_q \mid q\in Q}_{gr.} \\
	& = \gen{\A_q \cup \A_q^{-1} \mid q\in Q}_{sem.}
\end{align*}
is the group generated by $\A_q, q \in Q$.

\ronerthree{The group} $G(\A)$ acts on $T = T(X)$ by automorphisms and for each $q\in Q$ the wreath recursion \eqref{eq:wreath recursions} becomes
\[\A_q = (\A_{\pi(q_0,x_1)}, \hdots ,\A_{\pi(q,x_d)}) \sigma_g. \]
Therefore the groups generated by \ronerthree{the} states of the invertible automaton are self-similar.

The opposite is also true, any self-similar group can be realized as $G(\A)$ for some invertible automaton, only the automaton could be infinite. \ronertwo{Self-similar groups generated by finite automata, called fractal groups (see \cite{BGN03}), constitute an interesting class of groups.} Study of fractal groups in many cases leads to study of fractal objects, as explained for instance in \cite{BGN03,Nek05,GNS15}.


\section{Self-Similar Algebras}\label{sec:self similar algebras}

The definition of self-similar algebras resembles the definition of self-similar groups. It is based on the important property of infinite dimensional Hilbert space $H$ to be isomorphic to the direct sum of $d$ (for $d\geq 2)$ copies of it. A $d$-fold similarity of $H$ is an isomorphism
\[\psi \colon H \to H^d = H \oplus \hdots  \oplus H.\]
There are many such isomorphisms and they are in a natural \ronertwo{bijection} with the $*$-representations of the Cuntz algebra $\OO_d$ as observed in \cite{GN07}. The Cuntz algebra is given by the presentation
\begin{equation}\label{eq:Cuntz relations}
\OO_d \cong \left< a_1, \hdots, a_d \mid a_1 a_1^* + \hdots + a_d a_d^* = 1, a_i^* a_i = 1, i = 1,\hdots,d \right>
\end{equation}
by generators and relations that we will call Cuntz relations.

\begin{Theorem}[Proposition 3.1 from \cite{GN07}]
	The relation putting into correspondence to a $*$-representation $\rho \colon \OO_d \to \B(H)$ \ronertwo{of $\OO_d$} into \ronerthree{the} $C^*$-algebra $\B(H)$ of bounded operators on a separable infinite dimensional Hibert space $H$, the map $\tau_\rho = (\rho(a_1^*), \hdots, \rho(a_d^*))$, where $a_1, \hdots,a_d$ are generators of $\OO_d$\ronertwo{,} is a bijection between the set of representations of $\OO_d$ on $H$ and the set of $d$-fold self-similarities.
	
	The inverse of this bijection puts into correspondence to a $d$-similarity $\psi \colon H \to H^d$ the $*$-representation of $\OO_d$ given by $\rho(a_k) = T_k$, for
	\begin{equation}\label{eq:correspondence star rep to d similarity}
	T_k(\xi) = \psi^{-1}(0,\hdots,0,\xi,0,\hdots,0),
	\end{equation}
	where $\xi$ in the right hand side is at the $k$-th coordinate of $H^d$.
\end{Theorem}

A natural example of a $d$-similarity comes from the $d$-regular rooted tree $T$ and its boundary $\partial T$ supplied by uniform Bernoulli measure $\mu$. That is, $\partial T \cong X^\N$, where $X = \{x_1, \hdots, x_d\}$ and $\mu$ is the uniform distribution on $X$. Then $H = L^2(\partial T, \mu)$ decomposes as
\[ \bigoplus_{x\in X} L^2(\partial T_x, \mu_x),\]
where $T_x$ is the subtree of $T$ with the root at  the vertex $x$ of the first level and $\mu_x = \mu |_{\partial T_x}$. Then ${L^2(\partial T_x, \mu_x) \cong L^2(\partial T, \mu)}$ via the isomorphism given by the operator $U_x \colon L^2(\partial T_x, \mu_x) \to L^2(\partial T, \mu)$,
\[U_x f(\xi) = \frac{1}{\sqrt{d}} f(x\xi), \quad \xi \in \partial T_i, \quad x\in X.\]

Another example associated with the self-similar subgroup $G < \Aut T$ would be to consider a countable self-similar subset $W \subset \partial T$, i.e., a subset $W$ such that $W = \sqcup_{x\in X} xW$. Such a set can be obtained by including the orbit $G\xi$, $\xi \in \partial T$ into the set $W$ that is self-similar closure of $G\xi$. Then
\[\ell^2(W) = \bigoplus_{x\in X} \ell^2(xW)\]
and $\ell^2(xW)$ is isomorphic to $\ell^2(W)$ via the isomorphism $U_x \colon \ell^2(xW) \to \ell^2(W)$ given by
\[U_x(f)(w) = f(xw).\]

Let $G$ be a self-similar group acting on the $d$-regular tree $T=T(X)$, $|X| =d$, and $\psi\colon H \to H^X$ be a $d$-fold similarity. The unitary representation $\rho$ of $G$ on $H$ is said to be self-similar with respect to $\psi$ if for all $g \in G$ and for all $x \in X$
\begin{equation}\label{eq:self-similat unitary rep}
\rho(g) T_x = T_y \rho(h),
\end{equation}
where $T_x$ is the operator defined by \eqref{eq:correspondence star rep to d similarity}, the element $h$ is the section $g|_x$ of $g$ at the vertex $x$ of the first level, and $y = g(x)$, for each $x\in X$.

The meaning of the relation \eqref{eq:self-similat unitary rep} comes from the wreath recursion \eqref{eq:wreath recursions} and its generalization represented by the relation 
\[ g(xw) = g(x) g|_x(w), \quad \forall w \in W.\]

Examples of self-similar representations are the Koopman representation $\kappa$ of $G$ in $L^2(\partial T, \mu)$ (that is, $\left( \kappa (g) f \right) (x) = f(g\inv x)$ for $f \in L^2(\partial T, \mu)$) and permutational representations in $\ell^2(W)$ given by the action of $G$ on the self-similar subset $W \subset \partial T$.

The papers \ronertwo{\cite{GN07,Nek04,Nek09} introduce and discuss} a number of self-similar operator algebras associated with self-similar group\ronertwo{s}. They are denoted by $\A_{max}, \A_{min}$ and if $\A_\rho$ is the algebra obtained by the completion of the group algebra $\CC[G]$ with respect to a self-similar representation $\rho$, then there are natural surjective homomorphisms
\[\A_{max} \to \A_\rho \to \A_{min}.\]
The definition of $\A_{max}$ involves a general theory of Cuntz-Pimsner $C^*$-algebras, which was developed in \cite{Pat99}.

Study of the algebra $\A_\rho$ is based on the matrix recursions. A matrix recursion on an associative algebra $A$ is a homomorphism
\[\varphi \colon A \to M_d(A),\]
where $M_d(A)$ is the algebra of $d\times d$ matrices with entries in $A$.

Wreath recursions \eqref{eq:wreath recursions} associated with a self-similar representation $\rho$ of a self-similar group  $G$ naturally lead to a matrix recursion $\varphi$ for the group algebra $\CC[G]$. Define $\varphi$ on group elements $g \in G$ by
\begin{equation}\label{eq:gp element to matrix}
\varphi (g) = \left( A_{y,x} \right)_{x,y \in X},
\end{equation}
where
\begin{equation}\label{eq:matrix to gp element}
A_{y,x} = 
\begin{cases}
\rho(g|_x) & \text{ if } g(x) = y \\
0 & \text{ otherwise}
\end{cases},
\end{equation}
and extended to the group algebra $\CC[G]$ and its closure $\A_\rho = \overline{\rho(\CC[G])}$ linearly.

In terms of the associated representation $\rho$ of the Cuntz algebra we have the relations
\[ g|_x = T^*_{g(x)} g T_x, \]
where $T_x = \rho(a_x)$, and $a_x$ are generators of $\OO_d$ from presentation \eqref{eq:Cuntz relations}.

\begin{Example}
	In the case of the group $\GG = \gen{a,b,c,d}$, given by presentation \eqref{eq:1st gp presentation}, acting on $T_2$ via the automaton in Figure \ref{subfig:gri}, the recursions for the Koopman representation are, 
	\begin{align}\label{eq:matrix map 1st gp}
		a & = 
		\left(
		\begin{array}{cc}
			0 & 1 \\
			1 & 0 \\
		\end{array}
		\right), &
		b & = 
		\left(
		\begin{array}{cc}
			a & 0 \\
			0 & c \\
		\end{array}
		\right), &
		c & = 
		\left(
		\begin{array}{cc}
			a & 0 \\
			0 & d \\
		\end{array}
		\right), &
		d & = 
		\left(
		\begin{array}{cc}
			1 & 0 \\
			0 & b \\
		\end{array}
		\right),
	\end{align}
	where $1$ stands for the identity operator. Here (and henceforth), we abuse the notation and write $g$ in place of $\kappa(g)$, for any group element $g$.
\end{Example}

\begin{Example}
	In the case of the Basilica group $\B = \gen{a,b}$, the recursions for the Koopman representation are,
	\begin{equation*}\label{eq:matrix map Basilica}
		a = 
		\left(
		\begin{array}{cc}
			1 & 0 \\
			0 & b \\
		\end{array}
		\right),
		\quad
		b  = 
		\left(
		\begin{array}{cc}
			0 & 1 \\
			a & 0 \\
		\end{array}
		\right).
	\end{equation*}
\end{Example}

The minimal self-similar algebra $\A_{min}$ is defined using the algebra generated by permutational representation in $\ell^2(W)$ for an arbitrary self-similar subset $W \subset \partial T$ spanned by the orbit of any $G$-regular point \cite{GN07}. A point $\xi \in \partial T$ is $G$-regular if $g\xi \neq \xi$ or $g\zeta = \zeta$ for all $\zeta \in U_\xi$ for some neighborhood $U_\xi$ of $\xi$. Points that are not $G$-regular are called $G$-singular. In the case of ${\GG = \gen{a,b,c,d}}$, $\GG$-singular points constitute the orbit $\GG(1^\infty)$. The set of $G$-regular points is co-meager (i.e., an intersection of a countable family of open dense sets). This notion was \ronertwo{introduced} in \cite{GNS00} and now play an important role in numerous studies.

Another example of a self-similar algebra is $\A_{mes}$ generated by Koopman representation $\kappa$ on $L^2(\partial T, \mu)$. The representation $\kappa$ is the sum of finite dimensional representations and $\A_{mes}$ is residually finite dimensional \cite{BG00,Gri05,BGH03}. The algebra $\A_{mes}$ has a natural self-similar trace $\tau$, i.e., a trace that satisfies 
\[\tau(a) = \frac{1}{d} \sum_{i=1}^d \tau(a_{ii})\]
for $a \in \A_{mes}$, where $\varphi(a) = (a_{ij})_{i,j = 1}^d$.

This trace was used in \cite{GZ01} to compute the spectral measure associated with the Laplace operator on the Lamplighter group. The range of values of $\tau(g)$ for $g\in \GG$ is $\frac{1}{7} \{\frac{m}{2^n} \mid m,n \in \N, 0 \leq m \leq 2^n \}$ \cite{Gri11}.

A group $G$ is said to be just-infinite if $G$ is infinite and every proper quotient of $G$ is finite. An algebra $\mathcal{C}$ is just-infinite dimensional if it is infinite dimensional but every proper quotient is finite dimensional. Infinite simple groups and infinite dimensional simple $C^*$-algebras are examples of such objects. Infinite cyclic group $\Z$, infinite dihedral group $D_\infty$, and $\GG$ are examples of just-infinite groups. There is a natural partition of the class of just-infinite groups into the class of just-infinite branch groups, hereditary just-infinite groups, and near-simple groups \cite{Gri00}. The proof of this result uses a result of J. Wilson from \cite{Wil71}. Roughly speaking, a branch group is a group acting in a branch way on some spherically homogeneous rooted tree $T_{\overline{m}}$, given by a sequence $\overline{m} = \{m_n \}_{n=1}^\infty$ (where, $m_n \geq 2$ for all $n$) of integers (the integer $m_n$ is called the branching number for vertices of $n$-th level). The group $\GG$ is an example of a just-infinite branch group and the algebra $\varphi(\CC[\GG])$, which sometimes (following S. Sidki \cite{Sid97}) is called the ``thin'' algebra, is just-infinite dimensional \cite{Bar06}.
\begin{Problem}
	Is the $C^*$-algebra $\mathcal{C}_\kappa^*(\GG)$, generated by Koopman representation of $\GG$ in $L^2(\partial T, \mu)$, just-infinite dimensional?
\end{Problem}

There is also a natural partition of separable just-infinite dimensional $C^*$-algebras into three subclasses by the structure of its space of primitive ideals which can be of one of the types $Y_n$, $0\leq n \leq \infty$, where type $Y_0$ means a singleton and corresponds to the case of simple $C^*$-algebras, the type $Y_n$, $1\leq n < \infty$ corresponds to an essential extension of a simple $C^*$-algebra by a finite dimensional $C^*$-algebra with $n$-simple summands, and in the $Y_\infty$ case the algebras are residually finite dimensional. If $\mathcal{C}_\kappa^*(\GG)$ happens to be just-infinite dimensional, this would be a good addendum to the examples presented in \cite{GMR18}, \ronertwo{where} the above trichotomy for $C^*$-algebras is proven.

Given $G$, a self-similar group acting on $T(X)$, $|X|=d$, the associated universal Cuntz-Pimsner $C^*$-algebra $\OO_G$, denoted as $\A_{max}$,  is defined as the universal $C^*$-algebra generated by $G$ and $\OO_d$ satisfying the following relations:
\begin{enumerate}[leftmargin=*,labelsep=4.9mm]
	\item Relations of $G$,
	\item Cuntz relations \eqref{eq:Cuntz relations},
	\item $ga_x = a_y h $ for $g,h \in G$ and $x,y \in X$ if $g(xw) = y h(w)$ for all $w \in X^*$ (i.e., if $g(x) = y$ and $h = g|_x$ is a section).
\end{enumerate}

A self-similar group $G$ is said to be contracting if there exists a finite set $\mathcal N \subset G$ such that for all $g \in G$, there exists $n_0 \in \N$ with $g|_w \in \mathcal N$ for all words $w\in X^*$ of length greater than $n_0$. The smallest set $\mathcal N$ having this property is called the nucleus. Examples of contracting groups are the adding machine $\alpha$ given by the relation $\psi(\alpha) = (1, \cdots, 1, \alpha) \sigma$, where $\sigma$ is a cyclic permutation of $X$, group $\GG$, Basilica, Hanoi groups, and $IMG(z^2+i)$. The Lamplighter group $\mathcal{L}$ presented by the automaton in Figure \ref{subfig:lamplighter} as well as the examples given by automata from Figures \ref{subfig:free3}, \ref{subfig:bellaterra} are not contracting.

The contracting property of the group is \ronerthree{a} tool used to prove subexponentiality of the growth. But not all contracting groups grow subexponentially. For instance, Basilica is contracting but has exponential growth. For a contracting group $G$ with the nucleus $\mathcal N$, the Cuntz-Pimsner algebra $\OO_G$ has the following presentation by generators and relations:
\begin{enumerate}
	\item Cuntz relations,
	\item relations $\displaystyle{g = \sum_{x \in X} a_{g(x)} g|_x a^*_x}$ for $g \in \NN$,
	\item relations $g_1g_2g_3=1$, when $g_1,g_2,g_3 \in \NN$ and relations $gg^* = g^*g =1$ for $g \in \NN$ \cite{Nek09}.
\end{enumerate}
Thus, for contracting groups, the algebra $\A_\infty$ is finitely presented.

The nucleus of the group $\GG$ is $\NN = \{1,a,b,c,d\}$ hence $\A_\infty(\Gr)$ is given by the presentation
\begin{align*}
	\A_\infty(\GG) = \left< a_0,a_1,a,b,c,d \right. \mid  &  1 = a_0 a_0^* + a_1 a_1^* = a_0^* a_0 = a_1^* a_1 = aa^* = bb^* = cc^* =dd^* = bcd,  \\
	& \left. b = a_0aa_0^* + a_1ca_1^*, c = a_0aa_0^* + a_1da_1^*, d = a_0a_0^* + a_1 b a_1^* \right>. 
\end{align*}

For contracting level transitive groups with the property that for every element $g$ of the nucleus the interior of the set of fixed points of $g$ is closed\ronertwo{,} all self-similar completions of $\CC[G]$ are isomorphic\ronerthree{,} and the isomorphisms $\A_{max} \cong \A_{mes} \cong \A_{min}$ hold. This condition holds for Basilica but \ronertwo{does} not hold for the group $\GG$ and one of the groups $G_\omega$, $\omega \in \{0,1,2\}^\N$ from \cite{Gri84}, presented by the sequence $(01)^\infty$ (and studied by A. Erschler \cite{Ers04}). In the latter case $\A_{max} \ncong \A_{min}$. It is unclear at the moment if for $\GG$ we have the equality $\A_{max} \cong \A_{min}$.

\ronertwo{Representations and characters of self-similar groups of branch type are considered in \cite{DG18,DG17Irreducibility}. On self-similarity, operators, and dynamics see also \cite{MT03}.}


\section{Joint Spectrum and Operator Systems}\label{sec:joint spectrum}

Given a set $\{M_1, \hdots, M_k\}$ of bounded operators in a Hilbert (or more generally Banach) space $H$, one can consider \ronerthree{the} pencil of operators $M(\ol z) = z_1 M_1 + \hdots + z_k M_k$ for $\ol z = (z_1, \hdots,z_k) \in \CC^k$, and define the joint (projective) spectrum $\jspec(M(\ol z))$ as the set of parameters $\ol z \in \CC$ for which $M(\ol z)$ is not invertible. Surprisingly, such a simple concept did not attract much attention when the number of non-homogeneous parameters is greater or equal to two, until publication \ronerthree{of} \cite{Yan09}. Sporadic examples of analysis of the structure of $\jspec$ and in some cases of computations are presented in \cite{Ant84,BG00,GZ01,GS06,GS07,GS19,AL83,Vin88,ZKK75}. Important examples of pencils come from $C^*$-algebras associated with self-similar groups that were discussed in the previous section.

Recall that given a unital $C^*$-algebra $\A$, a closed subspace $\Sy$ containing the identity element 1 is called an operator system. One can associate to each subspace $\M \subset \A$ an operator system via $\Sy = \M + \M^* + \CC 1$. Such systems are important for the study of completely bounded maps \cite{Pau02,Pis90,Pis01}.

If $G = \gen{a_1, \hdots, a_m}$ is a contracting self-similar group with (finite) nucleus $\NN = \{ n_1, \hdots, n_k\}$ and $\A^*$ is a self-similar $C^*$-algebra associated with a self-similar unitary representation $\rho$ of $G$, then the identity element belongs to $\NN$ and a natural operator space and a self-similar pencil of operators are: $\Sy = span\{\rho(n_1), \hdots, \rho(n_k), \rho(n_1^*), \hdots, \rho(n_k^*) \} \subset \A$ and 
\[ M(\ol z) = \sum_{i=1}^k (\rho(n_i) + \rho(n_i^*)),\]
respectively. The main examples considered in this article come from the group $\Gr$ and the overgroup $\wt \Gr$ that possess the nuclei $\{1,a,b,c,d\}$ and $\{1,a,b,c,d,\wt a,\wt b, \wt c, \wt d \}$, respectively.


\section{Graphs of Algebraic Origin and Their Growth}\label{sec:graphs}

A graph $\Gamma = (V,E)$ consists of a set $V$ of vertices and a set $E$ of edges. The edges are presented by the map $e \colon V\times V \to \N_0$ (\ronerthree{here} $\N_0$ denotes the set of non negative integers), where $e(u,v)$ represents the number of edges connecting the vertex $u$ to the vertex $v$. If $u=v$, then the edges are loops. So what we call a graph in graph theory usually is called a \ronerthree{directed multi-graph or an oriented multi-graph. Depending on the situation, graph can be non-oriented (if the edges are independent of the orientation, i.e., $e(u,v) = e(v,u)$ and edges $(u,v)$ and $(v,u)$ are identified)} and labeled (if edges are colored  by elements of a certain alphabet). We only consider connected locally finite graphs (the later means that each vertex is incident to a finite number of edges). The degree $d(u)$ of the vertex $u$ is the number of edges incident to it \ronerthree{(where each edge from or to $u$ contributes 1 to the degree and each loop contributes 2 to the degree)}. A graph is of uniformly bounded degree if there is a constant $C$ such that $d(v) \leq C$ for all $v \in V$, and $\Gamma$ is a regular graph if all vertices have the same degree.

\begin{figure}[htb]
	\centering
	\begin{subfigure}[H]{.4\textwidth}
		\centering
		\includegraphics[width=0.9\textwidth]{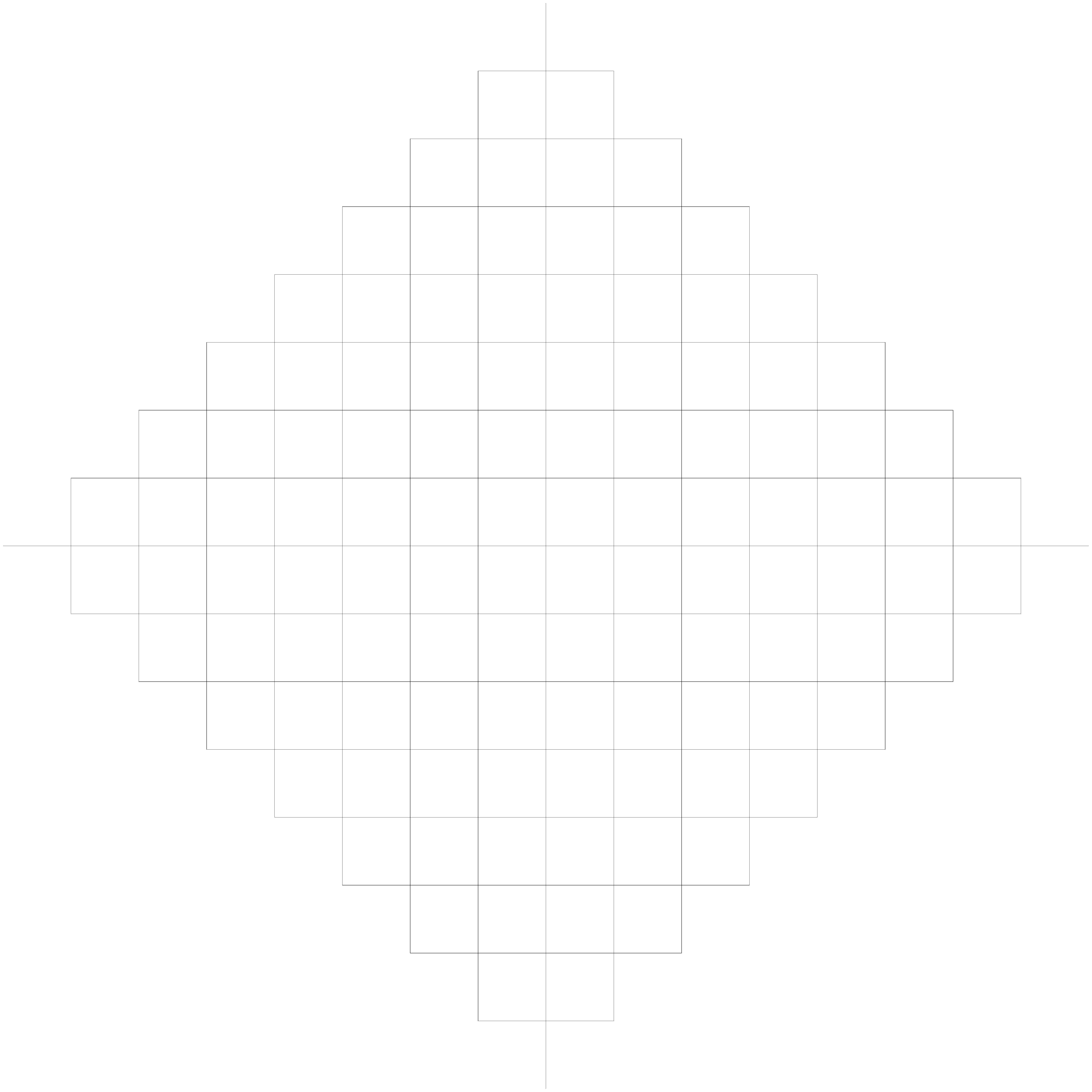}
		\caption{}\label{subfig:cayley_Z2}
	\end{subfigure}
	\qquad
	\begin{subfigure}[H]{.4\textwidth}
		\centering
		\includegraphics[width=0.9\textwidth]{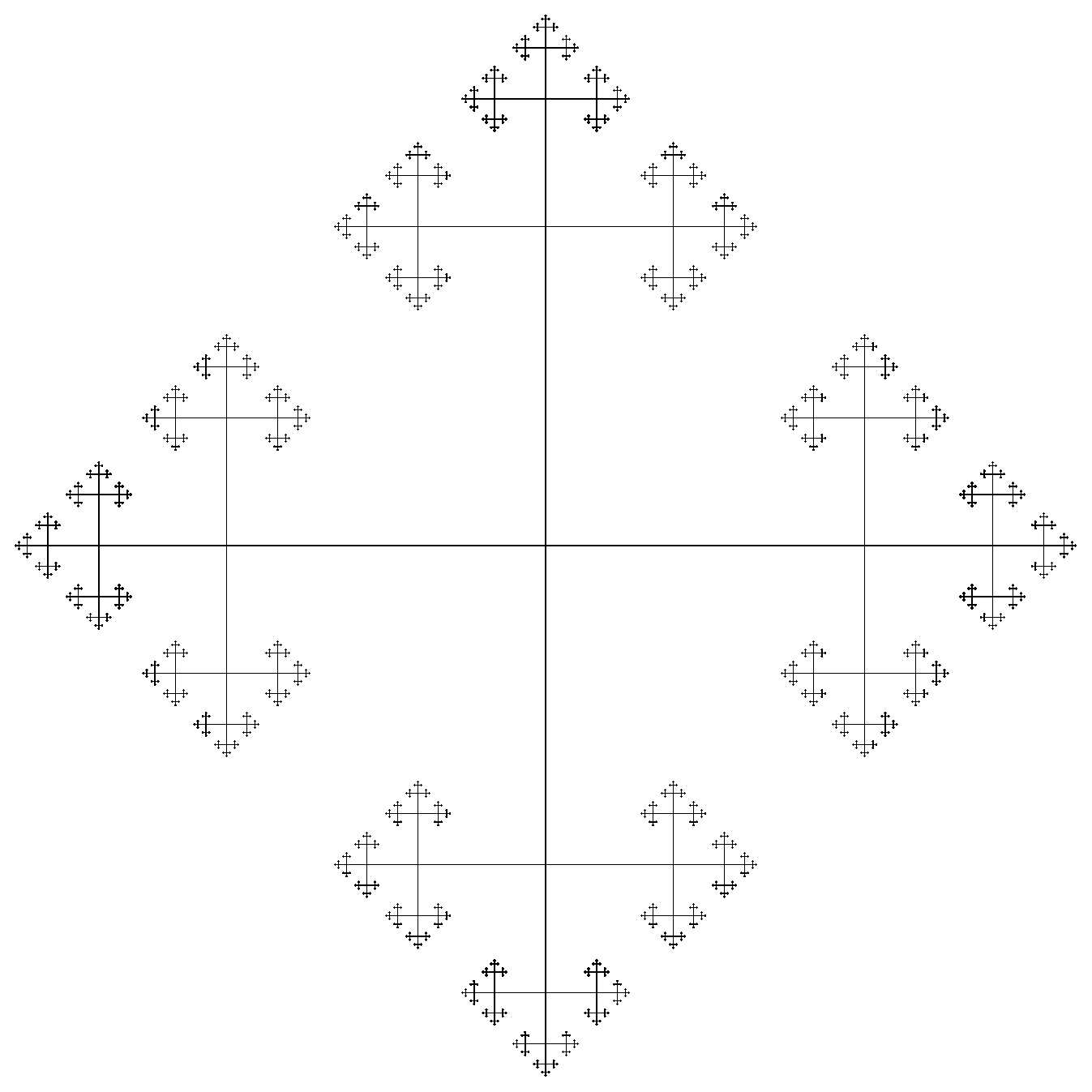}
		\caption{}\label{subfig:cayley_F2}
	\end{subfigure}
	\\ \vspace{.5cm}
	\begin{subfigure}[H]{.4\textwidth}
		\centering
		\includegraphics[width=0.9\textwidth]{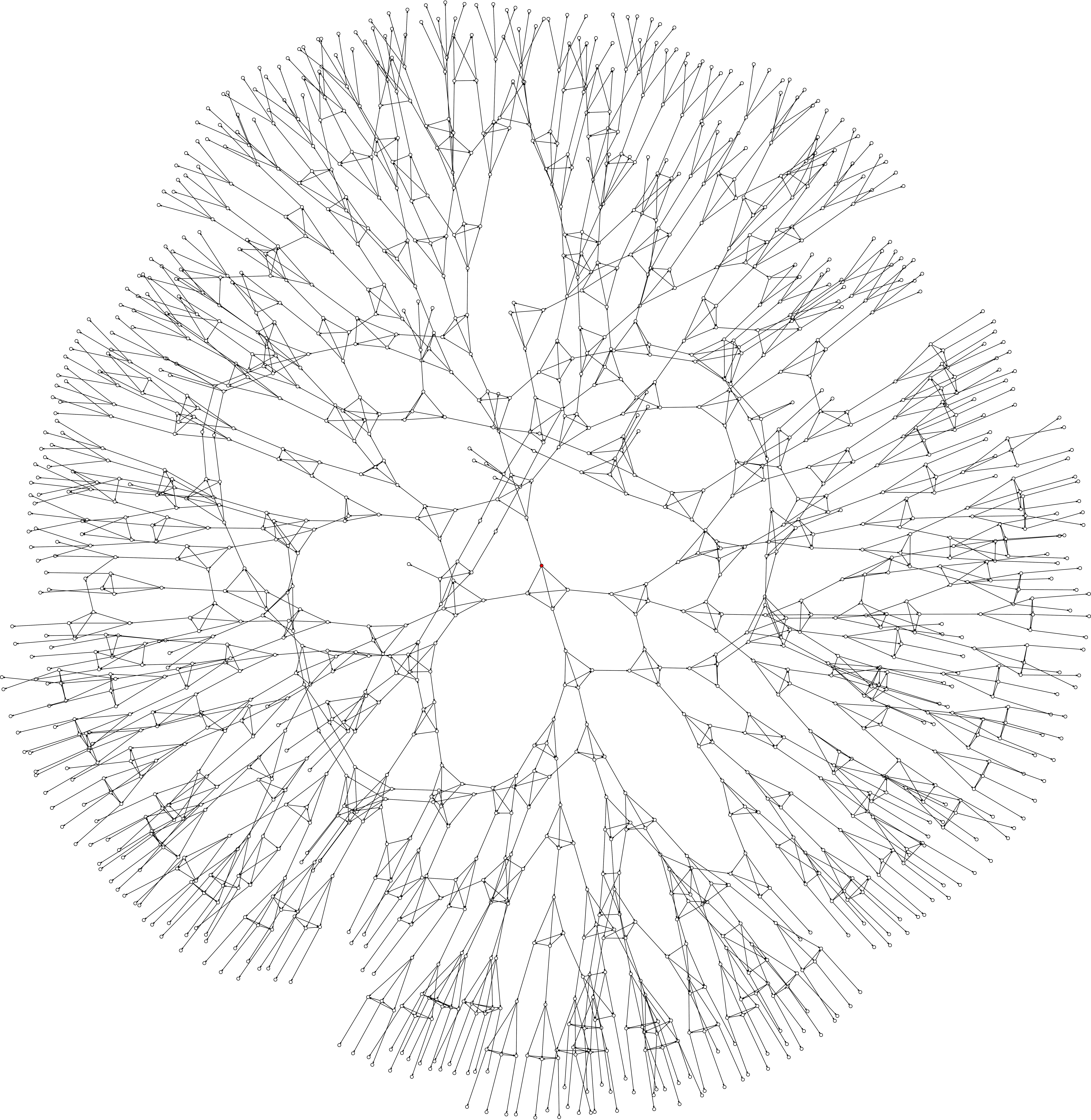}
		\caption{}\label{subfig:cayley_gri}
	\end{subfigure}
	\qquad
	\begin{subfigure}[H]{.4\textwidth}
		\centering
		\includegraphics[width=0.9\textwidth]{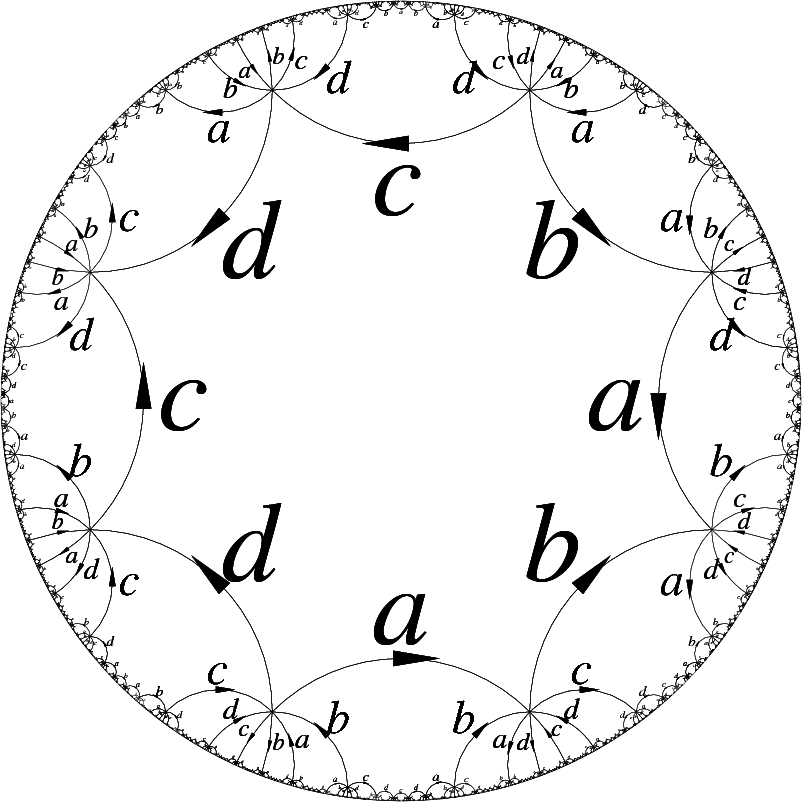}
		\caption{}\label{subfig:cayley_2tor}
	\end{subfigure}
	\caption{Cayley graphs of (\textbf{a}) $\Z^2$, (\textbf{b}) free group of rank 2, (\textbf{c}) group of intermediate growth  $\Gr$, (\textbf{d}) surface group of genus 2.}\label{fig:cayley}
\end{figure}

There is a rich source of examples of graphs coming from groups. Namely, given a marked group $(G,A)$ (i.e., a group $G$ with a generating set $A$, usually we assume that $|A| < \infty$ and therefore the group is finitely generated), one defines the directed graph $\Gamma = \Gamma_l(G,A)$ with $V =G$ and $E = \{(g,ag) \mid a \in A\cup A^{-1} \}$, where $g$ is the origin and $ag$ is the end of the edge $(g,ag)$. This is \ronertwo{the} left Cayley graph. Similarly, one can define \ronertwo{the} right Cayley graph $\Gamma_r(G,A)$, and there is a natural isomorphism $\Gamma_l(G,A) \cong \Gamma_r(G,A)$. Left and right Cayley graphs are vertex transitive, i.e., the group $\Aut(\Gamma)$ of automorphisms acts transitively on the set of vertices (right translations by elements of $G$ on $V=G$ induce automorphisms of $\Gamma_l(G,A)$). When speaking about Cayley graph, \ronerthree{we usually} keep in mind \ronerthree{the} left Cayley graph. Depending on the situation, Cayley graphs are considered as labeled graphs (the edge $e = (g,ag)$ has label $a$), or unlabeled (if labels do not play a role). Cayley graphs can also be converted into undirected graphs by identification of pairs $(g,ag), (ag,a^{-1}(ag)) = (ag,g)$ of mutually inverse pairs of edges. The examples of Cayley graphs are presented in Figure \ref{fig:cayley}. \ronertwo{Non-oriented Cayley graph of $(G,A)$ is $d$-regular with $d = 2\abs{A\setminus A_2}+\abs{A_2}$, where $A_2 \subset A$ is the set of generators whose order is two (involutions)}.

A Schreier graph $\Gamma = \Gamma(G,H,A)$ is determined by a triple $(G,H,A)$, where as before $A$ is a system of generators of $G$ and $H$ is a subgroup of $G$. In this case $V = \{gH \mid g\in G\}$ is a set of left cosets (for the left version of definition) and $E = \{(gH,agH) \mid g \in G, a \in A \cup A^{-1} \}$. Again, one can consider a right version of the definition, oriented or non-oriented, labeled or unlabeled versions of the Schreier graph \cite{NP20,DDMN10,BDN17}. 

\begin{figure}[H]
	\centering
	\begin{subfigure}[b]{.8\textwidth}
		\centering
		\includegraphics[width=\textwidth]{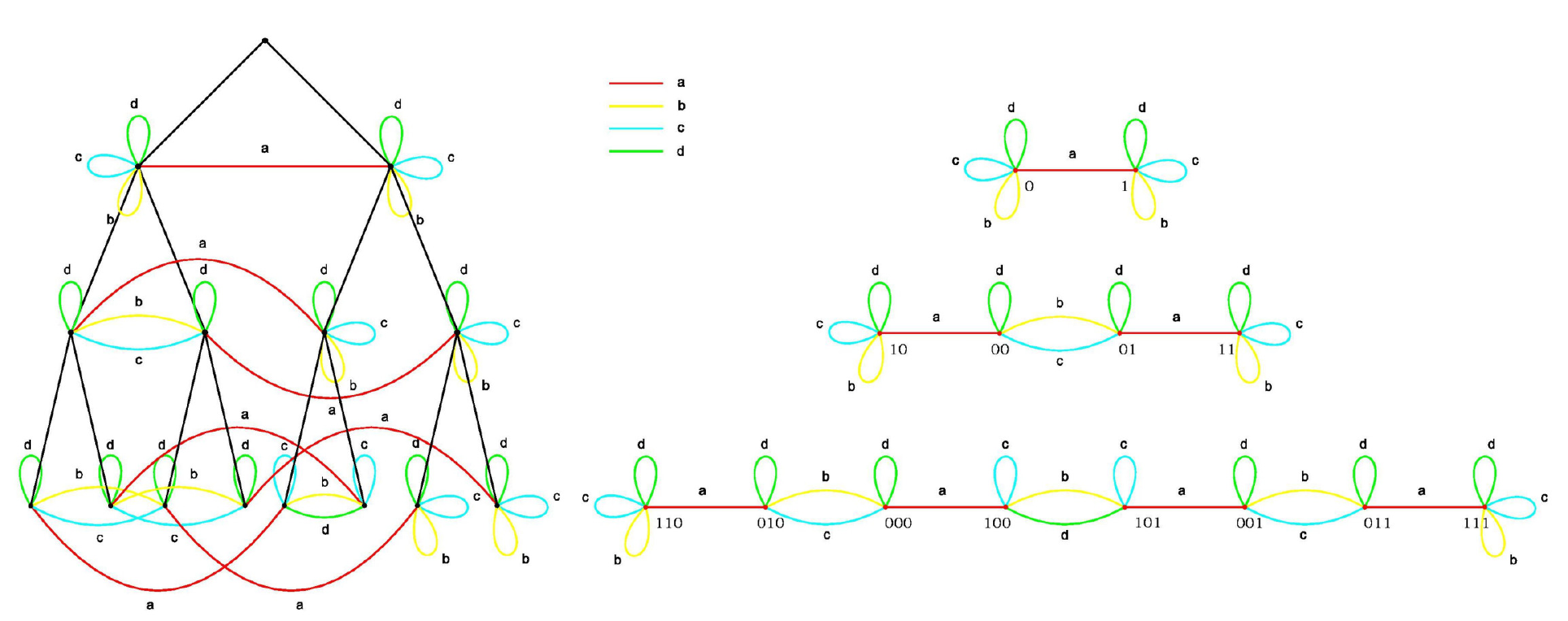}
		\caption{}\label{subfig:schreier_gri}
	\end{subfigure}
	\\ \vspace*{8pt}
	\begin{subfigure}[b]{.8\textwidth}
		\centering
		\includegraphics[trim={0 0 0 5.5cm},clip,width=\textwidth]{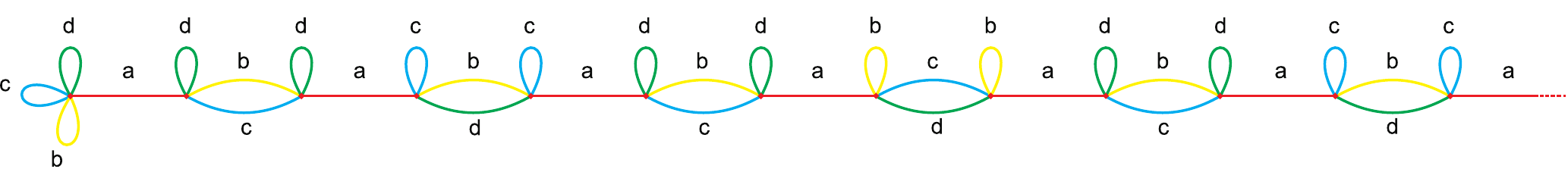}
		\\ \vspace*{8pt}
		\includegraphics[trim={0 0 0 5.5cm},clip,width=\textwidth]{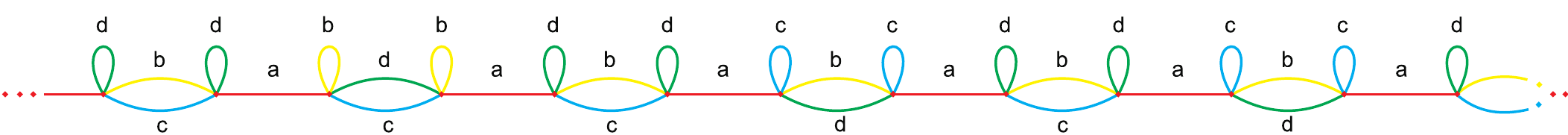}
		\caption{}\label{subfig:schreier_gri_boundary}
	\end{subfigure}
	\\ \vspace*{8pt}
	\begin{subfigure}[b]{\textwidth}
		\centering
		\includegraphics[width=.2\textwidth]{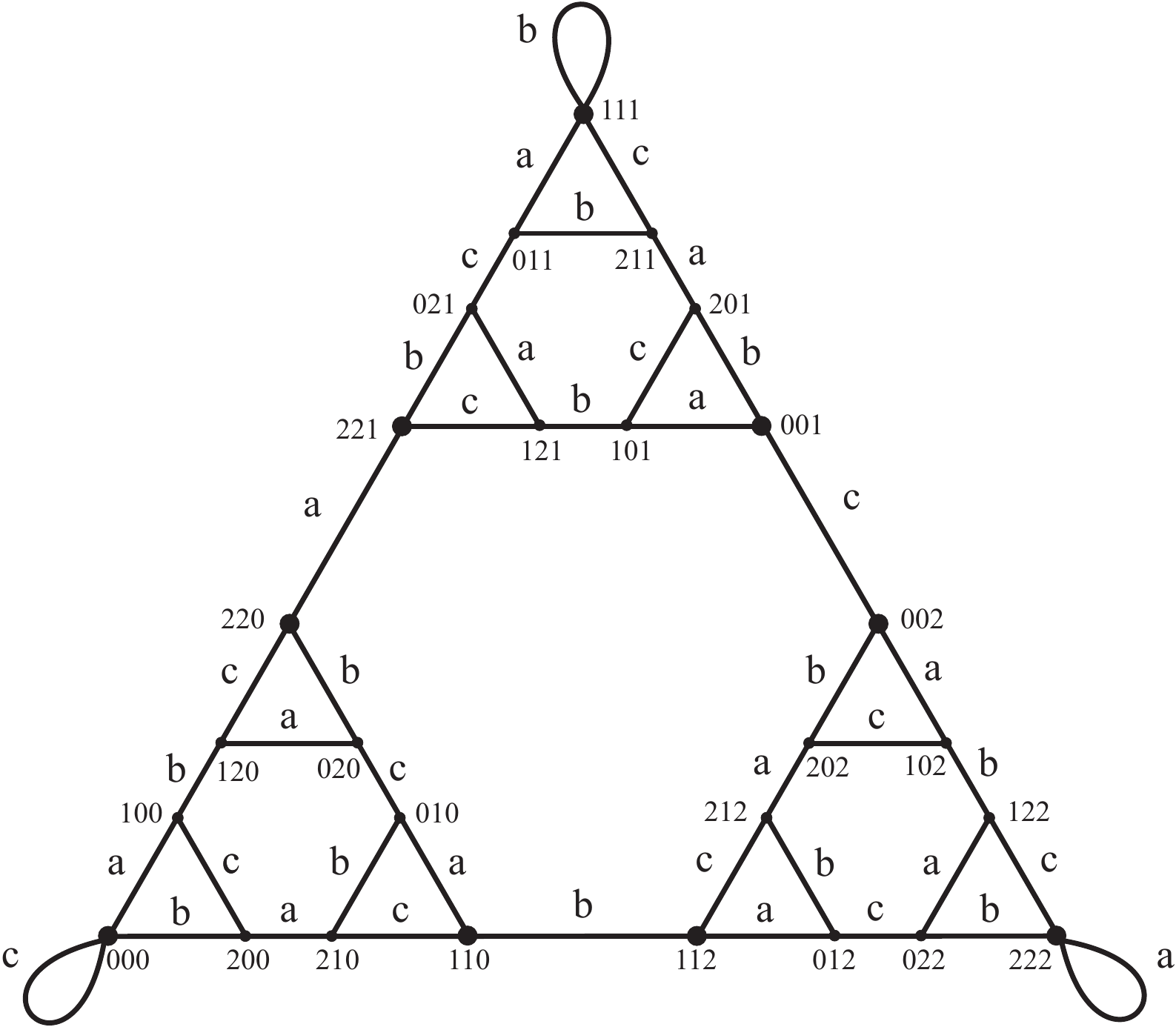}
		\caption{}\label{subfig:schreier_hanoi3}
	\end{subfigure}
	\\ \vspace*{8pt}
	\begin{subfigure}[b]{.6\textwidth}
		\centering
		\includegraphics[width=\textwidth]{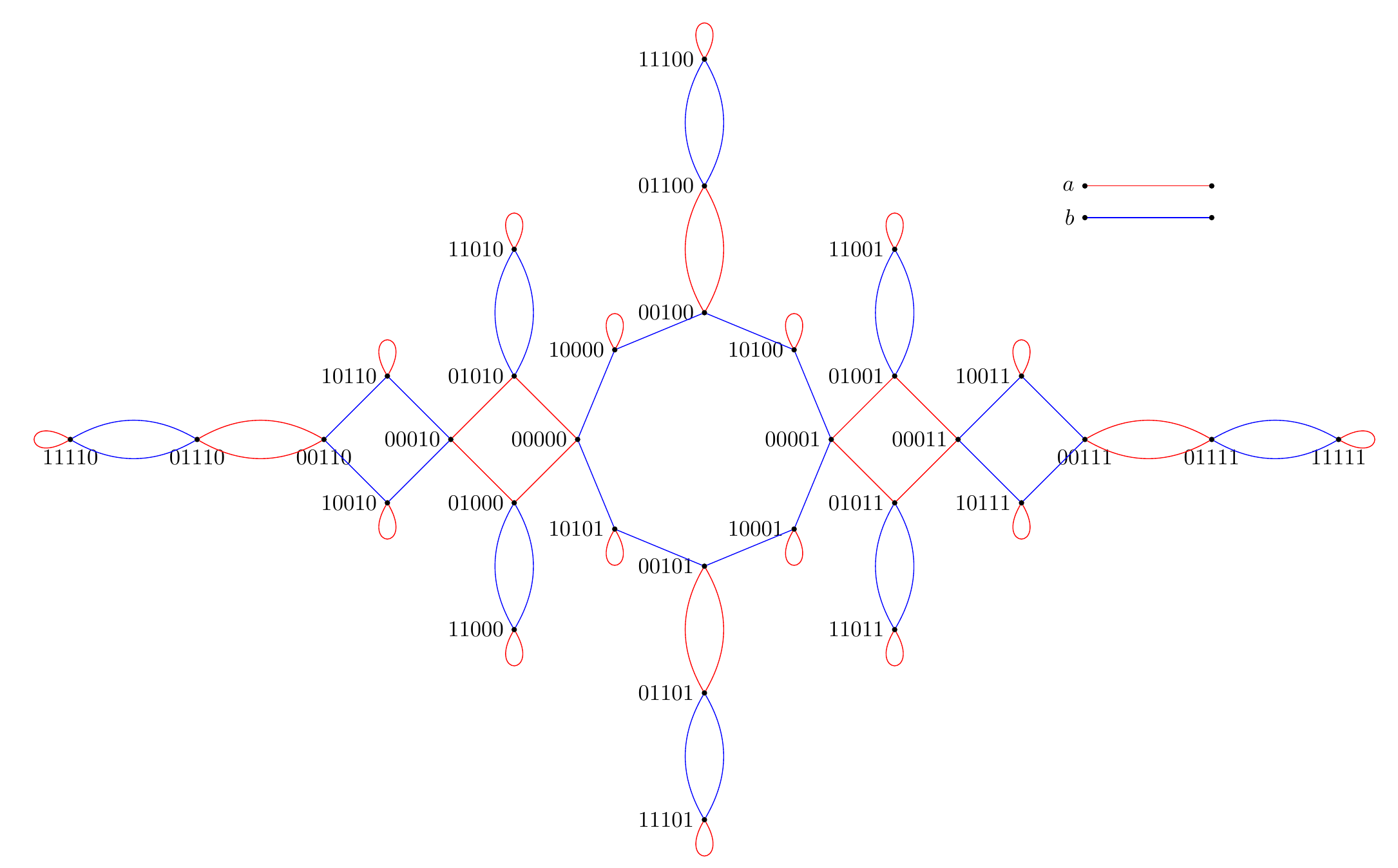}
		\caption{}\label{subfig:schreier_basilica}
	\end{subfigure}
	\vspace*{8pt}
	\caption{Schreier graphs of (\textbf{a}) $\Gr$ (finite), (\textbf{b}) $\Gr$ (infinite and bi-infinite) (\textbf{c}) $\Hh^{(3)}$, (\textbf{d}) Basilica.}\label{fig:schreier}
\end{figure}

Cayley graph $\Gamma(G,A)$ is isomorphic to the Schreier graph $\Gamma(G,H,A)$ when $H = \{1\}$ is the trivial subgroup. \ronertwo{Non-oriented Schreier graphs are also $d$-regular with $d$ given by the same expression as above}, but in contrast with Cayley graphs, they may have a trivial group of automorphism. Examples of Schreier graphs are presented in the Figure \ref{fig:schreier}.

We have the following chains of classes of graphs:
\begin{gather*}
	\{\text{locally finite}\} \supset \left\{\begin{array}{c}\text{bounded}\\ \text{degree}\end{array}\right\}  \supset \{\text{regular}\} \supset \left\{\begin{array}{c}\text{vertex}\\ \text{transitive}\end{array}\right\} \supset \{\text{Cayley}\},\\
	\{\text{regular}\} \supset  \{\text{Schreier}\} \ronerthree{\supset \{\text{Cayley}\}}.
\end{gather*}

\ronertwo{In fact the class of $d$-regular graphs of even degree $d =2m$ \ronertwo{coincides} with the class of Schreier graphs of the free group $F_m$ of rank $m$} (for finite graphs this was observed by Cross \cite{Lub95} and for the general case see \cite{Har00} and Theorem 6.1 in \cite{Gri11}).
For \ronerthree{an} odd degree, the situation is slightly more complicated, but there is clear understanding on which of them are Schreier graphs \cite{Lee20}.

Schreier graphs have much more applications in mathematics being able to provide a geometrical-combinatorial representation of many objects and situations. In particular, they are used to approximate fractals, Julia sets, study the dynamics of groups of iterated monodromy, Hanoi Tower Game on $d$ pegs for $d \geq 3$, etc.

Growth function of a graph $\Gamma = (V,E)$ with distinguished vertex $v_0 \in V $ is the function \[ \gamma(n) = \gamma_{\Gamma,v_0}(n) = \#\{ v \in V \mid d(v_0,v) \leq n \}, \]
where $d(u,v)$ is the combinatorial distance given by the length of a shortest path connecting two vertices $u$ and $v$. Its rate of growth when $n \to \infty$, \ronerthree{defines} the rate of growth of the graph at infinity (if $\Gamma$ is an infinite graph). It does not depend on the choice of $v_0$ (in case of connected graphs), and is bounded by the exponential function $C d^n$ when $\Gamma$ is of uniformly bounded degree $\leq d$.

The growth of Schreier graph  can be of power function $n^\alpha$, $\alpha >0$ type, even with irrational $\alpha$ \cite{BG00}, of the type $a^{(\log n)^\alpha}$ \cite{GS06}, and of many other unusual types of growth.

On the other hand, the growth of \ronerthree{a} Cayley graph (or what is the same the growth of the corresponding group) is much more restrictive. It is known that if it is of the power type $n^\alpha$, then $\alpha$ is a positive integer (in which case the group is said to be of polynomial growth) and the group is virtually nilpotent (i.e., contains nilpotent subgroup of finite index) \cite{Gro81}. The Cayley graphs of virtually solvable groups or of linear groups (i.e., groups presented by matrices over a field) either have polynomial or exponential growth \cite{Mil68,Mil68note,Wol68,Tit72}. The question about existence of groups of intermediate growth was raised by J.~Milnor \cite{Mil68Problem} and got the answer in \cite{Gri83,Gri84} \ronerthree{using} the group $\Gr$.

The Cayley graph of the group of intermediate growth $\Gr$ is presented by Figure \ref{subfig:cayley_gri}. The group $\GG$ is a representative of an uncountable family of groups $\GG_\omega = \gen{a,b_\omega,c_\omega,d_\omega}$, $\omega \in \{0,1,2\}^\N = \Omega$, mostly consisting of groups of intermediate growth \cite{Gri84} (definition is given by \eqref{eq:random recursions} and \eqref{eq:gen_gri_gp}). Moreover, there are uncountably many of different rates of growth in this family, where  by the rate (or degree) of growth of a group $G$  we mean the dilatational equivalence class of $\gamma_G(n)$ (two functions $\gamma_1(n), \gamma_2(n)$ are equivalent, $\gamma_1(n) \sim \gamma_2(n)$, if there is $C$ such that $\gamma_1(n) \leq \gamma_2(Cn)$ and $\gamma_2(n) \leq \gamma_1(Cn)$). This gives the first family of cardinality $2^{\aleph_0}$ of continuum of finitely generated groups with pairwise non quasi-isometric Cayley graphs. As shown in  \cite{Bar98,EZ20}, for any $\epsilon >0$,
\[ e^{n^{\beta -\epsilon}} \leq \gamma_{\GG}(n) \leq e^{n^\beta},\]
where $\beta \approx 0.767$ is $1/\log_2 \beta_0$ and $\beta_0$ is the (unique) positive root of the polynomial $x^3-x^2-2x-4$. In \cite{BE14} the group $\GG$ is used to show that for each $\delta$ such that $\beta < \delta <1$, there is a group with growth equivalent to $e^{n^\delta}$.

Surprisingly, so far there is no example of a group with super-polynomial growth but slower than the growth of $\GG$. There is a conjecture \cite{Gri91} that there is a gap in the scale of growth degrees of finitely generated groups between polynomial growth and growth of the partition function $P(n) \sim e^{\sqrt{n}}$ (i.e., if $\gamma_G(n)  \prec e^{\sqrt{n}}$, then $G$ is virtually nilpotent and hence has a polynomial growth). The conjecture \ronerthree{has been} confirmed to be true for the class of groups approximated by nilpotent groups. More on the gap conjecture for group growth and other asymptotic characteristics of groups see \cite{Gri14Onthegap,Gri14Mil}. 

\begin{Problem}
	Is there a finitely generated group with super polynomial growth smaller than the growth of the group $\Gr$?
\end{Problem}


\section{Space of Groups and Graphs and Approximation} \label{sec:space of groups}

A pair $(G,A)$ where $G$ is a group and $A = \{ a_1, \hdots, a_m \}$ is an ordered system of (not necessarily distinct) generators is said to be a marked group. There is a natural topology in the space $\M_m$ of $m$-generated marked groups introduced in \cite{Gri84}. Similarly, there is a natural topology in the space $\mkdschreier m$ of marked Schreier graphs associated with marked triples $(G,H,A)$ (a marked graph is a graph with distinguished vertex viewed as the origin). For every $m \geq 2$, \ronerthree{the spaces} $\M_m$ and $\mkdschreier m$ are totally disconnected compact metrizable spaces whose structure is closely related to various topics of groups theory and dynamics. For instance the closure $\X$ of $\{ \Gr\om \}_{\omega \in \Omega}$ in $\M_4$ consists of a Cantor set $\X_0$ and a countable set $\X_1$ of isolated points accumulating to $\X_0$ and consisting of virtually metabelian groups containing a direct product of copies of the Lamplighter group $\LL = \Z_2 \wr \Z$ \cite{Gri84,BG14} as a subgroup of finite index. The closure of $\{ \wt\Gr\om \}_{\omega \in \Omega}$ is described in \cite{Sam20} and has a more complicated structure. More on spaces $\M_m$ and $\mkdschreier m$ see \cite{Gri05,Cha,OM}.

One of the fundamental questions about these spaces is finding the Cantor-Bendixson rank  (for the definition see \cite{KM04}) characterized by the first ordinal when taking of Cantor-Bendixson derivative does not change the space.

A more general notion than Schreier graph is the notion of orbital graph. Given an action $\alpha$ of marked group $(G,A)$ on \ronerthree{the} space $X$, one can build a graph $\Gamma_\alpha(G,A)$ with the set of vertices $V=X$ and \ronerthree{the set of} edges $E = \{(x,ax) \mid a \in A\cup A\inv \}$. Connected components of this graph are Schreier graphs $\Gamma(G,H_{x_i},A)$, $i \in I$, where $H_{x_i}$ is the stabilizer of point $x_i \in X$ and $\{x_i \}_{i \in I}$ is the set of representatives of orbits. If action is transitive, then orbit graph is a Schreier graph.

\begin{figure}[h]
\centering
\includegraphics[height=.3\textwidth]{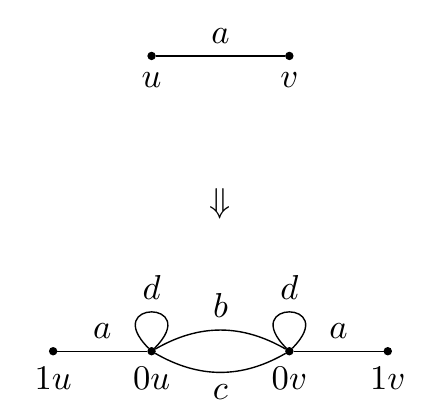}
\includegraphics[height=.3\textwidth]{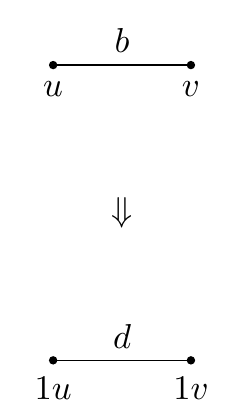}
\includegraphics[height=.3\textwidth]{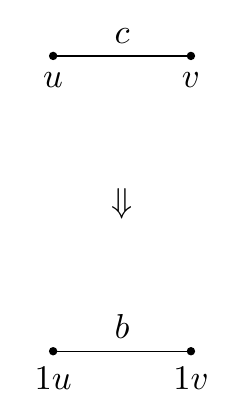}
\includegraphics[height=.3\textwidth]{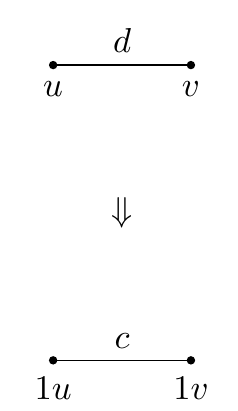}
\vspace*{8pt}
\caption{Graph substitution to obtain $\Gamma_{n+1}$ from $\Gamma_n$ of the group $\Gr$.}\label{fig:graph_sub_gri}
\end{figure}

Given a level transitive action of marked group $(G,A)$ by automorphisms on a $d$-regular rooted tree $T$, one can consider the covering sequence $\{ \Gamma_n \}_{n=1}^\infty$ of graphs where $\Gamma_n$ is orbital graph for action on $n$-th level of the tree and $\Gamma_{n+1}$ covers $\Gamma_n$. Additionally, for every point $\xi \in \partial T$ (the boundary of $T$) one can associate a Schreier graph $\Gamma_\xi = \Gamma(G,H_\xi,A)$ built on the orbit $G(\xi)$ of $\xi$ ($H_\xi = \stab_G(\xi)$). If $v_n$ is a vertex of level $n$ that belongs to the path representing $\xi$, then 
\begin{equation}\label{eq:limit of orbital graphs}
(\Gamma_\xi,\xi) = \lim_{n\to \infty} (\Gamma_n,v_n)
\end{equation}
(the limit is taken in the topology of the space of marked Schreier graphs). The relation \eqref{eq:limit of orbital graphs} \ronerthree{allows an} approximation of infinite graphs by finite graphs that leads also to the approximation of their spectra as shortly explained in the next section.

The example of $\Gamma_n$ and $\Gamma_\xi$ associated with the group $(\Gr, \{a,b,c,d\})$ is given by Figures \ref{subfig:schreier_gri} and \ref{subfig:schreier_gri_boundary}. In this example $\Gamma_{n+1}$ can be obtained from $\Gamma_n$ by substitution rule given by Figure \ref{fig:graph_sub_gri} that \ronerthree{mimics} the substitution $\sigma$ used in presentation \eqref{eq:1st gp presentation}.

Similar property holds for $\Gamma_n, \Gamma_\xi$ associated with the overgroup $\wt\Gr\om$ and and graphs associated with $\Gr\om$ and $\wt\Gr\om$ if instead of a single substitution, to use three substitutions and iterate them accordingly to ``oracle'' $\omega$.

\begin{figure}[h]
	\centering
	\hspace*{-1.5cm}
	\includegraphics[width=0.5\textwidth]{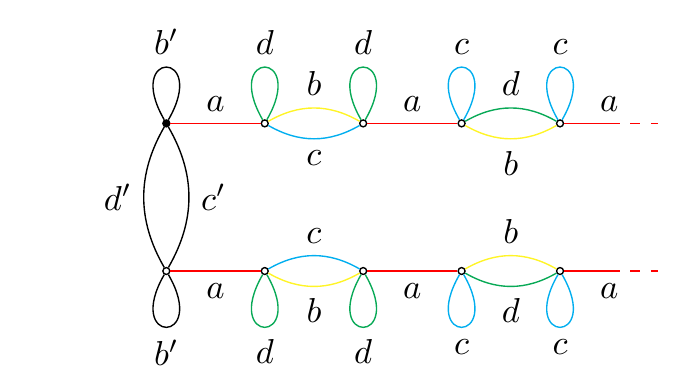}
	\vspace*{8pt}
	\caption{Limit graphs of $\varphi(\partial T)$.}\label{fig:connected Schreier}
\end{figure}

The correspondence $\partial T \ni \xi \overset{\varphi} \longrightarrow (\Gamma_\xi,\xi)$ gives a map from $\partial T$ to $\mkdschreier{4}$ with the image ${\varphi(\partial T) = \{ (\Gamma_\xi,\xi) \mid \xi \in \partial T \}}$. The set of vertices $V_\xi$ of $\Gamma_\xi$ is the orbit $\Gr \xi$ and $\Gr$ acts on $\varphi(\partial T)$ by changing the root vertex $(\Gamma_\xi,\xi) \overset{g} \arr (\Gamma_\xi, g\xi)$. The graphs $(\Gamma_\xi,\xi) , \xi \in \Gr 1^\infty$ are one-ended and are isolated points in $\varphi(\partial T)$ and also in the closure $\overline{\varphi(\partial T)}$ in $\mkdschreier 4$. Deletion of them from $\overline{\varphi(\partial T)}$ gives the set $X = \overline{\varphi(\partial T)} \setminus \{\text{isolated points} \}$ which is a union of $X_0 = \{(\Gamma_\xi,\xi) \mid \xi \in \partial T, \xi \notin \Gr 1^\infty \}$ and the  countable set $X_1$ consisting  of  the  limit points  of  $X_0$  that  do not belong  to  $X_0$.  The  set $X_1$ consists of  pairs $(\Gamma,  v)$,  where  $\Gamma$  is  one  of  the  three  graphs given  by   Figure \ref{fig:connected Schreier}  (where, $(b',c',d')$  is any cyclic permutation of $(b,c,d)$),  and  v  is an  arbitrary  vertex of $\Gamma$ \cite{Vor12}.

Now $X$ also is homeomorphic to a Cantor set, $\Gr$ acts on $X$ and the action is minimal and uniquely ergodic and the map $X_0 \ni (\Gamma_\xi,\xi) \arr \xi \in \partial T$ extends to a continuous factor map
\[\Phi \colon X \to \partial T\]
which is one-to-one except in a countable set of points, where it is three-to-one \cite{Vor12}.


\section{Spectra of Groups and Graphs}\label{sec:spectra}

Let $\Gamma = (V,E)$ be a $d$-regular (non-oriented) graph. \ronerthree{The} Markov operator $M$ acts \ronertwo{on the} Hilbert space $\ell^2(V)$ (which we denote by $\ltwo$) and is defined as
\[(Mf)(x) = \frac{1}{d} \sum_{y \sim x} f(y),\]
where $f \in \ltwo$ and $x \sim y$ is the adjacency relation. The operator $L = I -M$ where $I$ is the identity operator is called \ronerthree{the} discrete Laplace operator. Operators $M$ and $L$ can be defined also for non-regular graphs as it is done for instance in \cite{MW89,Chu97}. The Markov operator $M$ is a self-adjoint operator with the norm $\norm{M} \leq 1$ and the spectrum $\spec(M) \subset [-1,1]$. The name ``Markov'' comes from the fact that $M$ is \ronerthree{the} Markov operator associated with \ronertwo{the} random walk on $\Gamma$ in which \ronerthree{a} transition $u \to v$ occurs with probability $p = \frac{1}{d}$, if $u$ and $v$ are adjacent vertices. Random walks on graphs are special case of Markov chains.

A graph $\Gamma$ of uniformly bounded degree is called amenable if $1\in \spec(M)$ ($\iff \norm{M} = 1$). Such definition comes from the analogy with von-Neumann -- Bogolyubov theory of amenable groups, i.e., groups with invariant mean \cite{Gre69}. By Kesten's criterion \cite{Kes59}, amenable groups can be characterized as groups for which the spectral radius $r = \overline{\lim}_{n\to \infty} \sqrt[n]{P_{1,1}^{(n)}}$ \ronertwo{is equal to one}, where $P_{1,1}^{(n)}$ is the probability of return to identity element in $n$ steps  of the simple random walk on the Cayley graph.

By a spectrum of a graph (or a group), we mean the spectrum of $M$. A more general concept is \ronerthree{used} when graph is \emph{weighted}, in the sense that a weight function $w \colon E \to \R$ on edges is \ronertwo{given} and the ``weighted Markov'' operator $M_w$ is defined in $\ltwo$ as 
\[ (M_w f)(x) = \sum_{y \sim x} w(x,y)f(y). \]
A special case of such situation is given by a marked group $(G,A)$ (i.e., group $G$ together with its generating set $A$) and symmetric probability distribution $P(g)$ on $A \cup A^{-1}$: $P(a) = P(a^{-1})$, for all $a \in A$ and $\sum_{a\in A} P(a) = 1/2$. Then Markov operator $M_P$ acts as 
\[(M_P f)(g) = \sum_{a \in A \cup A^{-1}} P(a) f(ag)\]
and $M_P$ is the operator associated with a random walk on the (left) Cayley graph $\Gamma(G,A)$, where transition $g \to ag$ holds with probability $P(a)$.

The case of uniform distribution on $A \cup A^{-1}$ (i.e., of a simple random walk) is called isotropic case, while non-uniform distribution corresponds to the anisotropic case.

Basic questions about spectra of infinite graphs are:
\begin{enumerate}[leftmargin=*,labelsep=4.9mm]
	\item What is the shape (up to homeomorphism) of $\spec(M)$?
	\item What can be said about spectral measures $\mu_\varphi$ associated with functions $\varphi \in \ltwo$, in particular, with delta functions $\delta_v$, $v \in V$?
\end{enumerate}

As usual in mathematical physics, the gaps in the spectrum, discrete and singular continuous parts of the spectrum are of special interest. Also \ronerthree{an} important case is when graph $\Gamma$ has a subgroup of \ronerthree{the} group of automorphisms acting on the set of vertices freely and co-compactly (i.e., with finitely many orbits). The case of vertex transitive graphs and especially of the Cayley graph is of special interest and is related to many topics in abstract harmonic analysis, operator algebras, asymptotic group theory and theory of random walks. Among open problems, let us mention the following.

\begin{Problem}
	Can the spectrum of a Cayley graph of a finitely generated group be a Cantor set? (i.e., homeomorphic to a Cantor set). The problem is open in both isotropic and anisotropic cases.
\end{Problem}

\begin{Problem}
	Can a torsion free group have a gap in spectrum?
\end{Problem}

If the answer to the last question is affirmative, then this would give a counterexample to the Kadison - Kaplanski conjecture on idempotents.

Spectral theory of graphs of algebraic origin is a part of spectral theory of convolution operators in $\ltwoG$ or $\ltwoG \oplus \hdots \oplus \ltwoG$ given by elements of a group algebra $\CC[G]$ or $n \times n$ matrices with entries in $\CC[G]$ and is closely related to many problems on $L^2$-invariants, including $L^2$-Betti numbers, Novikov-Slubin invariants, etc.

The state of art of the above \ronerthree{problem is roughly} as follows. Spectra of Euclidean grids (lattices~$\cong \Z^d$, $d \geq 1$) \ronerthree{and} of their perturbations is a classical subject based on the use of Bloch-Floquet theory, representation theory of abelian groups and classical methods. The main facts include finiteness of the number of gaps in spectrum, band structure of the spectrum, absence of singular continuous spectrum, finiteness (and in many cases) absence of the discrete part in the spectrum \cite{BK12}.

\ronertwo{Another important case is trees and tree like graphs, such as, the Cayley graphs of free groups and free product of finite groups.} For these groups the use of representation theory is limited but somehow possible, the structure of the spectrum is similar to the case of graphs with co-compact $\Z^d$-action, although the methods are quite different, see \cite{Kes59,Car72,CC98,FN91,FP82,FP83,KLW13,KLW15,KPT84,Woe09,Woe86}.

\ronerthree{One more important case constitute graphs associated with classical and non-classical self-similar fractals or self-similar groups \cite{BG00,MT95,SS20}.}
In particular, in \cite{BG00} it is shown that

\begin{Theorem}
	Spectrum of the Schreier graph of self-similar group can be a Cantor set or a Cantor set and a countable set of isolated points accumulated to it.
\end{Theorem}

Recall that in Section \ref{sec:space of groups}, for a group acting on rooted tree $T$, we introduced a sequence $\{ \Gamma_n\}_{n=1}^\infty$ of finite graphs and family $\{ \Gamma_\xi\}_{\xi \in \partial T}$ of infinite graphs. In the next result $M_n$ is a Markov operator associated with $\Gamma_n$ and $\kappa$ is the Koopman representation.

\begin{Theorem}\label{thm:spec limit}
	Let $G$ be a group acting on rooted tree $T$ and $\Sigma = \spec(\kappa(M))$, where $M = \sum_{a \in A} (a+a\inv) \in \Z[G]$. Then,
	\begin{enumerate}[leftmargin=*,labelsep=4.9mm]
		\item $\Sigma= \overline{ \cup_n sp(M_n)}. \hfill \refstepcounter{equation}\left(\theequation\right)\label{eq:spectrum approximation}$
		\item If \ronerthree{the} action is level transitive and $G$ is amenable, then $\spec(\Gamma_\xi)$ does not depend on the point $\xi \in \partial T$ and is equal to $\Sigma$.
		\item The limit
		\begin{equation}\label{eq:measure limit}
		\mu_* = \lim_{n\to \infty} \mu_n
		\end{equation}
		of counting measures
		\[ \mu_n = \frac 1 {|\Gamma_n|} \sum_{\lambda \in \spec(M_n)} \mass{\lambda} \]
		(summation is taken with multiplicities) exists. \ronerthree{It is called the} density of states.
	\end{enumerate}
\end{Theorem}
This result is a combination of observations made in \cite{BG00Hecke} and \cite{GZ04}. In fact, the relation \eqref{eq:spectrum approximation} and the fact about the existence of limit hold not only for elements $M = \sum_{a\in A} (a + a\inv)$ of the group algebra (i.e., after normalization corresponding to a simple random walk, i.e., isotropic case) but for arbitrary self-adjoint element of the group algebra.

In \cite{DG17}, the following result related to Theorem \ref{thm:spec limit} is proved.
\begin{Theorem}
	Let $(X,\mu)$ be a measure space and a group $G$ act by transformations preserving the class of measure $\mu$ (i.e., $g_*\mu< \mu$ for all $g \in G$). Let $\kappa \colon G \to U(L^2(X,\mu))$ be the Koopman representation: $\kappa(g)f(x) = \sqrt{\dv{g_*\mu}{\mu}(x)} f(g\inv x)$ and $\pi$ be groupoid representation. Then for any $M \in \CC[G]$,
	\begin{equation}\label{eq:spec measure}
	\spec(\kappa(M)) \supset \spec(\Gamma_x) = \spec(\pi(M))
	\end{equation}
	$\mu$-almost surely ($\Gamma_x$ is the orbital graph on $Gx$) and if the action $(G,X,\mu)$ is amenable (i.e., partition on orbits is hyperfinite) then in \eqref{eq:spec measure} we have equality $\mu$-almost surely instead of inclusion.
\end{Theorem}

\ronertwo{For definition of groupoid representation, we direct reader to \cite{DG17}. Amenable actions are discussed in \cite{Kai97,KM04}}.

The idea of approximation of groups in the situation of group action was explored by A.~Stepin and A.~Vershik in the 1970s \cite{Ste84}. Now the corresponding group property is called local \ronerthree{embeddability} into finite group (LEF). This is a weaker form of the classical residual finiteness of groups. \ronertwo{For instance, topological full groups \cite{GM14} mentioned in Section \ref{sec:dynamical system} are LEF}.

Approximation of Ising model in infinite residually finite groups by Ising models on finite quotients is suggested in \cite{GS84}. \ronertwo{The Ising model on self-similar Schreier graphs by finite approximations was studied in \cite{DDN11} and the dimer model on these graphs was studied in \cite{DDN12}.}

A recent result of B. Simanek and R. Grigorchuk \cite{GS19} shows that,
\begin{Theorem}
	\ronertwo{There are Cayley graphs with infinitely many gaps in the spectrum.}
\end{Theorem}

It is a direct consequence of the next theorem.

\begin{Theorem}[\cite{GS19}]
	Let $\LL$ be the lamplighter group $\Z_2 \wr \Z$. There is a system of generators $\{a,b\}$ of $\LL$ such that the convolution operation $M_\mu$ in $\ell^2(\LL)$ determined by the element $a+a^{-1}+b+b^{-1}+\mu c$ of the group algebra $\R[\LL]$ where $c=b^{-1}a$, $\mu \in \R$ has pure point spectrum. Moreover,
	\begin{enumerate}[leftmargin=*,labelsep=4.9mm]
		\item if $|\mu|\leq1$, the eigenvalues of $M_{\mu}$ densely pack the interval $[-4-\mu,4-\mu]$.
		\item If $|\mu|>1$, the eigenvalues of $M_{\mu}$ form a countable set that densely packs the interval $[-4-\mu,4-\mu]$ and also has an accumulation point $\mu + \dfrac{2}{\mu} \notin [-4-\mu,4-\mu]$.
		\item The spectral measure $\nu_{\mu}$ of the operator $M_{\mu}$ is discrete  and  is given by
		\[
		\nu_{\mu}=\frac{1}{4}\delta_{\mu}+\sum_{k=2}^{\infty}\left[\frac{1}{2^{k+1}}\sum_{\{s:G_{k}(s,\mu)=0\}}\delta_s\right],
		\]
		where
		\[
		G_{k}(z,\mu)=2^k\left[U_k\left(\frac{-z-\mu}{4}\right)+\mu U_{k-1}\left(\frac{-z-\mu}{4}\right)\right],
		\]
		$U_k$ is the degree $k$ Chebyshev polynomial of the second kind.
	\end{enumerate}
\end{Theorem}

Observe that if $\mu \geq 2$ is \ronerthree{an} integer, then the spectrum of the Cayley graph of $\LL$ \ronertwo{built} using the system of generators $a,b,c_1, \hdots, c_\mu$ with $c_1 = \hdots =c_\mu =c$ coincides with the spectrum of $M_\mu$ and hence has infinitely many gaps. This is the first example of a Cayley graph with infinitely many gaps in the spectrum.

In the case when $\mu = 0$, this result was obtained by A. Zuk and R. Grigorchuk \cite{GZ01} and used in \cite{GLSZ00} to answer \ronerthree{a} question of M. Atiyah and to give a counterexample to \ronerthree{a} version of the strong Atiyah conjecture \ronertwo{known in 2001}. This was done by constructing a 7-dimensional closed Reimannian manifold with third $L^2$-Betti number $\beta_2^{(3)} = \frac{1}{3}$. For more on spectra of Lamplighter type groups \ronertwo{and their finite approximations} see \cite{KSS06,GLN16}.

The Schreier spectrum of the Hanoi Tower group $\Hh^{(3)}$ is a Cantor set and a countable  set  of  isolated  points  accumulating  to  it (see Figure \ref{fig:hanoi spectrum}) as follows from the following result.

\begin{Theorem}[\cite{GS06}]\label{thm:Hanoi_spectrum}
	The $n$-th level spectrum $\spec(M_n)$ ($n \geq 1$), as a set, has $3 \cdot 2^{n-1}-1$ points and is equal to
	\[ \{3\} \ \cup \ \bigcup_{i=0}^{n-1} f^{-i}(0) \ \cup \ \bigcup_{j=0}^{n-2} f^{-j}(-2). \]
	The multiplicity of the $2^i$ level $n$ eigenvalues in $f^{-i}(0)$, $i=0,\dots,n-1$ is $a_{n-i}$ and the multiplicity of the $2^j$ eigenvalues in $f^{-j}(-2)$, $j=0,\dots,n-2$ is $b_{n-j}$, where $f(x) = x^3-x-3$ and ${a_i = \dfrac{3^{i-1}+3}{2}}, \quad {b_j = \dfrac{3^{j-1}-1}{2}}$. Moreover, the Schreier spectrum of $\Hh^{(3)}$ (i.e., $\spec(\Gamma_x), x \in \partial T$) is equal to
	\[\overline{\bigcup_{i=0}^\infty f^{-i}(0)}.\]
	It consists of a set of isolated points $\Sigma_0 = \bigcup_{i=0}^\infty f^{-i}(0)$ and its set of accumulation points $\Sigma_1$. The set $\Sigma_1$ is a Cantor set and is the Julia set of the polynomial $f$. The density of states is discrete and concentrated on the set $\bigcup_{i=0}^\infty f^{-i}\{0,-2\}$. Its mass at each point of $f^{-i}\{0,-2\}$ is $\frac 1 {6 \cdot 3^i}$, $i \geq 0$.
\end{Theorem}

\ronertwo{
	Grabowski and Virag \cite{Gra1511,GV15} observed that the Lamplighter group has a system of generators such that the spectral measure is purely singular continuous measure. An infinite family of Schreier graphs with a non-trivial singular continuous part of the spectral measure is constructed in \cite{Per20,GNP}.
}


\section{Schur Complements}\label{sec:Schur}

Schur complement is \ronertwo{a} useful tool in linear algebra, networks, differential operators, applied mathematics, etc. \cite{Cot74}

Let $H$ be a Hilbert space (finite or infinite dimensional) decomposed \ronertwo{into} a direct sum ${H = H_1 \oplus H_2}$, $H_i \neq \{0\}$, $i=1,2$. Let $M\in \B(H)$ be a bounded operator and 
\[M = \left(
\begin{array}{cc}
A & B \\
C & D \\
\end{array}
\right)\]
be a matrix representation of $M$ by block matrices corresponding to this decomposition. Thus
\[A \colon H_1 \to H_1, \qquad B \colon H_1 \to H_2, \qquad C\colon H_2 \to H_1, \qquad D \colon H_2 \to H_2.\]

Two partially defined maps $S_1 \colon \B(H) \to \B(H_1)$ and $S_2 \colon \B(H) \to \B(H_2)$ are defined \ronertwo{by,
	\begin{align*}
		S_1(M) & = A - B D\inv C, \\
		S_2(M) & = D - C A\inv B,
	\end{align*}
	for any $M \in \B(H)$. Note that $S_1(M)$ is defined when $D$ is invertible, and $S_2(M)$ is defined when $A$ is invertible}. Maps $S_1, S_2$ are said to be the Schur complements and the following fact holds.

\begin{Theorem}[\cite{GN07}]
	Suppose $D$ is invertible. Then $M$ is invertible if and only if $S_1(M)$ is invertible and 
	\[ M\inv =
	\left( 
	\begin{array}{cc}
	S_1\inv & -S_1\inv B D\inv \\
	-D\inv C S_1\inv & D\inv C S_1\inv B D\inv + D\inv
	\end{array}
	\right),\]
	where $S_1 = S_1(M)$.
\end{Theorem}

A similar statement holds for $S_2(M)$. The above expression for $M\inv$ is called Frobenius formula.
In the case $\dim H < \infty$, the determinant $|M|$ of matrix $M$ satisfies
\[|M| = |S_1(M)| |D| \]
and the latter relation is attributed to Schur.

There is nothing special in decomposition of $H$ into a direct sum of two subspaces. If ${H = H_1 \oplus \hdots \oplus H_d}$ and 
\[M = \left(
\begin{array}{ccc}
\M_{11} & \dots & M_{1,d} \\
\vdots & \ddots & \vdots \\
M_{d1} & \dots & M_{dd}
\end{array}
\right)\]
for $M_{ij} \colon H_i \to H_j$ and $H = H_1 \oplus H_1^\perp$, where $H_1^\perp = H_2 \oplus \hdots \oplus H_d$, then we are back in the case $d =2$. By change of the order of the summands (putting $H_i$ on the first place) one can define the $i$-th Schur complement $S_i(M)$, for each $i =1,\hdots,d$.

If $\dim H = \infty$ and $\psi \colon H \to H^d$ is a $d$-similarity, then $S_i(M) = (T_i^* M\inv T_i )\inv$, where $T_i = \rho (a_i)$ is the image of the generator $a_i$ of the Cuntz algebra $\OO_d$ under the representation $\rho$ associated with $\psi$ (as explained in Section \ref{sec:self similar algebras}).
Therefore, for each $d \geq 2$, one can define  $\mathscr{S}_d^*$ \ronertwo{the} semigroup generated by the Schur transformations $S_i$, $1 \leq i \leq d$ with the operation of composition. \ronerthree{We will call $\mathscr{S}_d^*$ the Schur semigroup}. For a general element of this semigroup, we get the following expression,
\[S_{i_1} \circ \hdots \circ S_{i_k} (M) = ( (T_{i_k} \hdots T_{I_1} )^* M\inv (T_{i_k} \hdots T_{i_1} ) )\inv\]
(see Corollary 5.4 in \cite{GN07}).

The Schur semigroup $\mathscr{S}_d^*$ consists of partially defined transformations on the infinite dimensional space $\B(H)$. There are  examples of finite dimensional subspaces $L \subset \B(H)$ invariant with respect to $\mathscr S_d^*$, where the restrictions of Schur complements to $L$ generates a semigroup of rational transformations on $L$. An example of this sort is the case of 3-dimensional subspace in $\B(H)$, for $H = L^2(\partial T, \mu)$, where $T$ is \ronertwo{the} binary tree and $\mu$ is \ronertwo{the} uniform Bernoulli measure. It comes from the group $\Gr$ and $L$, the space generated by three operators $\kappa(a), \kappa(b+c+d), I$, where $\kappa$ is the Koopman representation and $I$ is the identity operator. If 
\[x\kappa(a) + y\kappa(u) + zI \in L,\]
where $u = b+c+d$, then in coordinates $x,y,z$, the Schur complements  $S_1,S_2$ are given by
\begin{align}
	S_1 & = 2y \pi(a) + \dfrac{x^2y}{(z+3y)(z-y)} \pi(u) +\left(z+y-\dfrac{z+2y}{(z+3y)(z-y)} \right) I, \label{eq:Schur complement 1}\\
	S_2 & = \dfrac{2x^2y}{(z+3y)(z-y)} \pi(a) + y \pi(u) +\left(z-\dfrac{x^2(z+y)}{(z+3y)(z-y)} \right) I, \label{eq:Schur complement 2}
\end{align}
and the semigroup $\left< S_1,S_2 \right>$ generated by them is isomorphic to the semigroup $\gen{F,G}$ generated by maps \eqref{eq:F map 1stGri}, \eqref{eq:G map 1stGri} (as \eqref{eq:Schur complement 1}, \eqref{eq:Schur complement 2} are homogeneous realizations of $F,G$).
Study of properties of the semigroups $\mathscr S_d^*$ and its restriction on finite dimensional invariant subspaces is a challenging problem.


\section{Self-Similar Random Walks coming From Self-Similar Groups and the M\"{u}nchhausen Trick}

Let $G$ be a self-similar group acting on $T_d$ and $\mu \in \mathcal{M}(G)$ (\ronerthree{where} $\mathcal{M}(G)$ denotes the simplex of probability measures on $G$) be a probability measure whose support generates $G$ (we call such $\mu$ non-degenerate). Using $\mu$ one can define a (left) random walk that begins at $1 \in G$ and transition $g \to h$ holds with probability $\mu(hg^{-1})$. Study of random \ronertwo{walks} on groups is a \ronerthree{large} area initiated by H.~Kesten \cite{Kes59} (see \cite{FS94,GW86,NW02,Saw78,Woe00,Woe87,Woe81} for more on random walks on groups and trees). The main topics of study in random walks are: the asymptotic behavior of the return probabilities $P_{1,1}^{(n)}$ when $n \to \infty$, the rate of escape, the entropy, the Liouville property and the spectral properties of the Markov operator $M$  acting in $\ell^2(G)$ by
\[\left(Mf\right)(g) = \sum_{h\in G} {\mu(h) f(hg)},\]
as was discussed in Section \ref{sec:spectra}.
The case when the measure $\mu$ is symmetric (i.e., $\mu(g)=\mu(g^{-1})$) is of special interest as in this case $M$ is self-adjoint.

A remarkable progress in the theory of random walks on groups was made by L. Bartholdi and B.~Virag. They showed that the self-similarity of a group can be converted into a self-similarity of a random walk on the group and used it for proving \ronerthree{the} amenability of the group. In such a way it was shown that the Basilica is amenable \cite{BS06}. The idea of Bartholdi and Virag was developed by V.~Kaimanovich in terms of entropy and interpreted as a kind of mathematical implementation of the legendary ``M\"{u}nchhausen trick'' \cite{Kai05}.

Let us briefly describe the idea of the self-similarity of random walks. Recall, that a self-similar group is determined by the Mealy invertible automaton, or equivalently, by the wreath recursion \eqref{eq:wreath recursions} coming from the embedding $\psi \colon G \to G \wr_X \Sy_d$.

If $Y^{(n)}$ is a random element of $G$ at the moment $n \in \mathbb{N}$ associated with the random walk determined by $\mu$, then 
\begin{equation}\label{eq:random walk wreath recursion}
\psi(Y^{(n)}) = (Y_1^{(n)}, Y_2^{(n)},\hdots, Y_d^{(n)}) \sigma^{(n)}.
\end{equation}
Let $H = H_i = st_G(i)$ be the stabilizer of $i \in X$. The index $[G\colon H] \leq d$ is finite and hence the random walk hits $H$ with probability 1. Denote by $\mu_H$ the distribution on $H$ given by the probability of the first hit:
\[\mu_H(h) = \sum_{n=0}^\infty {f_{1,h}^{(n)}}\]
where $f_{1,h}^{(n)}$ is the probability of hitting $H$ at the element $h$ for the first time at time $n$.

\ronerthree{Now let us construct transformations $K_i$, $i=1, \hdots,d$ on the space of bounded operators $\B(H)$ in a Hilbert space $H = \ell^2(G)$ when a $d$-similarity $\varphi \colon H \to H^d = H \oplus \hdots \oplus H$ is fixed. Let $S_i$ as before be the $i$-th Schur complement $1\leq i \leq d$ associated with $\varphi$, $\J$ be the map in $\B(H)$, given by $\J(A) = A + I$ (where, $I$ is the identity operator). We define $K_i = \J S_i \J^{-1}$, so if
	\[M = \left(\begin{array}{cc}
	A & B \\
	C & D
	\end{array} \right), \]
	where $A$ is an operator acting on the first copy of $H$ in $H \oplus \hdots \oplus H$, then
	\[K_1(M) = A + B (I-D)^{-1} C. \]
	In the case when $M$ is the Markov operator of the random walk on $G$ determined by the measure $\mu$, this leads to the \ronertwo{analogous} map on simplex $\M(G)$ which we denote by $k_i$. The measure $k_i(\mu)$ is called the $i$-th probabilistic Schur complement.}

\begin{Theorem}[\cite{Kai05,GN07}]\label{thm:measure and probabilistic Schur}
	Let $p_i \colon H \to G$ be the $i$-th projection map $h \mapsto h|_i$ (where $H = st_G(i)$ and $h|_i$ is the  section of $h$ at vertex $i$ of the first level) and let $\mu_i$ be the image of $\mu_H$ under $p_i$. Then
	\[ \mu_i = k_i(\mu). \]
\end{Theorem}

In the most interesting cases, the group $G$ \ronertwo{acts} level transitively (in particular, transitively on the first level) and its action of the first level $V_1$ is a \ronerthree{free transitive} action of some subgroup $R < S_d$, for instance of $\Z_d$ (the latter always \ronertwo{holds} in the case of binary tree as $\Sy_2 = \Z_2$). In this case for each $i \in X$,  $H = st_G(i) = st_G(1)$ (where $st_G(1)$ is the stabilizer of the first level of $T$), there is a random sequence of hitting times $\tau(n)$ of the subgroup $H$ so that $\sigma^{\tau(n)} = 1$ in \eqref{eq:random walk wreath recursion} and 
\[\psi(Y^{\tau(n)}) = (Y_1^{\tau(n)}, Y_2^{\tau(n)},\hdots, Y_d^{\tau(n)}). \]
Moreover, the random process $Z_i^{(n)} = Y_i^{\tau(n)}$ is a random walk on $G$ determined by the measure $\mu_i$,~$1 \leq i \leq d$. We call $\mu_i$ the section (or the projection) of $\mu$ at vertex $i$.

The maps $k_i \colon \M(G) \to \M(G)$ have the property that they enlarge the \ronerthree{\emph{weight}} $\mu(1)$ of the identity element and hence cannot have  fixed points. \ronertwo{Now let us resolve this \ronerthree{difficulty} by following \cite{Kai05} (see also \cite{GN07}). A measure $\mu$ on \ronerthree{a} self-similar group is said to be self-similar (or self-affine) at position $i, 1 \leq i \leq d$ if for some $\alpha >0$
	\begin{equation}\label{eq:self affine measure}
	\mu_i = (1-\alpha) \delta_e + \alpha \mu,
	\end{equation}
	where $\delta_e$ is the delta mass at the identity element. Observe that $\mu$ is self-similar at position $i$ if and only if it is a fixed point of the map
	\begin{equation}\label{eq:normalize measure with no mass at identity}
	\wt{k}_i \colon \mu \mapsto \frac{k_i(\mu) - k_i(\mu)(e)}{1-k_i(\mu)(e)}.
	\end{equation}}
Thus $\wt{k}_i$ is a modification of $k_i$: we delete from the measure $k_i(\mu)$ the mass at the identity element and normalize. \ronertwo{Note that $\wt{k}_i$ is defined everywhere, except $\delta_e$. We can extend it by assigning $\wt{k}_i(\delta_e) = \delta_e$, which makes $\wt{k}_i$ a continuous map} $\M(G) \to \M(G)$ for the weak topology on $\M(G)$.
We are interested in non-degenerate fixed points because of the following theorem:

\begin{Theorem}\label{thm:non-degenerate measure gives amenability}
	If a self-similar group $G$ has a non-degenerate symmetric self-similar probability measure, then
	\begin{enumerate}[leftmargin=*,labelsep=4.9mm]
		\item \cite{BV05} the rate of escape
		\[ \theta = \lim_{n\to \infty}{\frac{1}{n} Y^{(n)}} = 0, \]
		\item \cite{Kai05} the entropy 
		\[h = \lim_{n\to \infty}{\frac{1}{n} {\sum_{g\in G} {\mu_n(g) \log{\mu_n(g)}}}} = 0 \]
	\end{enumerate}
	(where $\mu_n = \mu * \hdots * \mu$ is the $n$-th convolution of $\mu$ determining  the distribution of random walk at time $n$).
	Hence the group $G$ is amenable in this case.
\end{Theorem}

In Section \ref{sec:probabilistic_Schur} we provide an example of self-similar measure in the case of the group $\Gr$.


\section{Can One Hear the Shape of a Group?}

One of interesting directions of studies in spectral theory of graphs is finding of iso-spectral but not isomorphic graphs. It is inspired by the famous question of M. Kac, ``Can you hear the shape of a drum'' \cite{Kac66}. It attracted a lot of attention of researchers, and after several preliminary results, starting with the result of J. Milnor \cite{Mil64}, the negative answer was given in 1992 by C. Gordon, D.~Webb and S.~Walpert, who constructed a pair of plane regions that have different shapes but identical eigenvalues \cite{GW84}. The regions are concave polygons, their construction uses group theoretical result of T.~Sunada \cite{Sun85}.

In 1993, A. Valete \cite{Val} raised the following question; ``Can one hear the shape of a group?'', which means ``Does the spectrum of the Cayley graph determines it up to isometry?''. The answer is immediate, and is \emph{no} as the spectra of all grids $\Z^d$, $d\geq1$ are the same, namely the interval $[-1,1]$. Still the question has some interest. The paper \cite{DG20} shows that the answer is \emph{no} in a very strong sense.
\begin{Theorem}[\cite{DG20}]
	~\\
	\vspace{-\baselineskip}
	\begin{enumerate}[leftmargin=*,labelsep=4.9mm]
		\item Let  $\mathcal G_{\omega}=\langle S_{\omega}\rangle,  \omega  \in \Omega=\{0,1,2\}^{\mathbb N},  S_{\omega }=\{a,b_{\omega }c_{\omega},d_{\omega }\}$  be a  family  of  groups  of  intermediate  growth  between  polynomial  and  exponential. Then  for  each  $\omega  \in  \Omega$  the  spectrum  of the  Cayley  graph  $\Gamma_{\omega }=\Gamma(\mathcal G_{\omega},S_{\omega})$   is   the  union	
		\[\Sigma= [-\frac{1}{2},0]\cup [\frac{1}{2},1].\]
		\item Moreover,  for  each  $\omega  \in  \Omega$   that  is not eventually  constant  sequence  the  group $\mathcal G_{\omega}$  has  uncountably  many  covering amenable  groups  $\wt G=\langle \wt S \rangle$ (i.e.,  there is a  surjective  homomorphism $\wt G \twoheadrightarrow \mathcal G_{\omega}$)  generated  by $\wt S=\{\wt a, \wt b,\wt c, \wt d\}$ such  that  the  spectrum  of the  Cayley  graphs $\Gamma(\wt G, \wt S)$   is   the  same  set  $\Sigma= [-\frac{1}{2},0]\cup [\frac{1}{2},1]$.
	\end{enumerate}
\end{Theorem}

The proof uses the Hulanicki theorem \cite{Hul64} on characterization of amenable groups in terms of weak containment of \ronerthree{the} trivial representation \ronerthree{in the} regular representation, and a weak Hulanicki type theorem for covering graphs. More examples of this sort are in \cite{GNP20}. The above theorem is for the isotropic case. In the anisotrophic case by the result of D. Lenz, T. Nagnibeda and first author \cite{GLN17aCombinatorics,GLN17bSchreier}, we know only that $\spec(M_P)$ contains a Cantor subset of the Lebesgue measure 0, which is a spectrum of a random Schr\"odinger operator, whose potential is \ronerthree{ruled} by the substitutional dynamical system generated by the substitution $\sigma$ used in presentation \eqref{eq:1st gp presentation}.

In the case of \ronerthree{a} vertex transitive graph, in particular Cayley graph, a natural choice of a spectral measure  is the  spectral measure $\nu$ associated with delta function $\delta_w$ for $w \in V$. The moments of this measure are the probabilities $P_{w,w}^{(n)}$ of return.
\begin{Problem}
	Does the spectral measure $\nu$ determine Cayley graph \ronertwo{of an infinite finitely generated group} up to isometry?
\end{Problem}


\section{Substitutional and Schreier Dynamical System}\label{sec:dynamical system}

Given an alphabet $A = \{a_1, \hdots ,a_m\}$ and a substitution $\rho \colon A \to A^*, \rho(a_i) = A_i(a_\mu)$, assuming that for some distinguished symbol $a \in A$, $a$ is a prefix of $\rho(a)$, we can consider the sequence of iterates
\[ a \arr \rho(a) \arr \rho^2(a) \arr \hdots \arr \rho^n(a) \arr \hdots, \]
where application of $\rho$ to a word $W \in A^*$ means the replacement of each symbol $a_i$ in $A$ by $\rho(a_i)$. If we denote $W_n =\rho^n(a)$, then $W_n$ is a prefix of $W_{n+1}, n = 1,2,\hdots$ and there is a natural limit 
\[W_\infty = \lim_{n\to \infty} W_n.\]
\ronerthree{This limit} $W_\infty$ is an infinite word over $A$ and the words $W_n$ are prefixes of $W_\infty$. Also $W_\infty$ is a fixed point of $\rho$: $\rho(W_\infty) = W_\infty$. Using $W_\infty$, we can now define subshifts of the full shifts $(A^\N, T)$ and $(A^\Z,T)$ (where $T$ is a shift map in the space of sequences). Let us do this for the bilateral shift.

Let $\mathcal{L}(\rho) = \{ W \in A^* \mid W \text{ is a subword of } W_\infty\}$. Equivalently, $\mathcal{L}(\rho)$ consists of words that appear as a subword of some $W_n$ (and hence in all $W_k, k \geq n$).

Now let $\Omega_\rho$ be the set of sequences $\omega = (\omega_n)_{n \in \Z}$ that are unions of words from $\mathcal{L}(\rho)$, where $\omega_n \in A$. In other words, $\omega \in \Omega_\rho$ if and only if for all $m <0$, $n>0$ there exist $M<m$, $N>n$ such that the subword $\omega_M \hdots \omega_N$ of $\omega$ belongs to \supun{$\mathcal{L}(\rho)$}. Obviously, $\Omega_\rho$ is shift invariant closed subset of $A^\Z$. The dynamical system $(\Omega_\rho, T)$ with the shift map $T$ restricted to $\Omega_\rho$ is a substitutional dynamical system generated by $\rho$.

The most important case is when such system is minimal, i.e., for each $x \in \Omega_\rho$ the orbit $\{ T^nx \}_{n=-\infty}^\infty$ is dense in $\Omega_\rho$. \ronertwo{For instance, this} is the case when the substitution $\rho$ is primitive, which means that there exists $K$ such that for each $i,j, 1 \leq i, j \leq m$ \ronertwo{the} symbol $a_i$ occur in the word $\rho^K(a_j)$.

By Krylov-Bogolyubov theorem, the system $(\Omega_\tau, T)$ has at least one $T$-invariant probability measure and the invariant  ergodic measures (i.e., extreme points of the simplex of $T$-invariant probability measures) are of special interest. Another important case is when the system $(\Omega_\rho, T)$ is uniquely ergodic, i.e., there is only one invariant probability measure (necessarily ergodic).

A subshift $(\Omega_\rho,T)$ is called linearly repetitive (LR) if there exists a constant $C$ such that any word $W \in \mathcal{L}(\rho)$ occurs in any word $U \in  \mathcal{L}(\rho)$ of length $\geq C|W|$. This is a stronger condition than minimality. The following result \ronertwo{goes back to M. Boshernitzan \cite{Bos84}} (see also \cite{Dur00}).

\begin{Theorem}\label{thm:linearly_repetitive}
	Let $(\Omega,T)$ be a linearly repetitive subshift. Then, the subshift is uniquely ergodic.
\end{Theorem}
Here, it is not necessary for the subshift to be generated by a substitution.
It is known that subshifts associated with primitive substitutions are linearly repetitive \cite{DZ00,DHS99}.
Theorem 1 of \cite{DL06} shows that linear repetitivity in fact holds for subshifts associated to any substitution provided minimality holds. Unique ergodicity is then a direct consequence of linear repetitivity due to Theorem \ref{thm:linearly_repetitive}.

The classical example of \ronertwo{a} substitutional system is the Thue-Morse system determined by the substitution $0 \arr01, 1 \arr 10$ over binary alphabet \cite{BLRS09}.

Following \cite{GLN17aCombinatorics,GLN17bSchreier,GLN18} we consider the substitution $\sigma \colon a \arr aca, b \arr d, c \arr b, d \arr c$ over alphabet $\{a,b,c,d\}$ and system $(\Omega_\sigma,T)$ generated by it. Despite \ronerthree{$\sigma$ not being} primitive, the system $(\Omega_\sigma,T)$ satisfies the linear repetivity property (in fact the same system can be generated by a primitive substitution $\sigma' \colon a \arr ac, b \arr ac, c \arr ad, d \arr ab$).

An additional property of the fixed point $\eta = \lim_{n \to \infty} \sigma^n(a)$ is that it is a Toeplitz sequence. i.e., for each entry $\eta_n$ of $\eta = (\eta_n)_{n =0}^\infty$ there is period $p = p(n)$ such that all entries with indices of the form $n+pk, k = 0,1, \hdots$ contain the same symbol $\eta_n$. In our case the periods have the form $2^l, l \in \N$. More on combinatorial properties of $\eta$ and associated system see \cite{GLN17aCombinatorics}.

Our interest \ronertwo{in} \ronerthree{the} substitution $\sigma$ and \ronerthree{the} associated subshift \ronertwo{comes} from \ronertwo{following four} facts:
\begin{enumerate}[leftmargin=*,labelsep=4.9mm]
	\item \ronerthree{The substitution} $\sigma$ \ronertwo{appears in} the presentation \eqref{eq:1st gp presentation} of the group $\Gr = \gen{a,b,c,d}$ of intermediate growth by generators and relations, so it determines the group modulo finite set of \ronertwo{relators}.
	\item The system $(\Omega_\sigma,T)$ gives a model for a Schreier dynamical system (in terminology of \cite{Gri11}) determined by the action of $\Gr$ on the boundary $\partial T$ of binary tree.
	\item The latter property allows \ronerthree{us} to translate the spectral properties of Schreier graphs $\Gamma_x, x \in \partial T$ into the spectral properties of the corresponding random \ronertwo{Schr\"odinger} operator and conclude that in anisotropic case the spectrum is a Cantor set of Lebesgue measure zero \cite{GLN17bSchreier,GLN18}.
	\item The group $\Gr$ embeds into \ronerthree{the} topological full group $[[\sigma]]$ associated with the subshift $(\Omega_\sigma,T)$.
\end{enumerate}

For any minimal action $\alpha$ of a group $G$ on a Cantor set $X$, one can define a topological full group (TFG in short) $[[\alpha]]$ as a group consisting of homeomorphisms $h \in Homeo(X)$ that locally act as elements of $G$. \ronerthree{If $G$ is the infinite cyclic group generated by a minimal homeomorphism of a Cantor set, the TFG is an invariant of the Cantor minimal system up to the flip conjugacy \cite{GPS99} and its commutator $[[\alpha]]'$ is a simple group. Moreover, $[[\alpha]]'$ is finitely generated if the system is conjugate to a minimal subshift over a finite alphabet.} It was conjectured by K.~Medynets and the first author, and proved by K.~Juschenko and N.~Monod \cite{JM13} that if $G = \Z$, then $[[\alpha]]$ is amenable. Thus \ronertwo{TFGs} are a rich source of non-elementary amenable groups, and satisfy many unusual properties \cite{GM14,Mat06}. N.~Matte~Bon observed that $\Gr$ embeds into $[[\sigma]]$, where $\sigma$ is the substitution from \eqref{eq:1st gp presentation} \cite{Mat15}. A similar result holds for overgroup $\wt\Gr$.

Study of substitutional dynamical systems and more generally of aperiodic order is a rich area of mathematics (see \cite{KLS15,Baa17} and references there for instance). A special attention is paid to the classical substitutions like Thue-Morse, Arshon \cite{Ars37}, and Rudin-Shapiro substitutions.

\begin{Problem}
	For which primitive substitutions $\tau$, the TFG $[[\tau]]$, contains a subgroup of intermediate growth? contains a subgroup of Burnside type (i.e., finitely generated infinite torsion group)? In particular, does the classical substitutions listed above have such properties?
\end{Problem}

Given a Schreier graph $\Gamma = \Gamma(G,H,A) \in \Sch{m}$, one can consider the action of $G$ on ${\{ (\Gamma,v) \mid v \in V \}}$ (i.e., on the set of marked graphs where $(\Gamma,v) \xrightarrow{g} (\Gamma,gv)$) and extend it to the action on the closure $\overline{\{ (\Gamma,v) \mid v \in V \}}$ in $\Sch{m}$. This is called in \cite{Gri11} a Schreier dynamical system. Study of such systems is closely related to the study of invariant random subgroups. In important cases, such systems \ronerthree{allows} to \ronerthree{recover} the original action $(G,X)$ if $\Gamma = \Gamma_x$, $x \in X$ is an orbital graph. In particular, this holds if the action is extremely non-free (i.e., stabilizers $G_x$ of different points $x \in X$ are distinct). The action of $\Gr$ and any group of branch type is extremely \ronertwo{non-}free. More on this is in \cite{Gri11}.


\section{Computation of Schur Maps for $\Gr$ and $\tilde{\Gr}$}\label{sec:schur_calculation}

Recall the matrix recursions between generators of $\Gr = \gen{a,b,c,d}$ (see \eqref{eq:matrix map 1st gp}),
\begin{align*}
	1 & = 
	\left(
	\begin{array}{cc}
		1 & 0 \\
		0 & 1 \\
	\end{array}
	\right), &
	a & = 
	\left(
	\begin{array}{cc}
		0 & 1 \\
		1 & 0 \\
	\end{array}
	\right), &
	b & = 
	\left(
	\begin{array}{cc}
		a & 0 \\
		0 & c \\
	\end{array}
	\right), &
	c & = 
	\left(
	\begin{array}{cc}
		a & 0 \\
		0 & d \\
	\end{array}
	\right), &
	d & = 
	\left(
	\begin{array}{cc}
		1 & 0 \\
		0 & b \\
	\end{array}
	\right).
\end{align*}
Let $M = xa+yb+zc+ud+v1$ be an element of the group algebra $\CC[\Gr]$. By using \eqref{eq:matrix map 1st gp}, we identify,
\begin{equation} \label{eq:first gp operator}
M = \left(
\begin{array}{cc}
(y+z)a+(u+v)1 & x \\
x & ub+yc+zd+v1 \\
\end{array}
\right).
\end{equation}

First we will calculate the first Schur complement $S_1(M)$, which is defined when $D = v1+ub+yc+zd$ is invertible. Since the group generated by $\{1,b,c,d\}$ is isomorphic to $\mathbb{Z}_2^2$ (via the identification $1,b,c,d$ with $(0,0), (1,0), (0,1), (1,1)$, respectively), by a direct calculation, we obtain that $D$ is invertible if and only if 
\begin{equation}\label{eq:1st Schur 1st gp denominator}
(v+u+y+z)(v-u+y-z)(v+u-y-z)(v-u-y+z) \neq 0,
\end{equation}
and if the condition in \eqref{eq:1st Schur 1st gp denominator} is satisfied, then $D^{-1}$ is given by,
\begin{align*}\label{eq:1st Schur 1st gp inverse}
	D^{-1} = & \frac{1}{4} \left( \frac{1}{v+u+y+z} + \frac{1}{v-u+y-z} + \frac{1}{v+u-y-z} + \frac{1}{v-u-y+z} \right) 1 \\
	& + \frac{1}{4} \left( \frac{1}{(v+u+y+z)} - \frac{1}{v-u+y-z} + \frac{1}{v+u-y-z} - \frac{1}{v-u-y+z} \right) b \\
	& + \frac{1}{4} \left( \frac{1}{(v+u+y+z)} + \frac{1}{v-u+y-z} - \frac{1}{v+u-y-z} - \frac{1}{v-u-y+z} \right) c \\
	& + \frac{1}{4} \left( \frac{1}{(v+u+y+z)} - \frac{1}{v-u+y-z} - \frac{1}{v+u-y-z} + \frac{1}{v-u-y+z} \right) d.
\end{align*}
Therefore,
\begin{align*}
	S_1(M) = & A - B D^{-1} C \\
	= & (y+z)a + (v+u)1 - x^2 D^{-1}\\
	= & (y+z)a \\
	& + \left( v+u -  x^2  \frac{2uyz - v(-v^2+u^2+y^2+z^2)}{(v+u+y+z)(v-u+y-z)(v+u-y-z)(v-u-y+z)}   \right) 1 \\
	& - x^2 \frac{2vyz - u(v^2-u^2+y^2+z^2)}{(v+u+y+z)(v-u+y-z)(v+u-y-z)(v-u-y+z)} b, \\
	& - x^2 \frac{2vuz - y(v^2+u^2-y^2+z^2)}{(v+u+y+z)(v-u+y-z)(v+u-y-z)(v-u-y+z)} c, \\
	& - x^2 \frac{2vuy - z(v^2+u^2+y^2-z^2)}{(v+u+y+z)(v-u+y-z)(v+u-y-z)(v-u-y+z)} d.
\end{align*}
This leads to the map $\wt {G} \colon \CC^5 \to \CC^5$ given in \eqref{eq:5dim G map 1stGri}.

Now we will calculate the second Schur complement $S_2(M)$ which is defined when $A = (y+z)a + (u+v)1$ is invertible. Since the group generated by $\{1,a\}$ is isomorphic to $\mathbb{Z}_2$ (via the identification $1,a$ with $0,1$, respectively), by a direct calculation, we obtain that $A$ is invertible if and only if 
\begin{equation}\label{eq:2nd Schur 1st gp denominator}
(v+u+y+z)(v+u-y-z) \neq 0,
\end{equation}
and if the condition in \eqref{eq:2nd Schur 1st gp denominator} is satisfied, then $A^{-1}$ is given by,
\begin{align*}
	A^{-1}  & =  \frac{1}{2} \left( \frac{1}{v+u+y+z} + \frac{1}{v+u-y-z} \right) 1  
	+
	\frac{1}{2} \left( \frac{1}{v+u+y+z} - \frac{1}{v+u-y-z} \right) a \\
	& =   \frac{v+u}{(v+u+y+z)(v+u-y-z)} 1 - \frac{y+z}{(v+u+y+z)(v+u-y-z)} a.
\end{align*}
Therefore,
\begin{align*}
S_2(M) = & D - C A^{-1} B \\
= & v1+ub+yc+zd - x^2 A^{-1}\\
= & \frac{x^2(y+z)}{(v+u+y+z)(v+u-y-z)} a +ub+yc+zd +\left(v - \frac{x^2(v+u)}{(v+u+y+z)(v+u-y-z)} \right)1.
\end{align*}

This leads to the map $\wt {F} \colon \CC^5 \to \CC^5$ given in \eqref{eq:5dim F map 1stGri}.

Now consider the case where $y=z=u=1$. Note that $\wt{F}$ fixes second, third and fourth coordinates and so we may restrict the map to first and fifth coordinates. Therefore we get $\CC^2 \to \CC^2$ map

\begin{equation*}
	\wt{F} \colon 
	\left(
	\begin{array}{c}
		x \\
		v \\
	\end{array}
	\right)
	\mapsto 
	\left(
	\begin{array}{c}
		\dfrac{2 x^2}{(v+3)(v-1)} \\[3mm]
		v -  \dfrac{x^2(v+1)}{(v+3)(v-1)}  \\
	\end{array}
	\right).
\end{equation*}

By the change of coordinates $(x,v) \to (-x,-1-y)$, we obtain $F$ given in \eqref{eq:F map 1stGri}.

Now note that second, third and fourth coordinates of $\wt{G}$ are the same and are equal to $\dfrac{x^2}{(v+3)(v-1)}$. By re-normalization (i.e., multiplying by $\frac{(v+3)(v-1)}{x^2}$) we obtain a map which fixes second, third and fourth coordinates. So we may restrict the map to first and fifth coordinates and get $\CC^2 \to \CC^2$ map
\begin{equation*}
	\wt{G} \colon 
	\left(
	\begin{array}{c}
		x \\
		v \\
	\end{array}
	\right)
	\mapsto 
	\left(
	\begin{array}{c}
		\dfrac{2(v+3)(v-1)}{x^2} \\[3mm]
		- 2 - v + (v+1)\dfrac{(v+3)(v-1)}{x^2}  \\
	\end{array}
	\right).
\end{equation*}
By the change of coordinates $(x,v) \to (-x,-1-y)$, we obtain $G$ given in \eqref{eq:G map 1stGri}.

Now consider the overgroup $\wt{\mathcal{G}} = \gen{a,\wt {b},\wt {c}, \wt {d}} \leq  \Aut(T_2)$, where $\wt {b},\wt {c}, \wt {d}$ satisfy matrix recursions given by \eqref{eq:matrix map overgp} and $a$ is a generator of $\Gr$. We have $b= \wt c \wt d$, $c = \wt b \wt d$, $d = \wt b \wt c$ and hence $\Gr$ is a subgroup of $\wt\Gr$. It will be convenient to consider $\wt\Gr$ as a group generated by eight elements $a,b,c,d,\wt {a},\wt {b},\wt {c}, \wt {d}$, where $\wt a$ satisfies the matrix recursion in \eqref{eq:matrix map overgp}.
\begin{align}\label{eq:matrix map overgp}
	\wt {a} & = 
	\left(
	\begin{array}{cc}
		a & 0 \\
		0 & \wt {a} \\
	\end{array}
	\right), &
	\wt {b} & = 
	\left(
	\begin{array}{cc}
		1 & 0 \\
		0 & \wt {c} \\
	\end{array}
	\right), &
	\wt {c} & = 
	\left(
	\begin{array}{cc}
		1 & 0 \\
		0 & \wt {d} \\
	\end{array}
	\right), &
	\wt {d} & = 
	\left(
	\begin{array}{cc}
		a & 0 \\
		0 & \wt {b} \\
	\end{array}
	\right).
\end{align}
$\mathcal{G}$ is a subgroup of $\wt {\mathcal{G}}$ and so $\mathbb{C}[\mathcal{G}]$ is a subalgebra of $\mathbb{C}[\wt {\mathcal{G}}]$. So we can use \eqref{eq:matrix map 1st gp} as the matrix recursions of $1,a,b,c,d$. 
Let $M = xa+yb+zc+ud+q\wt {a}+r\wt {b}+s\wt {c}+t\wt {d}+v1$. By using \eqref{eq:matrix map overgp} and \eqref{eq:matrix map 1st gp}, we obtain the matrix recursion of $M$ as,

\begin{equation} \label{eq:overgp operator}
M = \left(
\begin{array}{cc}
(y+z+q+t)a+(u+r+s+v)1 & x \\
x & ub+yc+zd+q\wt {a}+t\wt {b}+r\wt {c}+s\wt {d}+v1 \\
\end{array}
\right).
\end{equation}

We can calculate first and second Schur complements $S_1(M), S_2(M)$, and the multi-dimensional maps $\wt S_1, \wt{S}_2$ associated with them as we did for the case of $\Gr$. These maps are nine-dimensional and are given by 
\begin{equation*}\label{eq:1st Schur overgp transformation}
\wt{S}_1 \colon 
\left(
\begin{array}{c}
	x \\
	y \\
	z \\
	u \\
	q \\
	r \\
	s \\
	t \\
	v \\
\end{array}
\right)
\mapsto 
\left(
\begin{array}{c}
	y+z+q+t \\
	- \frac{x^2}{8} \left(\frac{1}{\hat{D}_{000}} - \frac{1}{\hat{D}_{100}} + \frac{1}{\hat{D}_{010}} + \frac{1}{\hat{D}_{001}} - \frac{1}{\hat{D}_{110}} - \frac{1}{\hat{D}_{101}} + \frac{1}{\hat{D}_{011}} - \frac{1}{\hat{D}_{111}} \right) \\
	- \frac{x^2}{8} \left(\frac{1}{\hat{D}_{000}} + \frac{1}{\hat{D}_{100}} - \frac{1}{\hat{D}_{010}} + \frac{1}{\hat{D}_{001}} - \frac{1}{\hat{D}_{110}} + \frac{1}{\hat{D}_{101}} - \frac{1}{\hat{D}_{011}} - \frac{1}{\hat{D}_{111}} \right) \\
	- \frac{x^2}{8} \left(\frac{1}{\hat{D}_{000}} - \frac{1}{\hat{D}_{100}} - \frac{1}{\hat{D}_{010}} + \frac{1}{\hat{D}_{001}} + \frac{1}{\hat{D}_{110}} - \frac{1}{\hat{D}_{101}} - \frac{1}{\hat{D}_{011}} + \frac{1}{\hat{D}_{111}} \right) \\
	- \frac{x^2}{8} \left(\frac{1}{\hat{D}_{000}} - \frac{1}{\hat{D}_{100}} - \frac{1}{\hat{D}_{010}} - \frac{1}{\hat{D}_{001}} + \frac{1}{\hat{D}_{110}} + \frac{1}{\hat{D}_{101}} + \frac{1}{\hat{D}_{011}} - \frac{1}{\hat{D}_{111}} \right) \\
	- \frac{x^2}{8} \left(\frac{1}{\hat{D}_{000}} + \frac{1}{\hat{D}_{100}} - \frac{1}{\hat{D}_{010}} - \frac{1}{\hat{D}_{001}} - \frac{1}{\hat{D}_{110}} - \frac{1}{\hat{D}_{101}} + \frac{1}{\hat{D}_{011}} + \frac{1}{\hat{D}_{111}} \right) \\
	- \frac{x^2}{8} \left(\frac{1}{\hat{D}_{000}} - \frac{1}{\hat{D}_{100}} + \frac{1}{\hat{D}_{010}} - \frac{1}{\hat{D}_{001}} - \frac{1}{\hat{D}_{110}} + \frac{1}{\hat{D}_{101}} - \frac{1}{\hat{D}_{011}} + \frac{1}{\hat{D}_{111}} \right) \\
	- \frac{x^2}{8} \left(\frac{1}{\hat{D}_{000}} + \frac{1}{\hat{D}_{100}} + \frac{1}{\hat{D}_{010}} - \frac{1}{\hat{D}_{001}} + \frac{1}{\hat{D}_{110}} - \frac{1}{\hat{D}_{101}} - \frac{1}{\hat{D}_{011}} - \frac{1}{\hat{D}_{111}} \right) \\
	\left((u+r+s+v) - \frac{x^2}{8} \left(\frac{1}{\hat{D}_{000}} + \frac{1}{\hat{D}_{100}} + \frac{1}{\hat{D}_{010}} + \frac{1}{\hat{D}_{001}} + \frac{1}{\hat{D}_{110}} + \frac{1}{\hat{D}_{101}} + \frac{1}{\hat{D}_{011}} + \frac{1}{\hat{D}_{111}} \right)\right)  \\
\end{array}
\right),
\end{equation*}
\begin{equation}\label{eq:2nd Schur overgp transformation}
\wt S_2 \colon 
\left(
\begin{array}{c}
x \\
y \\
z \\
u \\
q \\
r \\
s \\
t \\
v \\
\end{array}
\right)
\mapsto 
\left(
\begin{array}{c}
\frac{x^2(y+z+q+t)}{ (v+u+r+s+y+z+q+t)(v+u+r+s-y-z-q-t)} \\
u \\
y \\
z \\
q \\
t \\
r \\
s \\
v -  \frac{x^2(v+u+r+s)}{ (v+u+r+s+y+z+q+t)(v+u+r+s-y-z-q-t)}  \\
\end{array}
\right),
\end{equation}
where
\begin{align}\label{eq:1st Schur overgp denominator}
	\hat{D}_{000} & = v + u + y + s + z + r + t + q, \nonumber\\
	\hat{D}_{100} & = v - u + y + s - z - r + t - q, \nonumber\\
	\hat{D}_{010} & = v + u - y + s - z + r - t - q, \nonumber\\
	\hat{D}_{001} & = v + u + y - s + z - r - t - q, \nonumber\\
	\hat{D}_{110} & = v - u - y + s + z - r - t + q, \nonumber\\
	\hat{D}_{101} & = v - u + y - s - z + r - t + q, \nonumber\\
	\hat{D}_{011} & = v + u - y - s - z - r + t + q, \nonumber\\
	\hat{D}_{111} & = v - u - y - s + z + r + t - q.
\end{align}

\section{Probabilistic Schur Map for $\Gr$}\label{sec:probabilistic_Schur}

Recall  that  maps  $k_i(\mu)$  were  defined  by  \eqref{eq:normalize measure with no mass at identity}. Let $M = xa +yb+zc+ud \in \CC[\Gr]$ be the Markov operator of random walk determined by the measure $\mu = x\mass{a} +y\mass{b}+z\mass{c}+u\mass{d}$, supported on the generating set of $\Gr$ (i.e., $x,y,z,u > 0$ and $x+y+z+u=1$). By \eqref{eq:first gp operator}, \ronerthree{taking} $v=0$ gives the matrix recursion,
\begin{equation}
M = \left(
\begin{array}{cc}
(y+z)a+(u)1 & x \\
x & ub+yc+zd \\
\end{array}
\right).
\end{equation}
Recall that $K_i = \J S_i \J\inv$ (where $\J$ is the shift map and $S_i$ is the $i$-th Schur complement) and $k_i$ is the \ronerthree{analogous} map on the simplex of probability measures $\M(\Gr)$ (see Theorem \ref{thm:measure and probabilistic Schur}). Then,
\[k_i(\mu) = X_i \mass{a} + Y_i \mass{b} +  Z_i \mass{c} + U_i \mass{d} + V_i \mass{e}\]
for $i =1,2$, where
\begin{align}\label{Prob_Schur_Coordinates}
	X_1 & = y+z, \\
	V_1 & = u -  \frac{x^2(2uyz +u^2+y^2+z^2-1)}{(u+y+z-1)(-u+y-z-1)(u-y-z-1)(-u-y+z-1)}, \nonumber\\
	Y_1 & = \frac{x^2(2yz + u(1-u^2+y^2+z^2))}{(u+y+z-1)(-u+y-z-1)(u-y-z-1)(-u-y+z-1)}, \nonumber\\
	Z_1 & = \frac{x^2(2uz + y(1+u^2-y^2+z^2))}{(u+y+z-1)(-u+y-z-1)(u-y-z-1)(-u-y+z-1)}, \nonumber\\
	U_1 & = \frac{x^2(2uy + z(1+u^2+y^2-z^2))}{(u+y+z-1)(-u+y-z-1)(u-y-z-1)(-u-y+z-1)}, \nonumber
\end{align}
and
\begin{align}
	X_2 & = \frac{x^2(y+z)}{(1-u-y-z)(1-u+y+z)}, \nonumber\\
	V_2 & = \frac{x^2(1-u)}{(1-u-y-z)(1-u+y+z)}, \nonumber\\
	Y_2 & = u, \nonumber\\
	Z_2 & = y, \nonumber\\
	U_2 & = z. \nonumber
\end{align}

We are interested in self-similar measures (i.e., measures that satisfy \eqref{eq:self affine measure} or, which is the same, fixed points of the map \eqref{eq:normalize measure with no mass at identity}). A direct calculation shows that $\wt{k}_2$ defined by \eqref{eq:normalize measure with no mass at identity} has no fixed points and so $\mu$ is not self-similar at the second coordinate. Therefore we restrict the rest of the section to study the map $\wt{k}_1$.

In order to understand $\wt{k}_1$, we extend it to the map,
\begin{align*}
	\widehat{k}_1 \colon \Delta & \to \Delta \\
	(x,y,z,u) & \mapsto \left( \frac{X_1}{1-V_1}, \frac{Y_1}{1-V_1}, \frac{Z_1}{1-V_1}, \frac{U_1}{1-V_1}\right),
\end{align*}
where $\Delta$ is the 3-simplex $\{(x,y,z,u) \mid x+y+z+u=1, x,y,z,u \geq 0 \}$. Note that the vertices (three coordinates are 0), the edges (two coordinates are 0), and the face $x=0$ correspond to degenerate probability measures $\mu$, whereas the faces $y=0$, $z=0$, and $u=0$ correspond to non-degenerate measures in $\M(\Gr)$. This is due to the fact that the group $\Gr$ is in fact 3-generated and removing exactly one of the elements $b,c$ or $d$ form the set $\{a,b,c,d\}$, still generates $\Gr$.  A direct calculation yields the following proposition.

\begin{figure}[H]
\centering
\begin{subfigure}[b]{0.2\textwidth}
	\includegraphics[width=\textwidth,trim={0 0 0 11},clip]{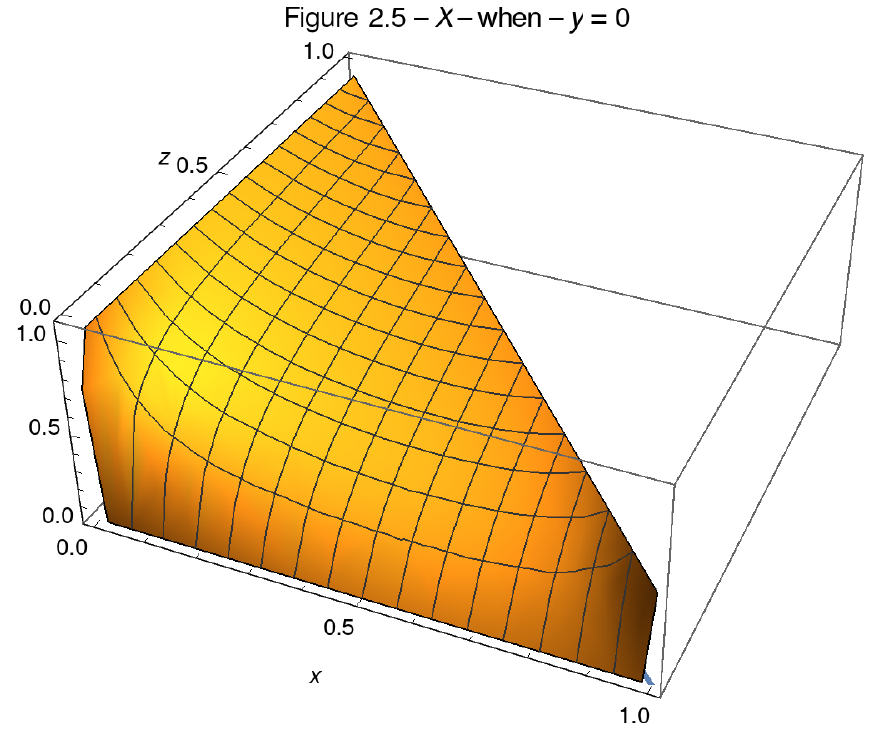}
	\caption{}
	\label{subfig:simplex_X}
\end{subfigure}
\qquad
\begin{subfigure}[b]{0.2\textwidth}
	\includegraphics[width=\textwidth,trim={0 0 0 11},clip]{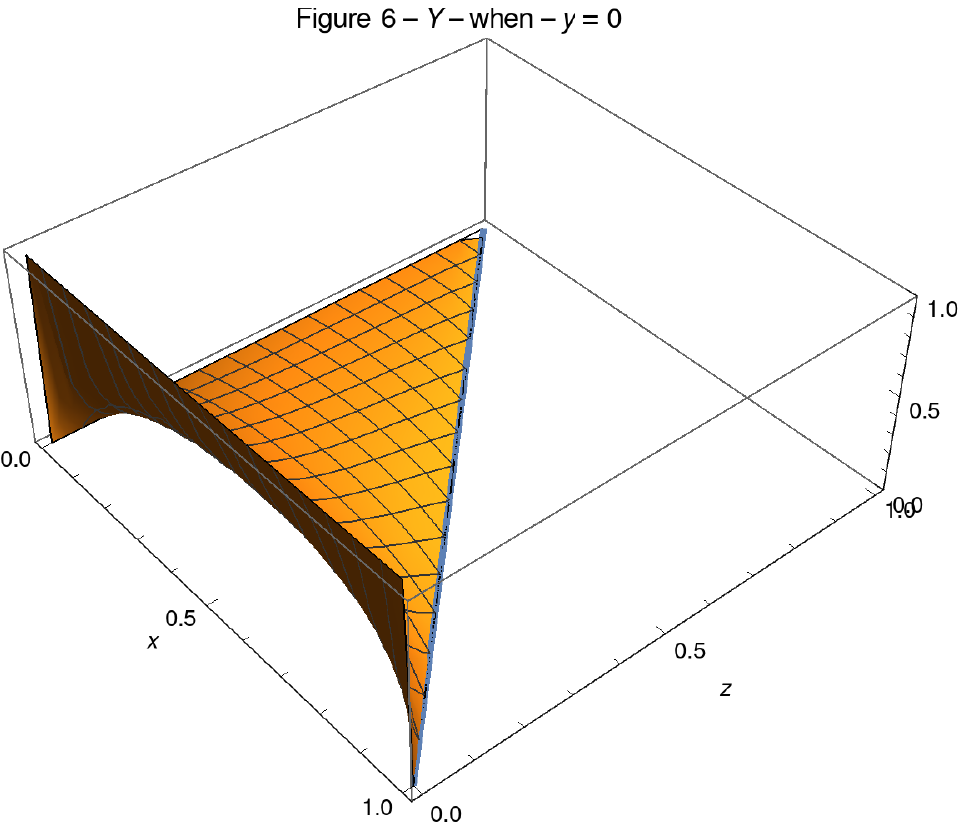}
	\caption{}
	\label{subfig:simplex_Y}
\end{subfigure}
\qquad
\begin{subfigure}[b]{0.2\textwidth}
	\includegraphics[width=\textwidth,trim={0 0 0 11},clip]{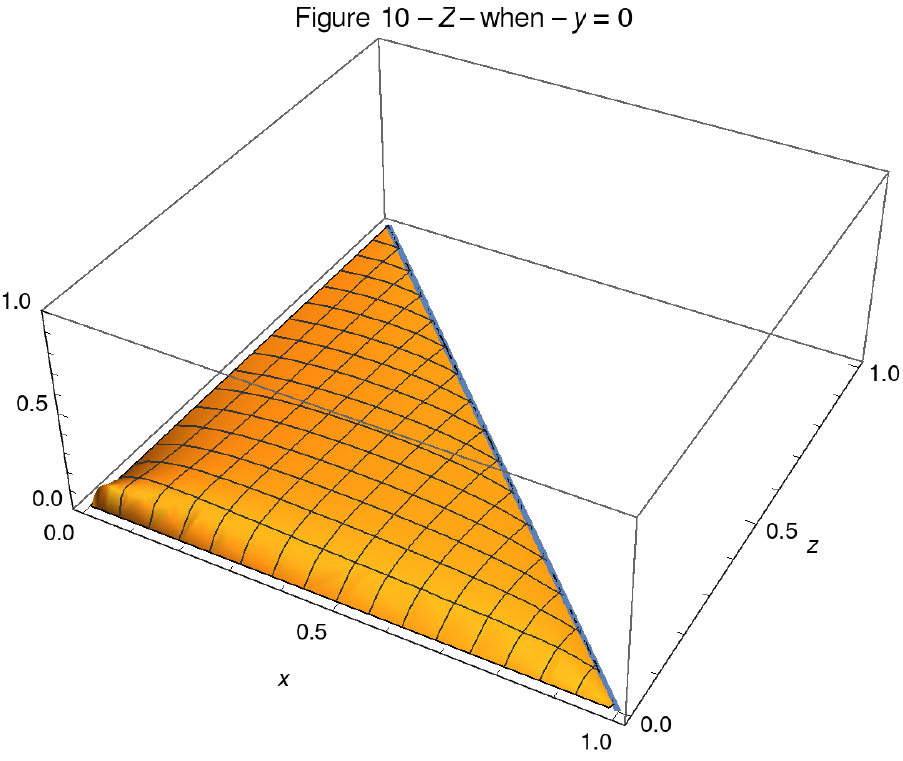}
	\caption{}
	\label{subfig:simplex_Z}
\end{subfigure}
\qquad
\begin{subfigure}[b]{0.2\textwidth}
	\includegraphics[width=\textwidth,trim={0 0 0 11},clip]{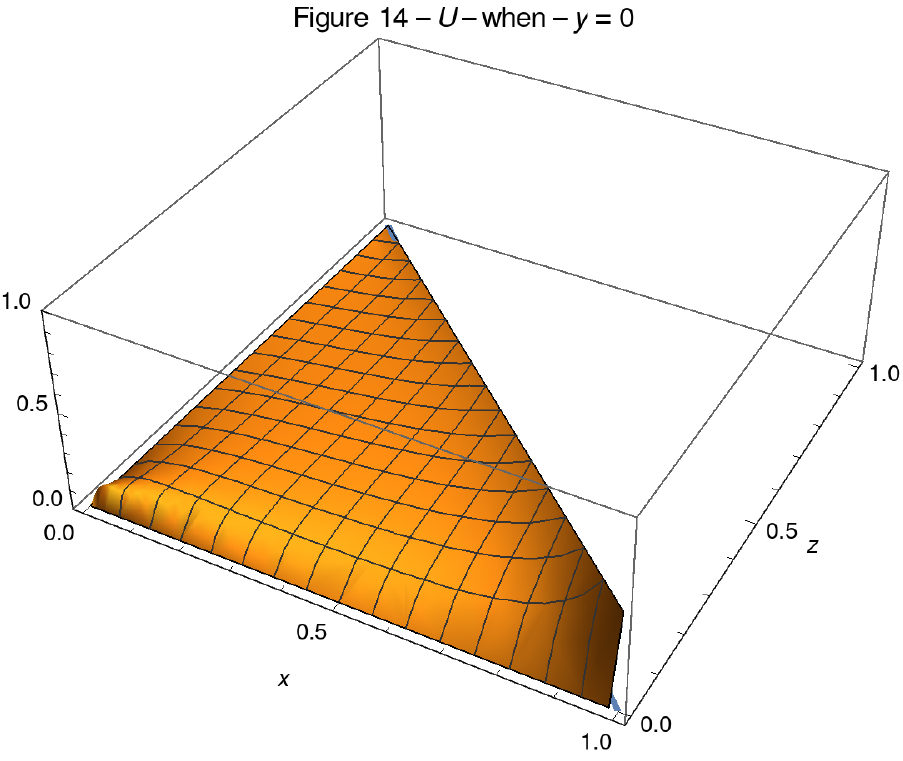}
	\caption{}
	\label{subfig:simplex_U}
\end{subfigure}
\caption{$\wt{k}_1$ values on the face $y=0$, where (\textbf{a}) $X$ values, (\textbf{b}) $Y$ values, (\textbf{c}) $Z$ values, and (\textbf{d}) $U$ values, plotted in $xz$ plane.}\label{fig:simplex}
\end{figure}

\begin{Proposition}
	Consider the map $\widehat{k}_1$ given above. Then;
	\begin{enumerate}[leftmargin=*,labelsep=4.9mm]
		\item All vertices are indeterminacy points.
		\item The edge $(x,y,0,0)$ maps to the edge $(X,0,Z,0)$, the edge $(x,0,z,0)$ maps to the edge $(X,0,0,U)$, and the edge $(x,0,0,0)$ maps to the vertex $(0,1,0,0)$.
		\item The face $x=0$ maps to the vertex $(1,0,0,0)$ and the faces $y=0, z=0$ and $u=0$ map to the interior of the 3-simplex $\Delta$. (See Figure \ref{fig:simplex}.)
		\item Interior of $\Delta$ maps to itself.
	\end{enumerate}
\end{Proposition}

$\wt{k}_1$ has a fixed point $\left(4/7, 1/7, 1/7, 1/7\right) \in \Delta$ and thus
\[\wt{k}_1(\mu) = \left(1-\frac{1}{2}\right)\mass{e} + \frac{1}{2} \mu. \]

\begin{Problem}
	Describe all fixed points of the maps $\wt{k}_1, \wt{k_2} \colon \M(\Gr) \to \M(\Gr)$.
\end{Problem}

\section{Random Groups $\{\Gr_\omega\}$ and Associated 4--Parametric Family of Maps}

Here we will introduce a family of subgroups of $\Aut T_2$, $\{ \Gr\om \mid \omega \in \Omega\}$, where $\Omega = \{0,1,2\}^{\N}$. Since each element in $\Aut T_2$ can be defined by wreath recursions \eqref{eq:wreath recursions}, for $\omega = \omega_0 \omega_1 \hdots \in \Omega$, we define recursively,
\begin{align}\label{eq:random recursions}
	b\om & = (b\omz,b\som), & c\om & = (c\omz,c\som), & d\om & = (d\omz,d\som), 
\end{align}
where $T$ is the left shift operator on $\Omega$, $b_0 = b_1 = c_0 = c_2 = d_1 = d_2 =a$, and $b_2 = c_1 = d_0 =1$. Here $a = (1,1)\sigma$ where $\sigma$ is the permutation of the symmetric group $\Sy_2$. Now define
\begin{equation}\label{eq:gen_gri_gp}
\Gr_\omega = \gen{a, b_\omega, c_\omega, d_\omega}.
\end{equation}
By \eqref{eq:random recursions} we obtain the recursions for the Koopman representation
\begin{align*}
	a & = 
	\left(
	\begin{array}{cc}
		0 & 1 \\
		1 & 0 \\
	\end{array}
	\right), &
	b\om & = 
	\left(
	\begin{array}{cc}
		b_{\omega_0} & 0 \\
		0 & b\som \\
	\end{array}
	\right), &
	c\om & = 
	\left(
	\begin{array}{cc}
		c_{\omega_0} & 0 \\
		0 & c\som \\
	\end{array}
	\right), &
	d\om & = 
	\left(
	\begin{array}{cc}
		d_{\omega_0} & 0 \\
		0 & d\som \\
	\end{array}
	\right).
\end{align*}
For ${M = xa+yb\om+zc\om+ud\om+v1}$, an element of the group algebra $\CC[\Gr\om]$, using above matrix recursions, we identify,
\begin{equation}\label{eq:matrix_gen_gri}
M = \left(
\begin{array}{cc}
y b\omz + z c\omz + u d\omz +v1 & x \\
x & yb\som+zc\som+ud\som+v1 \\
\end{array}
\right).
\end{equation}
We are interested in calculating \ronerthree{the Schur maps associated with $M$}. Direct calculation shows that the second Schur complement 
\[S_2(M) =
\begin{cases}
\frac{x^2(y+z)}{(v+u+y+z)(v+u-y-z)} a +yb\som+zc\som+ud\som +\left(v - \frac{x^2(v+u)}{(v+u+y+z)(v+u-y-z)} \right)1 & ;\omega_0 = 0\\
\frac{x^2(y+u)}{(v+z+u+y)(v+z-u-y)} a +yb\som+zc\som+ud\som +\left(v -  \frac{x^2(v+z)}{(v+z+u+y)(v+z-u-y)} \right)1 & ;\omega_0 = 1\\
\frac{x^2(z+u)}{(v+y+z+u)(v+y-z-u)} a +yb\som+zc\som+ud\som +\left(v -  \frac{x^2(v+y)}{(v+y+z+u)(v+y-z-u)} \right)1 & ;\omega_0 = 2
\end{cases}.
\]

Note that the middle three coefficients under $S_2$ are fixed and so independent of $\omega$ (or $\omega_0$). This allows us to reduce the second Schur map into two dimensional maps (i.e., $\CC \to \CC$) on three parameters $y,z,u$ and symbols $0,1,2$;
\begin{equation*}
	F_0 \colon 
	\left(
	\begin{array}{c}
		x \\
		v \\
	\end{array}
	\right)
	\mapsto 
	\left(
	\begin{array}{c}
		\dfrac{x^2(y+z)}{(v+u+y+z)(v+u-y-z)} \\
		v -  \dfrac{x^2(v+u)}{(v+u+y+z)(v+u-y-z)}  \\
	\end{array}
	\right),
\end{equation*}
\begin{equation*}
	F_1 \colon 
	\left(
	\begin{array}{c}
		x \\
		v \\
	\end{array}
	\right)
	\mapsto 
	\left(
	\begin{array}{c}
		\dfrac{x^2(y+u)}{(v+z+u+y)(v+z-u-y)} \\
		v -  \dfrac{x^2(v+z)}{(v+z+u+y)(v+z-u-y)}  \\
	\end{array}
	\right),
\end{equation*}
\begin{equation}\label{eq:random map}
F_2 \colon 
\left(
\begin{array}{c}
x \\
v \\
\end{array}
\right)
\mapsto 
\left(
\begin{array}{c}
\dfrac{x^2(z+u)}{(v+y+z+u)(v+y-z-u)} \\
v -  \dfrac{x^2(v+y)}{(v+y+z+u)(v+y-z-u)}  \\
\end{array}
\right).
\end{equation}
Thus for a given $\omega \in \Omega$, applying \ronerthree{the} second Schur complement $n$ times is equivalent to \ronerthree{taking} composition $F_{\omega_{n-1}} \circ \hdots \circ F_{\omega_0}$.

Consider the family of 2-dimensional (i.e., $\CC^2 \to \CC^2$) 4--parametric maps $\{{F}_{(\alpha,\beta,\gamma,\delta)} \mid \alpha,\beta,\gamma,\delta \in \CC \}$ given by
\begin{equation}\label{eq:4para map}
{F}_{(\alpha,\beta,\gamma,\delta)} \colon 
\left(
\begin{array}{c}
x \\
v \\
\end{array}
\right)
\mapsto 
\left(
\begin{array}{c}
\dfrac{\alpha x^2}{(v+\gamma)(v+\delta)} \\ [3mm]
v -  \dfrac{(v + \beta) x^2}{(v+\gamma)(v+\delta)}  \\
\end{array}
\right).
\end{equation}
The maps ${F}_0, {F}_1$ and ${F}_2$ belong to the above family and correspond to the case when $\gamma = \beta + \alpha$ and $\delta = \beta - \alpha$, where $\alpha, \beta$ are parameters depending on $y,z$ and $u$, according to \eqref{eq:random map}.

Therefore, ${F}_0, {F}_1$ and ${F}_2$ belong to the 2--parametric family $\{F_{(\alpha,\beta)} \}$, where $F_{(\alpha,\beta)} = {F}_{(\alpha,\beta,\beta + \alpha,\beta - \alpha)}$.
Similar to \eqref{eq:random map}, maps can be written for generalized overgroups $\wt\Gr\om$. They also fit in the 2--parametric family $\{ F_{(\alpha,\beta)} \}$.

Dynamical pictures of composition of the above maps for some sequences $\omega \in \Omega$ are shown in Figure~\ref{fig:random_maps}. The first Schur maps are much more complicated and so we have restricted this discussion to the second Schur map.


\section{Random Model and Concluding Remarks}

As explained above, the spectral problem associated with groups and their Schreier graphs in many important examples could be converted into study of invariant sets and dynamical properties of multi-dimensional rational maps. Some of these maps, like \eqref{eq:F map 1stGri}, \eqref{eq:G map 1stGri}, \eqref{map lamplighter}, \eqref{map Hanoi} demonstrate strong integrability features explored in \cite{BG00,BG00Hecke,GS06,GS07}. The roots of their integrability are comprehensively investigated in \cite{DGL20}. The examples given by \eqref{eq:5dim F map 1stGri}, \eqref{eq:5dim G map 1stGri}, \eqref{eq:map Basilica}, \eqref{eq:map img}, \eqref{eq:map img con}, \eqref{eq:2nd Schur overgp transformation} are much more complicated. They have \ronerthree{an} invariant set of fractal nature, and computer simulations demonstrate their chaotic behavior, shown by dynamical pictures given by Figures \ref{fig:basilica map} and \ref{fig:random_maps}.

The families of groups $\Gr\om, \wt\Gr\om,$ $\omega \in \Omega = \{0,1,2\}^\N$ (and many other similar families can be created) can be viewed as a random group if $\Omega$ is supplied with a shift invariant probability measure (for instance, Bernoulli or more generally Markov measure). The first step in this direction is publication \cite{BGV14} where it is shown that for any ergodic shift invariant probability measure satisfying a mild extra condition (all Bernoulli measures satisfy it), there is a constant $\beta < 1$ such that the growth function $\gamma_{\Gr\om}(n)$ is bounded by $e^{n^\beta}$.

More general model would be to supply the space $\M = \bigcup_{k=1}^\infty \M_k$ of finitely generated groups or any of its subspaces $\M_k$ \ronerthree{with} a measure $\mu$ (finite, or infinite, invariant or quasi-invariant with respect to any reasonable group or semigroup of transformations of the space) and study the typical properties of groups with respect to $\mu$. The system $(\Gr\om, \Omega, T, \mu)$ (where $T$ denotes the shift) is just \ronertwo{one} example of this sort.  As suggested in \cite{Gri05}, it would be wonderful if one \ronerthree{could supply the space $\M$ with} a measure that is invariant (or at least quasi-invariant) with respect to the group of finitary Nielsen transformations defined over infinite alphabet $\{x_1, x_2, \hdots \}$.

Additionally to the randomness of groups, \ronerthree{one can associate with each particular group} a random family of Schreier graphs, like the family $\Gamma_\xi$, $\xi \in \partial T$ for a group $G \leq \Aut(T)$ using the uniform Bernoulli measure on the boundary (other choices for $\mu$ are also possible, especially if $G$ is generated by automorphisms of polynomial activity \cite{Sid00,Kra10,DG17}). Putting all this together, it leads to study of random graphs associated with random groups (or equivalently, of random invariant subgroups in random groups).

Finally, even if we fix a group, say $\Gr\om$, $\omega \in \Omega$ and a Schreier graph $\Gamma_{\omega,\xi}$, $\xi \in \partial T$, study of spectral properties of this graph is related to the study of iterations $F_{\omega_{n-1}} \circ \hdots \circ F_{\omega_1} \circ F_{\omega_0}$ of maps given by \eqref{eq:random map} as was mentioned above.

Recall a classical construction of skew product in dynamical systems. Given two spaces $(X,\mu),(Y,\nu)$, the measure $\nu$ preserving transformation $S\colon Y \to Y$, and for any $y \in Y$ the measure $\mu$ preserving transformation $T_y \colon  X \to X$, under the assumption that the map $X \times Y \to X$, $(x,y) \mapsto T_y x$ is measurable, one can consider the map $Q \colon X \times Y \to X \times Y$, $Q(x,y) = (T_y x, Sy)$, which preserves the measure $\mu \times \nu$. Natural conditions imply that $Q$ is ergodic if $S$ and $T_y$, $y \in Y$ are ergodic. If for $k = 1,2,\hdots$, \ronertwo{we} put
\[ T_y^{(k)} = T_{S^ky} \circ T_{S^{k-1}y} \circ \hdots \circ T_{Sy} \circ T_y, \]
then the random ergodic theorem of Halmos--Kakutani \cite{Hal56,Kak51} states that for $f \in L^1(X,\mu)$, $\nu$ almost surely the averages $\frac 1n \sum_{k=0}^{n-1} f(T_y^{(k)} x)$ converge $\mu$ almost surely to some function $f^*_y (x) \in L^1(X,\mu)$.  In the simplest case when $Y = \{1,\hdots , m\}$ and $\nu$ is given by \ronerthree{a} probability vector $(p_1, \hdots ,p_m)$, $p_i>0$, $i=1,\hdots ,m$, $\sum_{i=1}^m p_i = 1$, we have $m$ transformations $T_1, \hdots , T_m$ of $X$. The semigroup generated by them typically is a free semigroup and if $T_1, \hdots, T_m$ are invertible in a typical case, the group $\gen{T_1, \hdots, T_m}$ is a free group of rank $m$.

The question of whether the pointwise ergodic theorem of Birkhoff holds for actions of a free group was raised by V.~Arnold and A.~Krylov \cite{AK62} and answered affirmatively by the first author in \cite{Gri87}. A similar theorem is proven for the action of a free semigroup \cite{Gri00Ergodic}. The proofs are based on the use of the skew product when $S$ is the Bernoulli shift on $\Lambda = \{1,2,\hdots,m\}^\N$ and $\nu$ is a uniform Bernoulli measure on $\Lambda$. In fact these ergodic theorems hold for stationary measures.

By a different method, a similar result for the free group actions was obtained by A.~Nevo and E.~Stein \cite{NS94}. The method from \cite{Gri87,Gri00Ergodic} was used by A.~Bufetov to get ergodic theorems for a large class of hyperbolic groups \cite{Buf02}. Ergodic theorem for action of non-commutative groups \ronerthree{became popular} \cite{BN15,BN19}, but the case of semigroup action is harder and not so many results are known, especially in the case of stationary measure.

In fact, the product measure $\mu\times\nu$ in the construction of the skew product is invariant if and only if the measure $\mu$ is $\nu$-stationary, which means the equality
\[ \mu = \int_Y (T_y)_* \mu \dd{\nu(y)}.\]
In the case of $Y = \{1,\hdots,m\}$, $\nu = (p_1, \hdots,p_m)$, it means that
\[ \mu = \sum_{i=1}^m p_i (T_i)_* \mu \]
(i.e., the $\nu$-average of images of $\mu$ under transformations $T_1, \hdots , T_m$ is equal to $\mu$). The skew product approach for non-commutative transformations leads to \ronerthree{a} not well-investigated notion of entropy \cite{Gri99Ergodic,Gri12}. It would be interesting to compare this approach with the approach of L.~Bowen for the definition of entropy of free group action \cite{Bow18,Bow20}.

Going back to \ronerthree{the} transformations $F_0, F_1, F_2$ given by \eqref{eq:random map}, we could try to apply the idea of skew product to them and investigate the random model. Random dynamical model in the context of holomorphic dynamics is successfully considered in \cite{CD20} where stationary measures also play an important role. Each of these maps (as well as any map from the 2-parametric family $F_{(\alpha,\beta)}$) is semiconjugate to the Chebyshev map and has families of \emph{horizontal} and \emph{vertical} hyperbolas similar to the case of the map $F$ given by Figure \ref{fig:hyperbola_all_var_G}. But when parameters $y,z,u$ (the coefficients of $b,c,d$) are not equal, even in the case of periodic sequence $\omega = (012)^\infty$ (in which case the group $\Gr\om$ is just our main \emph{hero} $\Gr$), the iterations $T_{\omega_{n-1}} \circ \hdots \circ T_{\omega_0}$ demonstrate chaotic dynamics presented by Figure \ref{subfi:dynamical map 012}. Still, it is possible that more chaos could appear if additionally, $\omega$ is chaotic itself. Study of these systems and other topics discussed above is challenging and promising.

The notion of amenable group was introduced by J.~von Neumann \cite{vNeu29} for discrete groups and by N.~Bogolyubov \cite{Bog39} for general topological groups. The concept of amenability entered many areas of mathematics \cite{Gre69,Wag93,HR63,Edw65,GH17,TW16,AG99}. Groups of intermediate growth and topological full groups remarkably extended the knowledge about the class $AG$ of amenable groups \cite{Gri84,Gri98,JM13}. There are many characterizations of amenability: via existence of left invariant mean (LIM), existence of F\"olner sets, Kesten's probability criterion, hyperfiniteness \cite{KM04}, co-growth \cite{Gri80Random}, etc.

From dynamical point of view, \ronerthree{an} important approach is due to Bogolyubov \cite{Bog39}; if a topological group $G$ with left invariant mean acts continuously on a compact set $X$, then there is a $G$-invariant probability measure $\mu$ on $X$ (this is a far \ronerthree{reaching} generalization of the famous Krylov-Bogolyubov theorem). In fact, such property characterizes amenability.

Amenability was mentioned in this article several times. We are going to conclude with open questions related to \ronerthree{the} considered maps $F,G$ and \ronerthree{the} conjugates $F_{(\alpha,\beta)}$ of $F$. \ronertwo{Even though} invariant measures \ronertwo{seem not to} play an important role in the study of dynamical properties of multi-dimensional rational maps (where the harmonic measure or measures of maximal entropy like the Man\'e-Lyubich measure dominate), we could be interested in existence of invariant or stationary measures supported on invariant subsets of maps coming from Schur complements, as \ronerthree{it} was discussed above. This could include the whole Schur semigroup $S_d^*$ defined in Section \ref{sec:Schur}, its subsemigroups, or semigroups involving some relatives of these maps (like the map $H$ given by \eqref{eq:H map}). The concrete questions are:
\begin{Problem}\label{prob:map}
	~\\
	\vspace{-\baselineskip}
	\begin{enumerate}[leftmargin=*,labelsep=4.9mm]
		\item Is the semigroup $\gen{F,H}$ amenable from the left or right?
		\item Is there a probability measure on the cross $\K$, shown by Figure \ref{subfig:all_hyperbola_interval_var_G}, invariant with respect to the above semigroup?
	\end{enumerate}
\end{Problem}

By the last part of Theorem \ref{thm:spec limit}, we know that each horizontal slice of the cross $\K$ \ronertwo{possesses} a probability measure that is \ronertwo{the} density of states for the corresponding Markov operator. Integrating it along the vertical direction we get a measure $\nu$ on $\K$ which is somehow related to both maps $F$ and $G$. Is it related to the semigroup $\gen{F,G}$? $\gen{F,G}$ is the simplest example of the Schur type semigroup. One can consider other semigroups of interest, for instance, $\gen{F_0,F_1,F_2}$ or even semigroup generated by the maps $F_{(\alpha,\beta)}$ for $\alpha,\beta \in \R$, and look for invariant or stationary measures.

The example of semigroup $\gen{F,G}$ is interesting because of the relations $H \circ F = G$, $H \circ G = F$. By J. Ritt's result \cite{Rit22}, it is known that in the case of the relation of the type $A \circ B = C \circ D$ for maps $A,B,C,D$ given by polynomials in one variable, all its solutions can be described explicitly. In the paper \cite{GTZ08}, it is proved that for polynomials $P$ and $Q$, if there exists a point $z_0$ in the complex Riemann sphere $\ol \CC$ such that the intersection of the forward orbits of $z_0$ with respect to $P$ and $Q$ is an infinite set, then there are natural numbers $n,m$ such that $P^n = Q^m$. The result of C. Cabrera and P.~Makienko \cite{CM20} generalizes this to the rational maps and includes in the statement the amenability properties of the semigroup $\gen{P,Q}$.

In the case of maps $F,G$, we know that the semigroup $\gen{F,G}$ is rationally semiconjugate to the commutative (and hence amenable) semigroup $\N \times \N$. Whether $\gen{F,G}$ is amenable itself is an open question included in the Problem \ref{prob:map}.

\begin{figure}[h]
	\centering
	\includegraphics[width=.5\textwidth]{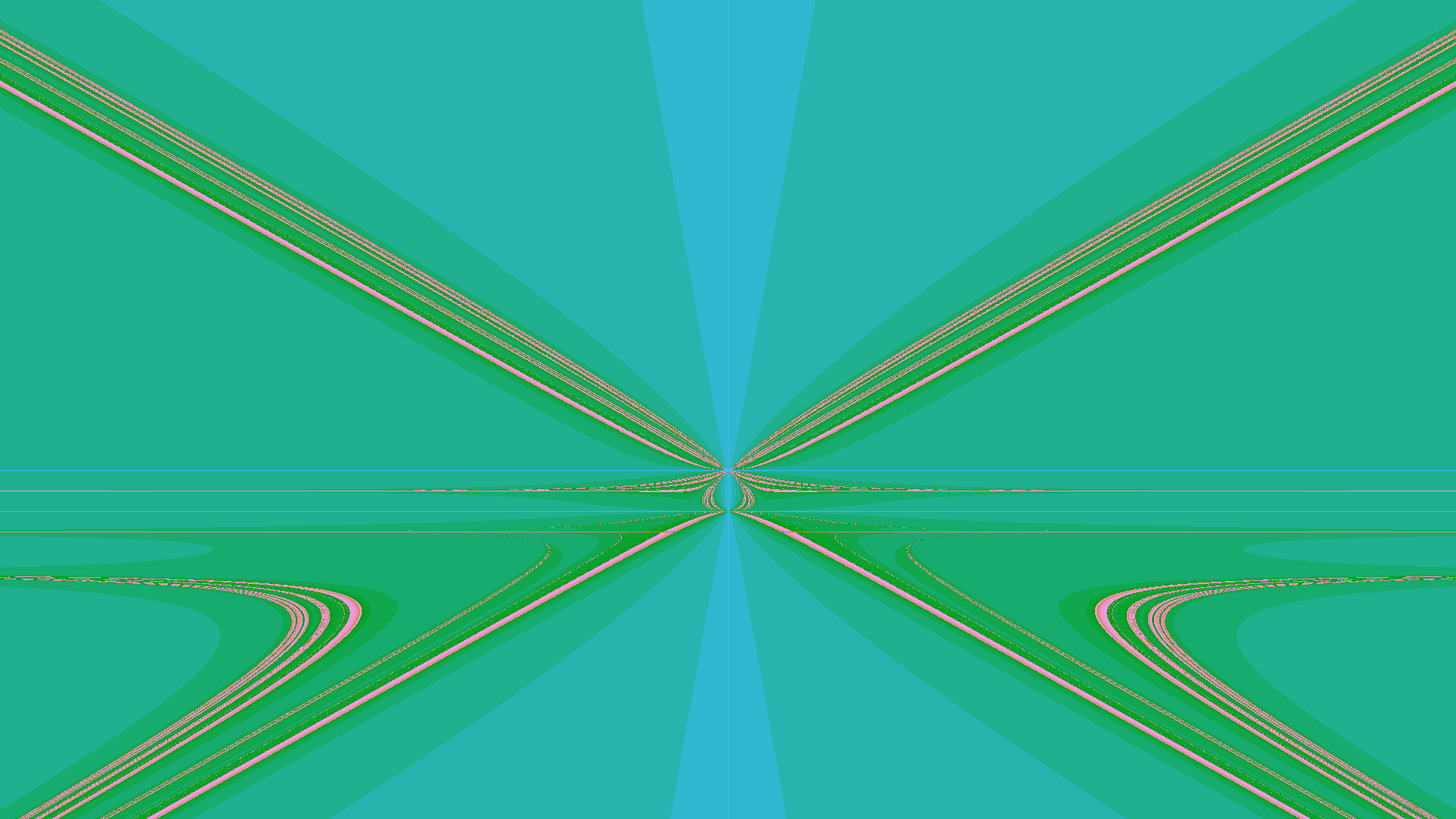}
	\vspace*{8pt}
	\caption{Dynamical pictures of ${F}_{(\alpha,\beta,\gamma,\delta)}$ for $(\alpha,\beta,\gamma,\delta) = (1,3,1.5,2.5)$.}\label{fig:4 para}
\end{figure}

A very interesting question is the question about the dynamical properties of maps $F_{(\alpha,\beta,\gamma,\delta)}$ given by \eqref{eq:4para map}. They are conjugate to the maps of the form
\begin{equation}\label{eq:3 para}
{F}_{(\alpha,\beta,\gamma)} \colon 
\left(
\begin{array}{c}
x \\
v \\
\end{array}
\right)
\mapsto 
\left(
\begin{array}{c}
\dfrac{\alpha x^2}{\gamma^2 -v^2} \\ [3mm]
v +  \dfrac{(v + \beta) x^2}{\gamma^2 -v^2}  \\
\end{array}
\right)
\end{equation}
via 
\[ 	S_{(\gamma,\delta)}  \colon 
\left(
\begin{array}{c}
x \\
v \\
\end{array}
\right)
\mapsto 
\left(
\begin{array}{c}
-x \\ [3mm]
- v -\dfrac {(\gamma+\delta)}2 \\
\end{array}
\right),
\]
and further simplification seems to be impossible. At the same time, maps $F_{(\alpha,\beta)}$ are conjugated to $F$ by 
\[ 	R_{(\alpha,\beta)} \colon 
\left(
\begin{array}{c}
x \\
v \\
\end{array}
\right)
\mapsto 
\left(
\begin{array}{c}
-\dfrac 2\alpha x \\ [3mm]
-\dfrac 2\alpha v -\dfrac 2\alpha \beta \\
\end{array}
\right).
\]
The dynamical picture for those that are outside the 2-parametric family $F_{(\alpha,\beta)}$ is presented by Figure~\ref{fig:4 para} and is quite different \ronerthree{of those dynamical pictures}  presented by Figure \ref{fig:random_maps}.

It is not clear at the moment if the maps presented in \eqref{eq:3 para} describe \ronerthree{the} joint spectrum of a pencil of operators associated with fractal groups. But the family \eqref{eq:3 para} itself could have interest for multi-dimensional dynamics and deserves to be carefully investigated, including semigroups generated by $F,H,R_{(\alpha,\beta)},S_{(\gamma,\delta)}$ in various combinations of the choice of generating set.


\acknowledgments{We thank Nguyen-Bac Dang and Mikhail Lyubich for stimulating discussions and valuable remarks. {Also, we are grateful to Tatiana Nagnibeda and Volodymyr Nekrashevych for their numerous suggestions and remarks that improved the exposition}. The first  author  is  grateful  to  Peter Kuchment  for  valuable  discussions  on Schur complement  and  density  of  states topics. The first author acknowledges a  partial  support from the Simons Foundation through the Collaboration Grant \#527814.
}


\reftitle{References}

\end{document}